\documentclass[11pt]{amsart}
\usepackage{amssymb, amscd}
\newcommand\al{\alpha}

\renewcommand\k{\kappa}

\newcommand\ld{\lambda}

\newcommand\x{\xi}

\newcommand\Ph{\Phi}

\newcommand\lieg{{\mathfrak g}}
\newcommand\liea{{\mathfrak a}}
\newcommand\liek{{\mathfrak k}}

\newcommand\liesl{{\mathfrak sl}}

\newcommand\liem{{\mathfrak m}}
\newcommand\lien{{\mathfrak n}}

\newcommand\CC{\mathbb C}

\newcommand\RR{\mathbb R}
\newcommand\ZZ{\mathbb Z}

\newcommand\HH{\mathbb H}

\newcommand\SL{{\mathrm{SL}}}

\newcommand\SU{{\mathrm{SU}}}
\newcommand\U{{\mathrm{U}}}

\newcommand\End{\operatorname{End}}
\newcommand\re{\operatorname{Re}}
\newcommand\im{\operatorname{Im}}

\makeatletter 
\newcommand\matc[4]{\left( {#1\@@atop #3}{#2\@@atop #4}\right)}
\newcommand\matr[4]{\left( {\hfill #1\@@atop\hfill #3}{\hfill
#2\@@atop\hfill #4}\right)}
\newcommand\matl[4]{\left( { #1\@@atop #3}{ #2\@@atop\hfill #4}\right)}
\makeatother
\newcommand\widearray[1]{\renewcommand\arraystretch{1.4}
\begin{array}{#1}}

\newcommand\lw[1]{\hbox{}_{#1}\!}

\newcommand\vz[1]{\mathchoice{\left\{ #1 \right\}}{\left\{ #1
\right\}}{\{#1
\}}{\{ #1 \}}}

\newcommand\vzm[2]{\mathchoice{\left\{\, #1 : #2 \,\right\}}{\{\, #1
:\allowbreak #2 \,\}}{\{ #1 :\allowbreak #2 \}}{\{ #1 :\allowbreak #2
\}}}

\theoremstyle{plain}
\newtheorem{thm}{Theorem}[section]
\newtheorem{lem}[thm]{Lemma}
\newtheorem*{propD}{Proposition \ref{Dexpresion}}
\newtheorem*{propE}{Proposition \ref{Eexpresion}}
\newtheorem*{propderxy}{Proposition \ref{dderxy}}

\newtheorem{prop}[thm]{Proposition}
\newtheorem{cor}[thm]{Corollary}
\newtheorem{conj}[thm]{Conjecture}
\theoremstyle{definition}
\newtheorem{defn}[thm]{Definition}

\theoremstyle{remark}
\newtheorem{remark}[thm]{Remark}

\newcommand\rmk[1]{\medskip\par\noindent{\em #1. }\ignorespaces}

\title[Spherical functions]{Matrix valued spherical
functions \\associated to\\
the three dimensional hyperbolic space}

\author{ F. A. Gr\"unbaum}
\address{Departament of Mathematics, University of California, Ber\-ke\-ley
CA 94705}
\email{grunbaum@math.berkeley.edu}

\author{I. Pacharoni}
\address{CIEM-FaMAF, Universidad Nacional de C\'or\-do\-ba,
C\'or\-do\-ba~5000, Argentina}
\email{pacharon@mate.uncor.edu}

\author{J. Tirao}
\address{CIEM-FaMAF, Universidad Nacional de C\'or\-do\-ba,
C\'or\-do\-ba~5000, Argentina}
\email{tirao@mate.uncor.edu}

\thanks{This paper is partially supported by NSF grants FD9971151 and
1-443964-21160 and by CONICET grant PIP655-98.}

\begin{document}
\begin{abstract}
The main purpose of this paper is to compute all
irreducible spherical functions
on $G=\SL(2,\CC)$ of arbitrary type $\delta\in \hat K$, where
$K=\SU(2)$. This is accomplished by
associating to a spherical
function $\Phi$ on $G$ a matrix valued function $H$ on the three dimensional
hyperbolic space $\HH=G/K$. The entries of $H$ are solutions of two coupled
systems of ordinary differential equations. By an appropriate twisting
involving Hahn polynomials we uncouple one of the systems and express
the entries of $H$ in terms of Gauss' functions $\lw2F_1$. Just as in the
compact instance treated in \cite{GPT} there is a useful role for
a special class of generalized hypergeometric functions $\lw{p+1}F_p$.
\end{abstract}
\maketitle

\section{Introduction and statement of results}\label{statements}

Let $G$ be a locally compact unimodular group and let
$K$ be a compact subgroup of $G$. Let $\hat K$ denote
the set of all
equivalence classes of complex finite dimensional
irreducible representations of $K$; for each
$\delta\in \hat K$, let
$\x_\delta$ denote the character of $\delta$,
$d(\delta)$ the degree of $\delta$, i.e. the dimension
of any representation in
the class $\delta$, and
$\chi_\delta=d(\delta)\x_\delta$. We shall choose once
and for all the Haar measure $dk$ on
$K$ normalized by $\int_K dk=1$.

We shall denote by $V$ a finite dimensional vector
space over the field $\CC$ of complex numbers and by
$\End(V)$ the space
of all linear transformations of $V$ into $V$.
Whenever we refer to a topology on such a vector space
we shall be talking
about the unique Hausdorff linear topology on it.

A spherical function $\Ph$ on $G$ of type $\delta\in
\hat K$ is a continuous function on $G$ with values in
$\End(V)$ such
that
\begin{enumerate} \item[i)] $\Ph(e)=I$. ($I$= identity
transformation).

\item[ii)] $\Ph(x)\Ph(y)=\int_K
\chi_{\delta}(k^{-1})\Ph(xky)\, dk$, for all $x,y\in
G$.
\end{enumerate}

The reader can find a number of general results in
\cite{T} and \cite{GV}. For our purpose it is
appropriate to recall

\begin{prop}\label{propesf}(\cite{T},\cite{GV})  If
$\Ph:G\longrightarrow \End(V)$ is a spherical function
of type $\delta$
then:
\begin{enumerate}
\item[i)] $\Ph(kgk')=\Ph(k)\Ph(g)\Ph(k')$, for all
$k,k'\in G$, $g\in G$.
\item[ii)] $k\mapsto \Ph(k)$ is a representation of
$K$ such that any irreducible subrepresentation
belongs to $\delta$.
\end{enumerate}
\end{prop}

Concerning the definition let us point out that the
spherical function $\Ph$ determines its type
univocally
(Proposition \ref{propesf}) and let us say that the
number of times that $\delta$ occurs in the
representation
$k\mapsto \Ph(k)$ is called the {\em height} of $\Ph$.

The three dimensional hyperbolic space $\HH=\CC\times
\RR^+$ can be realized as the homogeneous space $G/K$,
where $G=\SL(2,\CC)$
and $K={\mathrm{S}}\U(2)$. We are interested in
determining, up to equivalence, all irreducible
spherical functions,
associated to the pair $(G,K)$. If $(V,\pi)$ is a
finite dimensional irreducible representation of $K$
in the equivalence
class $\delta \in \hat K$, a spherical function on $G$
of type $\delta$ is characterized, see
\cite{T},\cite{GV}, by
\begin{enumerate}
\item[i)] $\Ph:G\longrightarrow \End(V)$ is analytic.
\item[ii)] $\Ph(k_1gk_2)=\pi(k_1)\Ph(g)\pi(k_2)$, for
all $k_1,k_2\in K$, $g\in G$, and $\Phi(e)=I$.
\item[iii)] $[\Omega\Ph ](g)=\tilde\ld\Ph(g)$,
$[\overline\Omega\Ph ](g)=\tilde\mu\Ph(g)$ for all
$g\in G$ and for
some $\tilde\ld,\tilde\mu\in \CC$.
\end{enumerate}
Here $\Omega$ and $\overline\Omega$ are two
algebraically independent generators of the
polynomial algebra $D(G)^G$ of all differential
operators on $G$ which are invariant under left and
right multiplication by
elements in $G$. A particular choice of these
operators is given in Proposition \ref{generadoresD(G)G}.

The set $\hat K$ can be identified with the set
$\ZZ_{\geq 0}$. If $A\in{\mathrm{S}}\U(2)$ then
$\pi(A)=\pi_\ell(A)=A^\ell$, where $A^\ell$ denotes
the $\ell$-symmetric power of $A$, defines an
irreducible representation of $K$
in the class $\ell\in\ZZ_{\geq 0}$.

The representation $\pi_\ell$ of $K={\mathrm{S}}\U(2)$
extends to a unique holomorphic representation of
$G=\SL(2,\CC)$, which we shall still denote by
$\pi_\ell$. Let $\Phi_\pi:G\longrightarrow
\End(V_\pi)$ be defined by
$$\Phi_\pi(g)=\pi_\ell(g).$$
Then applying the definition it follows that
$\Phi_\pi$ is a spherical function of type $\ell$.
These spherical functions
will play a crucial role in the rest of the paper.

We are now in a position to describe the plan of the paper.
Section \ref{prelim} contains a brief
review of standard facts.

To determine all spherical functions
$\Phi:G\longrightarrow \End(V_\pi)$ of type
$\pi=\pi_\ell$, we use the function $\Phi_\pi$ to
define a function
$H:G\longrightarrow \End(V_\pi)$ by
$$H(g)=\Phi(g)\, \Phi_\pi(g)^{-1}.$$
Then $H$ satisfies
\begin{enumerate}
\item [i)] $H(e)=I$.
\item [ii)]$ H(gk)=H(g)$, for all $g\in G, k\in K$.
\item [iii)] $H(kg)=\pi(k)H(g)\pi(k^{-1})$, for all
$g\in G, k\in K$.
\end{enumerate}

Property ii) says that $H$ may be considered as a
function on the hyperbolic space $\HH$.

The fact that $\Phi$ is an eigenfunction of $\Omega$
and $\overline\Omega$, makes $H$ into an eigenfunction
of certain
differential operators $D$ and $E$ on $\HH$. This is done in
Section \ref{redG/K} . For
completeness the explicit computation of these
operators is
carried out fully in the Appendix.

In Section \ref{onevariable} we take full advantage of the
$K$-orbit structure of $\HH$ combined with property
iii) of our
functions $H$. The $K$-orbits in $\HH$ are the spheres
with center in the positive axes $\vzm{(0,r)}{r>0}$
with north and south poles of the form $(0,s)$ and
$(0,\frac 1s)$, and the single point set $\vz{(0,1)}$.
Thus the set of $K$-orbits in $\HH$ is parametrized by
the interval $(0,1]$ or by $[1,\infty)$.

It follows that there exist ordinary differential operators
$\tilde D$ and $\tilde E$ acting on the space
$C^\infty((0,1))\otimes \End(V_\pi)$
such that $$(D\,H)(0,r)=(\tilde D\tilde H)(r)\, ,
\quad (E H)(0,r)=(\tilde E\tilde H)(r),$$
where $\tilde H(r)=H(0,r)$, $r\in(0,1)$.
These operators $\tilde D$ and $\tilde E$ are
explicitly given in Theorems \ref{mainth1} and
\ref{mainth2}. We need to
compute a number of second order partial derivatives
of the function $H:\HH\longrightarrow \End(V_\pi)$ at
the point
$(0,r)$. This detailed computation is broken down in a
number of lemmas included in the Appendix for the
benefit of
the reader.

Theorems \ref{mainth1} and \ref{mainth2} are given in
terms of linear transformations. The functions $\tilde
H$ turn out to
be diagonalizable (Lemma \ref{Hdiagonal}). Thus, in an
appropriate basis of $V_\pi$ we can write $\tilde
H(r)=(\tilde h_0(r),\cdots,
\tilde h_\ell(r))$. Then we give in Corollaries
\ref{sistema} and \ref{sistema2} the corresponding
statements of these theorems in
terms of the scalar functions $\tilde h_i$.

We also introduce here the variable $t=r^2$ which
converts the differential operators $\tilde D$ and
$\tilde E$
into new operators $D$ and $E$, and the functions
$\tilde H, \tilde h_i$ into the functions $H, h_i$
defined by
$H(t)=\tilde H(\sqrt t)$ and $h_i(t)=\tilde h_i(\sqrt
t)$. At this point there is a slight abuse of
notation.
The resulting equation $DH=\lambda H$ is a coupled
system of
$\ell+1$ second order differential equations in the
components
$(h_0,\cdots,h_\ell)$ of $H$. Fortunately the
coupling matrix $C_0+C_1$ of the system $DH=\lambda
H$ is symmetric, with
eigenvalues $-j(j+1)$, $0\le j\le\ell$, and
eigenvectors $u_j$ whose components $u_{i,j}$ are
given by certain Hahn orthogonal
polynomials,  see Proposition \ref{Hahn}. Let $U$ denote
the $(\ell+1)\times(\ell+1)$ matrix $(u_{i,j})$ which
will play a crucial
role. In fact, this $U$ allows us to decouple the
system above. More is true, the twisting
$\check{H}=U^{-1}H$ leads us to
Gauss' hypergeometric equation, a fact that we exploit
fully.

By well known reasons the smooth functions
$H:(0,1)\longrightarrow\CC^{\ell+1}$ that satisfy
$DH=\lambda H$ are analytic
functions and the dimension of the corresponding
eigenspace $V_\lambda$ is $2(\ell+1)$.

To get a handle on $V_\lambda$ we consider functions
$H\in V_\lambda$ of the form $H(t)=t^pF(t)$,
$p\in\CC$, with $F$
analytic at $t=0$ and $F(0)\ne 0$. This forces
$\lambda=4p(p-1)$. Observe that $p$ and $1-p$ give
rise to the same $\lambda$.
Indeed any statement valid for $\re(2p)\ge1$ has a
mirror image for $\re(2p)\le1$ obtained by changing $p$
into $1-p$.

Let $V(p)$ be the linear subspace of $V_\lambda$,
$\lambda=4p(p-1)$, of all $H\in V_\lambda$ of the form
$H(t)=t^pF(t)$ with $F$ analytic at $t=0$.

When $2p$ is not an integer we prove that
$V_\lambda=V(p)\oplus V(1-p)$. When $2p$ is an integer
the situation is more
complicated: if $2p\ge 1$ we have $V(p)\subset V(1-p)$
and $\dim V(p)=\ell+1$, $\dim V(1-p)=\min\{2(\ell+1),\ell+2p\}$.
The situation when $2p\le 0$ can be read off by
exchanging $p$ by $1-p$, see Proposition \ref{dimV(p)}.

When $2p\notin\ZZ$ we consider the linear map
$\eta:V(p)\longrightarrow\CC^{\ell+1}$ defined by
$\eta(H)=\lim_{t\rightarrow 0^+}(t^{-p}H(t))$. It
turns out that $\eta$ is an isomorphism and by using
it $E$ is described by an
$(\ell+1)\times(\ell+1)$ matrix $L(p)$, see
Proposition \ref{etadef1}. The eigenvalues of $L(p)$ are the
same as those of $L(1-p)$,
all of them with algebraic multiplicity one. Using
Proposition \ref{Eautovalores} combining the eigenvectors of $L(p)$
and $L(1-p)$
with the same eigenvalue $\mu$ we obtain all
simultaneous solutions of the system $DH=\lambda H$
and $EH=\mu H$. To get
our hands on the spherical functions we need to focus
on those $H$'s arising from {\it the} linear
combination of eigenvectors
of $L(p)$ and $L(1-p)$ with the same eigenvalue which
insures that $\lim_{t\rightarrow 1^-}H(t)=(1,\dots,1)$.

This requirement leads us to consider the linear subspace
$W_\lambda$ of $V_\lambda$ consisting of all $H$ such that
$\lim_{t\rightarrow 1^-}H(t)$
is finite. For any $p\in\CC$ with $\dim V(p)=\ell+1$,
$W_\lambda$ turns out to be a useful direct complement of
$V(p)$ in $V_\lambda$, see Propositions \ref{complementoV(p)} and
\ref{complementoV(p)2}. Furthermore $E$
preserves $W_\lambda$ and the linear isomorphism
$\eta:W_{\lambda}\longrightarrow\CC^{\ell+1}$ defined
by

\begin{equation}\label{etaP}
\eta(H)=
\begin{cases}
\lim_{t\rightarrow 0^+}t^{p-1}H(t)&\text{if
Re$(2p)>1$},\\
\lim_{t\rightarrow 0^+}t^{-p}H(t)&\text{if
Re$(2p)<1$},\\
\lim_{t\rightarrow 0^+}t^{p-1}P(H)(t) &\text{if Re$(2p)=1$
and $2p\ne 1$},\\
\lim_{t\rightarrow 0^+}(t^{\frac12}\log
t)^{-1}H(t)&\text{if $2p=1$},
\end{cases}
\end{equation}
reduces the determination of all spherical functions
to a linear algebra problem. More precisely, since we
establish in  Propositions \ref{etadiagr}, \ref{etadiagr1} and
\ref{etadiagr2},
that, in $W_\lambda$,
$E=\eta^{-1}L(1-p)\eta$ it suffices to find the
eigenvectors of
$L(1-p)$. See Proposition \ref{etadiagr2} for the definition
of the projection $P$ appearing in \eqref{etaP}.

Proposition \ref{Eautovalores} determines the eigen-structure
of $L(p)$ and Theorem \ref{punch} gives the main result of the paper,
namely an explicit expression for all spherical functions.
The section closes with a look at the relation between our bare hands
construction and the generalized (also known as nonunitary) principal series
representations of SL$(2,\CC)$.

Section \ref{ejemplos} displays explicit expressions of all
spherical functions for low values of $\ell$. The reader might
use these examples to get a better feeling for the results
in the paper. For instance we use these explicit expressions
to see that as a function of $p$, $H(t,p)$ is a meromorphic function
with a finite number of {\em removable} singularities. The values
at these singularities are exactly the cases when the
expression for $H(t,p)$
involves the term $\log t$. We also formulate a conjecture
expressing our spherical functions in terms of {\it generalized}
hypergeometric functions.

In Section \ref{unit} we discuss some generalities derived
from the definition of unitary spherical functions and we pick
among all spherical functions of our pair $(G,K)$ those that are
unitarizable, see Corollary \ref{unitarizables}.

In Section \ref{biespectral} we observe that a matrix valued
function $\Phi(t,p)$ put together from the functions
$H_0(t,p),\dots,H_\ell(t,p)$ alluded to in the comments preceding
Theorem \ref{punch} satisfies not only
differential  equations in $t$ but also a difference equation in $p$.

Finally in the Appendix we collect a number of explicit
computations for the benefit of the reader.

There is a point of contact between the present paper and
\cite{W}. By considering spherical functions associated to
generalized principal series representations A. Wang
derives two coupled systems of differential equations closely
related to the systems in Corollary \ref{DandEint}. From these
systems the author obtains a single second order differential
equation for the last component, and points out that,
in principle, this gives a way of finding the remaining
components.

Finally we remark that we have chosen this example as one of the simplest
to analyze among the noncompact symmetric spaces of rank one
leading to matrix valued spherical functions. We hope to deal with other
simple examples in the near future. The consideration of some more general
examples appears, at this point, to be a very interesting challenge.
\section{Preliminaries}\label{prelim}

We first review well known facts about the structure of the real Lie algebra
$\liesl(2,\CC)$ and of the center of its universal enveloping algebra.
We also collect here some useful facts about the quotient
$\SL(2,\CC)/\mathrm{S}\U(2)$.
A basis of $\lieg=\liesl(2,\CC)$ over $\RR$ is given by

$$H_1=\left[ \begin{matrix}
1&0 \\ 0&-1\end{matrix}\right], \quad
W_1=\left[\begin{matrix}
0&1 \\ -1&0\end{matrix}\right],\quad
V_1=\left[
\begin{matrix}
0&1 \\ 1&0\end{matrix}\right],$$

$$H_2=\left[\begin{matrix}
i&0 \\ 0&-i\end{matrix}\right], \quad
W_2=\left[\begin{matrix}
0&i \\ i&0\end{matrix}\right], \quad
V_2=\left[ \begin{matrix}
0&i \\ -i&0\end{matrix}\right]\;.$$

\medskip

We note that  $H_1, W_1$ and $V_1$ generate a Lie subalgebra isomorphic to
$\liesl(2,\RR)$, and that $H_2, W_1$ and $W_2$
form a basis of the Lie algebra of $K$.

\smallskip
\noindent The following are elements in the complexification  $\lieg_\CC$
of the Lie algebra $\lieg$.
$$\textstyle T=\frac 12 (H_1-iH_2)\, , \quad V=\frac 12 (V_1-iW_2)\, ,
\quad
W=\frac 12 (W_1-iV_2).$$

$$\textstyle \overline T=\frac 12 (H_1+iH_2)\, , \quad \overline V=\frac 12
(V_1+iW_2)\, , \quad \overline W=\frac 12 (W_1+iV_2).$$

\begin{prop} \label{generadoresD(G)G}
$D(G)^G$ as a polynomial algebra is generated by the
algebraically independent elements
$$\Omega= T^2+V^2 -W^2\; \quad \text{ and } \quad
\overline\Omega= {\overline T}^2+{\overline V}^2 -{\overline W}^2.$$
In the real basis of $\liesl(2,\CC)$ we have
$$4\Omega=H_1^2+V_1^2+V_2^2-(H_2^2+W_1^2+W_2^2)-2i(H_1H_2+V_1W_2-V_2W_1),$$
$$4\overline\Omega=H_1^2+V_1^2+V_2^2-(H_2^2+W_1^2+W_2^2)+2i(H_1H_2+V_1W_2-V_
2W_1).$$
\end{prop}

\medskip

We have a function from $G=\SL(2,\CC)$ into the space of all $2\times
2$
symmetric  positive definite matrices, given by
$g\mapsto gg^*$. This  function factors through  the quotient $G/K$,
because
$kk^*=I$.
Explicitly we have,
$$\left(\begin{matrix}a&b\\c&d\end{matrix} \right) \longmapsto
\left(\begin{matrix}
|a|^2+|b|^2 & a\bar c +b\bar d\\
\bar a c +\bar b d& |c|^2+|d|^2\end{matrix} \right).$$

\noindent If we put
\begin{equation}\label{variables}
r= |c|^2+|d|^2 \quad \text{and}\quad z=a\bar c +b\bar d,
\end{equation}
we have  $|a|^2+|b|^2= \frac {1+|z|^2}r$ (since $gg^*$ has determinant
one)
and
$$gg^*=  \left(\begin{matrix}\frac {1+|z|^2}r&z\\\bar z&r\end{matrix}
\right). $$
Then we can identify the quotient $G/K$ with
$$\HH=\vzm{(z,r)}{z\in \CC, r\in \RR_{>0}}.$$
The projection map $p:G\longrightarrow \HH$ is given by
$p(g)=(gg^*) \left (\begin{smallmatrix}
0\\1\end{smallmatrix}\right)=(z,r)$.

\noindent To simplify notation we identify row and column vectors.

\

\noindent The left multiplication on the group $G$ induces the action of
$G$
in $\HH$ given by
\begin{equation}\label{accion}
A \, (z,r)=A \matc{\frac {1+|z|^2}r}{z}{\bar z}{r} A^*\left
(\begin{smallmatrix}
0\\1\end{smallmatrix}\right)=(z^*, r^*).
\end{equation}
  More explicitly if $g=\matc{a}{b}{c}{d}$ and $p=(z,r)$ we have $g\cdot
p=(z^*,r^*)$
with
$$z^*=\frac{a\bar c+(az+br)(\bar c\bar z +\bar d r)}r  \quad \text{and}
\quad
r^*=\frac{|c|^2+|cz+dr|^2}r.$$

\section{Reduction to $G/K$}\label{redG/K}

Any irreducible finite dimensional representation of $K=\SU(2)$ is of
the form
$\pi_\ell(k)=k^\ell$, where $k^\ell$ denotes the $\ell$-symmetric power
of $k$. The representation $\pi_\ell$ extends to a unique holomorphic
representation of $\SL(2,\CC)$ and we also denote it by $\pi_\ell$.
In either case, we denote by $\dot \pi=\dot\pi_\ell$ the corresponding
derivative of $\pi$ at the identity. Note that $\dot\pi:\lieg
\longrightarrow \End(V_\pi)$ is $\CC$-linear.

We define for each representation $\pi=\pi_\ell$ of $K$ a
function $\Phi_\pi: G\longrightarrow \End(V_\pi)$ given by
$$\Phi_\pi(g)=\pi(g).$$

In order to determine all spherical functions $\Phi:G\longrightarrow
\End(V)$ we find it useful to introduce the function $H$ given by
$$H(g)=\Phi(g)\Phi_\pi(g)^{-1}.$$
This function $H$ satisfies

\begin{enumerate}
\item [i)] $H(e)=I$.
\item [ii)] $H(gk)=H(g)$, for $g\in G$, $k\in K$.
\item [iii)] $H(kg)=\pi(k)H(g)\pi(k^{-1})$, for $g\in G$, $k\in K$.
\end{enumerate}

\smallskip

Property ii) says that $H$ can be considered as a function on the
quotient
$G/K$. The fact that $\Phi$ is an
eigenfunction of $\Omega $ and $\overline \Omega$ makes the function
$H$,
introduced before, into an eigenfunction of certain differential
operators
on $\HH$, to be determined now.

Let
\begin{align}
\label{Ddefuniv} DH&=(H_1^2+V_1^2 +V_2^2)(H)\\
\begin{split}\label{Edefuniv}
EH&= H_1(H)\, \Phi_\pi \dot\pi(H_1)\Phi_\pi^{-1}+V_1(H)\, \Phi_\pi
\dot\pi(V_1)\Phi_\pi^{-1}\\ &\quad +V_2(H)\, \Phi_\pi
\dot\pi(V_2)\Phi_\pi^{-1}.
\end{split}
\end{align}

\

\begin{prop}\label{relacionautovalores}
  For any $H\in C^{\infty}(G)\otimes
\End(V_\pi)$ right invariant under $K$, the function $\Phi=H\Phi_\pi$
satisfies  $(4\overline \Omega) \Phi=\tilde\lambda \Phi$ and
$(\Omega-\overline\Omega)\Phi=\tilde\mu \Phi$
  if and only if $H$ satisfies $DH=\lambda H$ and $E H=\mu H$, with
$$\tilde\lambda =\lambda \quad \text{and}
\quad \tilde \mu =\mu +\ell(\ell+2),$$
when $\pi=\pi_\ell$.
\end{prop}

\begin{proof} If $X\in \liek$ then $X(H)=0$. In particular
\begin{equation*}
H_2(H)=W_1(H)=W_2(H)=0.
\end{equation*}

Therefore  $T(H)=\overline T(H)=\textstyle \frac 12 H_1(H)$, $V(H)=\overline
V(H)=\textstyle \frac 12 V_1(H)$,  $W(H)=-\overline W(H)=-\textstyle \frac
i2
V_2(H),$ and
$$ 4\Omega(H)=4\overline \Omega(H)= (H_1^2+V_1^2+V_2^2)(H)= DH.$$

We also have
\begin{align*}
\Omega(H\Phi_\pi)&= \Omega(H) \Phi_\pi+ H\Omega(\Phi_\pi)+
2T(H)T(\Phi_\pi)+2V(H)V(\Phi_\pi) \\&\quad -2W(H)W(\Phi_\pi)
\\ &=\Omega(H) \Phi_\pi+H\Phi_\pi \dot\pi(\Omega)+
H_1(H)\Phi_\pi\dot\pi(T)+V_1(H)\Phi_\pi\dot\pi(V)\\
&\quad +i V_2(H)\Phi_\pi\dot\pi(W)
\end{align*}
and
\begin{align*}
\overline \Omega(H\Phi_\pi)&= \overline \Omega(H) \Phi_\pi+H\Phi_\pi
\dot\pi(\overline \Omega)
+ H_1(H)\Phi_\pi\dot\pi(\overline T)+V_1(H)\Phi_\pi\dot\pi(\overline V)\\
&\quad -i V_2(H)\Phi_\pi\dot\pi(\overline W).
\end{align*}

Note that if $Y=iX$, for $X,Y\in \lieg$ then $\dot\pi(Y)=i\dot\pi(X)$. In
particular we have  $\dot\pi(\overline T)=\dot\pi(\overline
V)=\dot\pi(\overline W)=0$.
Therefore
\begin{align*}
4\overline\Omega (H\Phi_\pi)&= (DH)\Phi_\pi + 4H\Phi_\pi \dot\pi(\overline
\Omega)\\
(\Omega-\overline\Omega)(H\Phi_\pi)&= H\Phi_\pi \dot\pi(
\Omega-\overline \Omega) +EH \Phi_\pi
\end{align*}

Finally the representation theory of SL$(2,\CC)$ gives, for $\pi=\pi_\ell$,
$\dot\pi(\overline\Omega)=0$ and $\dot\pi(\Omega)=\ell(\ell+2)I$.
Now the proposition follows easily.
\end{proof}

\

Given $H\in C^\infty(\HH)\otimes\End(V_\pi)$ we shall also denote by $H\in
C^\infty(G)\otimes\End(V_\pi)$ the function defined by
$H(g)=H(p(g))$, $g\in G$. Moreover, if $F$ is a linear endomorphism of
$C^{\infty}(G)\otimes\End(V_\pi)$ which preserves the subspace
$C^{\infty}(G)^K\otimes\End(V_\pi)$ of all functions which are
right invariant by elements in $K$, then we shall also denote by $F$ the
endomorphism of $C^\infty(\HH)\otimes\End(V_\pi)$ which satisfies
$F(H)(p(g))=F(H)(g)$, $g\in G$, $H\in
C^\infty(\HH)\otimes\End(V_\pi)$.

\smallskip
\begin{lem} The differential operators $D$ and $E$  introduced
in \eqref{Ddefuniv} and \eqref{Edefuniv}, define differential operators $D$
and $E$ acting on $C^\infty(\HH)\otimes\End(V_\pi)$.
\end{lem}
\begin{proof}
The only thing we really need to prove is that $D$ and $E$
preserve the subspace $C^{\infty}(G)^K\otimes\End(V_\pi)$.

It is easy to see that $D=H_1^2+V_1^2 +V_2^2=-(H_2^2+W_1^2 +W_2^2)$ is a
multiple of the Casimir operator of $K$. Then $D$ preserves
$C^{\infty}(G)^K\otimes\End(V_\pi)$.

Let us now check that $E$ has the same property. Since $K$ is connected
this is equivalent to verifying that for any $X\in\liek$ and all
$H\in C^{\infty}(G)^K\otimes\End(V_\pi)$ we have $X(EH)=0$,
which in turns amounts to prove that
\begin{align*}
&X H_1(H)\Ph_\pi\dot\pi(H_1)+  H_1(H) X(\Phi_\pi) \dot\pi(H_1)-
H_1(H) \Phi_\pi \dot\pi(H_1)\dot\pi(X) \displaybreak[0]\\
&+ X V_1(H)\Ph_\pi\dot\pi(V_1)+  V_1(H) X(\Phi_\pi) \dot\pi(V_1)-
V_1(H) \Phi_\pi \dot\pi(V_1)\dot\pi(X) \displaybreak[0]\\
&+ X V_2(H)\Ph_\pi\dot\pi(V_2)+  V_2(H) X(\Phi_\pi) \dot\pi(V_2)-
V_2(H) \Phi_\pi \dot\pi(V_2)\dot\pi(X)=0
\end{align*}

\noindent Because $X(H)=0$ and $X(\Phi_\pi)=\Phi_\pi\dot\pi(X)$, this is
also equivalent to showing that
\begin{equation}\label{reduccion}
\begin{split}
&[X,H_1](H)\Ph_\pi \dot\pi( H_1)+ H_1(H) \Ph_\pi \dot\pi([X,H_1]) +
[X,V_1](H)\Phi_\pi \dot\pi(V_1)\displaybreak[0]\\
&+V_1(H) \Ph_\pi \dot\pi( [X,V_1])+ [X,V_2](H)\Phi_\pi \dot\pi(V_2)+ V_2(H)
\Ph_\pi \dot\pi( [X,V_2])
=0.
\end{split}
\end{equation}

\noindent If we put $X=H_1$ in \eqref{reduccion} and use
$[H_1,V_1]=2W_1$, $[H_1,V_2]=2W_2$, $W_2=iV_1$, $V_2=iW_1$   we get
\begin{align*}
2W_1(H)\Phi_\pi \dot\pi(V_1)+&2V_1(H)\Phi_\pi \dot\pi(W_1)+ 2W_2(H)\Phi_\pi
\dot\pi(V_2)\\ &+ 2V_2(H)\Phi_\pi \dot\pi(W_2)=0.
\end{align*}

\noindent If we substitute $X=\frac 12 (V_1+W_1)$ in \eqref{reduccion} and
use
$[H_1,X]=2X$,
$[X_,V_1]=H_1$, $[X_,V_2]=-iH_1$, we obtain
\begin{align*}
-V_1(H)& \Phi_\pi\dot\pi(H_1)-
W_1(H)\Phi_\pi\dot\pi(H_1)-H_1(H)\Phi_\pi\dot\pi(V_1)  -
H_1(H)\Phi_\pi\dot\pi(W_1)\\
& +H_1(H)\Phi_\pi\dot\pi(V_1) + V_1(H)\Phi_\pi\dot\pi(H_1)
- iH_1(H)\Phi_\pi\dot\pi(V_2) \\ &- i V_2(H)\Phi_\pi\dot\pi(H_1)=0,
\end{align*}
since $V_2=iW_1$, $W_2=iV_1$.
Now the representation theory of $\SL(2,\CC)$ tell us that $X(EH)=0$ for
all $X\in \liek$. This finishes the proof of the lemma.
\end{proof}

\smallskip
Now we give the expressions of the operators $D$ and $E$ on $\HH$ in the
real linear coordinates $(x,y,r)$ defined by $z=x+iy$. The detailed proofs
appear in the Appendix at the end of the paper.

\begin{prop}\label{Dexpresion} For any $H\in C^\infty(\HH)\otimes
\End(V_\pi)$
we have
\begin{align*}
\textstyle\frac 14 DH&= (\re^2 z +1)\,H_{xx}+(\im^2 z
+1)\,H_{yy}+r^2\,H_{rr}+2\re z\im z \, H_{xy} \\
& \quad
+2r\re z \,H_{xr} +2r\im z \, H_{yr}+3\re z\, H_x+3\im z \,H_y+3r
\,H_r.
\end{align*}
\end{prop}

\

\begin{prop}\label{Eexpresion} For any $H\in C^\infty(\HH)\otimes
\End(V_\pi)$
we have
\begin{align*}
EH&= H_x \;\dot\pi\matc{z-\overline
z}{\;\;\frac{2(1+|z|^2)}r}{2\,r}{\overline z-z} +
H_y \;\dot\pi\matc{i(z+\overline
z)}{\;\;\frac{2i(1+|z|^2)}r}{-2i\,r}{i(z+\overline z)}\\
&\quad +H_r \;\dot\pi\matc{-2r}{\;\;4z}{0}{\;\;2r}.
\end{align*}
\end{prop}

\smallskip

Summarizing the results of this section: to determine a simultaneous
$C^\infty$-eigenfunction $\Phi$ of $\Omega$ and $\overline\Omega$ on $G$
satisfying
$\Phi(k_1gk_2)=\pi(k_1)\Phi(g)\pi(k_2)$ for all $k_1,k_2\in K$, $g\in G$,
is equivalent to finding a $C^\infty$-eigenfunction $H=\Phi\Phi_\pi^{-1}$
of $D$ and $E$ on $\HH$ satisfying
$H(k\cdot p)=\pi(k)H(p)\pi(k^{-1})$ for all $k\in K$, $p\in\HH$. The
relation among the eigenvalues is given in Proposition
\ref{relacionautovalores}.

\section{Reduction to one variable}\label{onevariable}

We are interested in considering
the differential operators $D$ and $E$ acting on functions $H\in
C^\infty(\HH)\otimes \End(V_\pi)$ such that
$$H(k\cdot p)=\pi(k)H(p)\pi(k)^{-1}, \quad
\text{ for all $k\in K$ and $p\in\HH$. }$$
  This property of $H$ allows us to find ordinary differential operators
$\tilde D$
and $\tilde E$ acting on functions defined on the interval
$(0,1)$ such that
$$(D\,H)(0,r)=(\tilde D\tilde H)(r),\qquad (E\,H)(0,r)=(\tilde E\tilde
H)(r),$$
where $\tilde H(r)=H(0,r)$.

\smallskip
We need to know the $K$-orbit structure of $\HH$:
The orbits of $K$ in $\HH$ are the spheres with center in the positive axis
$\vzm{(0,r)}{r>0}$ with north and south poles of
the form $(0,s)$ and $(0,s^{-1})$, and the single point set
$\vz{(0,1)}$. The point $(0,1)$ will be denoted from now on by
$o=(0,1)$.
Thus the set of $K$-orbits in $\HH$ is parametrized by the interval
$(0,1]$ or by $[1,\infty)$.

We introduce now two important characters in the story, $\HH^{\times}$ and
$M$.
Let $\HH^\times=\HH-\{o\}$ and put $M=\{k\in K: k\cdot (0,r)=(0,r)\}$,
$0<r<1$.
For any such $r$ we get
$$M=\left\{\left(\begin{matrix}e^{i\theta}&0\\0&e^{-i\theta}\end{matrix}
\right):\theta\in\RR\right\}.$$
Notice that $M$ coincides with the centralizer in $K$ of the Abelian
subgroup $A$ of $G$
given by
$$A=\left\{\left(\begin{matrix}s&0\\0&s^{-1}\end{matrix}
\right):s>0\right\}.$$

Let $\psi:K/M\times (0,1)\longrightarrow\HH^\times$ be defined by
\begin{equation}\label{psi}
\psi(kM,s)=k\cdot(0,s).
\end{equation}
 From the orbit structure discussed above it follows that $\psi$ is a
$C^\infty$-bijection.
To give the expression of the operators $\tilde D$ and $\tilde E$,
obtained by restriction of $D$ and $E$ given in Propositions
\ref{Dexpresion} and \ref{Eexpresion}, we need to compute a number of first
and second order
partial derivatives of the function $H$ at the point $(0,r)$.

To begin with we note the following immediate results
$$H_{r}(0,r)=\frac{d\tilde H}{dr}(r) \quad \text{ and }\quad
H_{rr}(0,r)=\frac{d^2\tilde H}{dr^2}(r).$$
To compute the other partial derivatives of the function $H$ in $\HH^\times$
we need to define locally a pair of $C^\infty$-functions on $\HH^\times$,
$k=k(z,r)$
and $s=s(z,r)$ with values in $K$ and in the open interval $(0,1)$,
respectively,
such that $k\cdot(z,r)=(0,s)$.

In the open set $\HH-\{(0,r):r\ge1\}$ these functions will be explicitly
given below.
Similarly one can define a pair of functions $k$, $s$ in the open set
$\HH-\{(0,r):0<r\le1\}$.

None of these functions can be extended to all of $\HH^\times$ in a
continuous fashion.
Moreover there is no such a pair of continuous globally defined functions
$k$, $s$.
This would imply the existence of a continuous section of $K/M$, which would
make $K$ homeomorphic to $K/M\times M$. Now $K$ is simply connected and $M$
is not.

If $(z,r)\in\HH-\{(0,r):r\ge1\}$ we look for $k=\matc{a}{b}{-\bar b}{\bar
a}\in \SU(2)$
such that $k\cdot(z,r)=(0,s)$, $0<s<1$. This holds if and only if
$$ k\left(\begin{matrix}\frac {1+|z|^2}r&z\\\bar z&r\end{matrix}
\right)k^* =\left(\begin{matrix}s^{-1}&0\\0&s\end{matrix} \right).$$
Then $s$ and $s^{-1}$ are the eigenvalues of $\left(\begin{matrix}\frac
{1+|z|^2}r&z\\\bar z&r\end{matrix} \right)$.
The characteristic equation is
$$r x^2-x(1+|z|^2+r^2)+r=0,$$
thus
$$\vz{s,s^{-1}}=\vz{\tfrac 1{2r}\left(1+|z|^2+r^2 \pm
\sqrt{(|z|^2+(r+1)^2)(|z|^2+(r-1)^2)}\right)}.$$

\noindent Therefore the eigenvalue that satisfies $0<s<1$ is
\begin{equation*}
s=s(z,r)=\frac 1{2r}\left(1+|z|^2+r^2 - \sqrt{(|z|^2+(r+1)^2)(|z|^2+
(1-r)^2)}\right).
\end{equation*}
We have thus defined the function $s$. Now we can take as a function
$k=k(z,r)$
the one given by taking $k=\matc{a}{b}{-\bar b}{\bar a}\in \SU(2)$ with
\begin{equation*}
a=\frac{1-rs}{\sqrt{(1-rs)^2+(s\,|z|)^2}} \qquad \text{and}\qquad
b=\frac{s\,z}{\sqrt{(1-rs)^2+(s\,|z|)^2}}.
\end{equation*}
The explicit expressions above for $s$ and $k$ give as a pair of
$C^\infty$-functions
with the desired properties in the appropriate domain.

The inverse of the $C^\infty$-bijection $\psi$, defined in \eqref{psi}, in
$\HH-\{(0,r):r\ge1\}$
is given by
$$\psi^{-1}(z,r)=(k(z,r)^{-1}M,s(z,r))$$
making clearly that $\psi:K/M\times (0,1)\longrightarrow\HH^\times$ is a
diffeomorphism.

Returning now to the main goal of this section, let
$\left(C^\infty(\HH^\times)\otimes \End(V_\pi)\right)^K$ be the space of all
$C^\infty$-functions
$H$ on $\HH^\times$ with values in $\End(V_\pi)$ such that
$$H(k\cdot p)=\pi(k)H(p)\pi(k^{-1})\quad\text{for all}\quad k\in K,\,
p\in\HH^\times.$$
We will also need to consider the space $C^\infty((0,1))\otimes
\End_M(V_\pi)$ of all
$C^\infty$-functions $\tilde H$ on the open interval $(0,1)$ with values in
$\End(V_\pi)$
with the extra requirement that
$$\tilde H(s)=\pi(m)\tilde H(s)\pi(m^{-1})\quad\text{for all}\quad m\in M,\,
0<s<1.$$

Now it is clear that the restriction $H\longrightarrow\tilde H$ defines an
injective
linear map from $\left(C^\infty(\HH^\times)\otimes \End(V_\pi)\right)^K$
into
$C^\infty((0,1))\otimes \End_M(V_\pi)$. Moreover this map is a surjective
isomorphism
as follows from the following construction. Given
$\tilde H\in C^\infty((0,1))\otimes \End_M(V_\pi)$ define
$H(z,r)=\pi(k)\tilde H(s)\pi(k^{-1})$ if $(z,r)=\psi(kM,s)$.

The differential operators $D$, and $E$, being invariant
under the action of $K$ in $\HH^\times$, leave stable the subspace
$\left(C^\infty(\HH^\times)\otimes \End(V_\pi)\right)^K$. Thus we have the
following commutative
diagram

$$\begin{CD}
\left(C^\infty(\HH^\times)\otimes \End(V_\pi)\right)^K @ >
>>C^\infty((0,1))\otimes \End_M(V_\pi)\\ @ V
D,E VV @ VV \tilde D,\tilde E V \\\left(C^\infty(\HH^\times)\otimes
\End(V_\pi)\right)^K @> >>
C^\infty((0,1))\otimes \End_M(V_\pi)
\end{CD}$$

A spherical function $\Phi$ of type $\pi\in\hat K$ gives rise to a function
$H\in
\left(C^\infty(\HH^\times)\otimes \End(V_\pi)\right)^K$ which is an
eigenfunction of $D$
(also of $E$). In turn such an $H$ corresponds to an $\tilde H\in
C^\infty((0,1))\otimes \End_M(V_\pi)$
which is an eigenfunction of $\tilde D$ (also $\tilde E$) with the extra
condition
$\lim_{s\rightarrow1^-}\tilde H(s)=I$. Conversely if $\tilde H$ is such a
function  then it is the restriction of a unique eigenfunction
$H\in \left(C^\infty(\HH^\times)\otimes \End(V_\pi)\right)^K$ of $D$ (also
of $E$) such that
$\lim_{p\rightarrow o} H(p)=I$. From Proposition \ref{Dexpresion} it is
clear that
$D$  is an elliptic differential operator, thus $H$ can be extended to a
$C^\infty$-function on $\HH$.

Now that the overall strategy has been outlined we go on to give a variety
of explicit expressions of
$\tilde D$ and $\tilde E$. See Corollaries \ref{sistema}, \ref{sistema2},
\ref{DandEint} as well
as \eqref{DD} and \eqref{EE} below.

\

The proof of the following proposition is included, for completeness,
  in an Appendix at the end of the paper.

\begin{prop}\label{dderxy}
Let $J=\begin{pmatrix}0&1\\-1&0\end{pmatrix}$ and $T=\begin{pmatrix}
0&1\\ 1&0
\end{pmatrix}$. We have, for $0<r<1$,
\begin{align*}
H_{x}(0,0,r)&= -{\frac r{1-r^2}}\left( \dot\pi(J)\tilde H(r)-\tilde
H(r)\dot\pi(J)\right),\\
H_{y}(0,0,r)&= -{\frac {i\,r}{1-r^2}}\left( \dot\pi(T)\tilde H(r)-\tilde
H(r)\dot\pi(T)\right),
\end{align*}
and
\begin{equation*}
\begin{split}
H_{xx}(0,0,r)&= -\frac {2r}{1-r^2} \frac{d \tilde H}{d r} +
\frac{r^2}{(1-r^2)^2}\left(\dot\pi(J)^2\, \tilde H(r)+ \tilde
H(r)\,
\dot\pi(J)^2\right)\displaybreak[0]\\
&\quad -\frac{2r^2}{(1-r^2)^2}\dot\pi(J)\tilde H(r)\, \dot\pi(J), \\
H_{yy}(0,0,r)&= -\frac {2r}{1-r^2} \frac{d \tilde H}{d r}
- \frac{r^2}{(1-r^2)^2}\left(\dot\pi(T)^2\, \tilde H(r)+ \tilde H(r)\,
\dot\pi(T)^2\right)\displaybreak[0]\\
&\quad +\frac{2 r^2}{(1-r^2)^2}\dot\pi(T)\tilde H(r)\, \dot\pi(T).
\end{split}
\end{equation*}
\end{prop}

\

\begin{thm}\label{mainth1}
We have, for $0<r<1$,
\begin{align*}
\textstyle \frac 14 \tilde D\tilde H(r)&=r^2\frac{d^2 \tilde H}{dr^2}
-\frac{r(1+3\,r^2)}{1-r^2} \frac{d\tilde H}{dr}\displaybreak[0]
\\ &
+\frac {r^2}{(1-r^2)^2}\left( \dot\pi(J)^2\,\tilde H(r)\,+\tilde
H(r)\dot\pi(J)^2 -2\dot\pi(J)\tilde H(r)\dot\pi(J)\right)\\
& -\frac {r^2}{(1-r^2)^2}\left( \dot\pi(T)^2\,\tilde H(r)\,+\tilde
H(r)\dot\pi(T)^2 -2\dot\pi(T)\tilde H(r)\dot\pi(T)\right),
\end{align*}
where $J=\begin{pmatrix}0&1\\-1&0\end{pmatrix}$ and
$T=\begin{pmatrix}0&1\\1&0\end{pmatrix}$.
\end{thm}
\begin{proof} By Proposition \ref{Dexpresion} we have
$$\tfrac 14
DH(0,0,r)=H_{xx}(0,0,r)+H_{yy}(0,0,r)+r^2H_{rr}(0,0,r)+3rH_{r}(0,0,r).$$
Using  Proposition \ref{dderxy} the theorem follows.
\end{proof}

\medskip

\begin{thm}\label{mainth2}
We have
\begin{align*}
\textstyle  \tilde E\tilde H(r)&=
-2r  \frac{d\tilde H}{dr}\dot\pi(H_1)+\frac {4}{1-r^2}\left(
\dot\pi(X_2)\,\tilde H(r)\,-\tilde
H(r)\dot\pi(X_2)\right)\dot\pi(X_1)  \\
& - \frac {4r^2}{1-r^2}\left( \dot\pi(X_1)\,\tilde H(r)\,-\tilde
H(r)\dot\pi(X_1)\right)\dot\pi(X_2),
\end{align*}
where $H_1=\begin{pmatrix}1&0\\0&-1\end{pmatrix}$,
$X_1=\begin{pmatrix}0&1\\0&0\end{pmatrix}$ and
$X_2=\begin{pmatrix}0&0\\1&0\end{pmatrix}$.
\end{thm}

\begin{proof}
By Proposition \ref {Eexpresion} we have
\begin{align*}
EH(0,0,r)&= H_x \;\dot\pi\matc{0}{\;\;\frac{2}r}{2\,r}{0} +
H_y \;\dot\pi\matc{0}{\;\;\frac{2i}r}{-2i\,r}{0} +H_r
\;\dot\pi\matc{-2r}{0}{0}{\;\;2r}.
\end{align*}
Then using Proposition \ref{dderxy} we obtain
\begin{align*}
\tilde E\tilde H=&\textstyle -2r\frac{d\tilde
H}{dr}\dot\pi(H_1)-\frac{2r}{1-r^2}\left(\dot\pi(J)\tilde H(r)-\tilde
H(r)\dot\pi(J)\right) \dot\pi\matc{0}{\;\frac 1r}{r}{0}
\displaybreak[0]\\ &
\quad
+\textstyle \frac{2r}{1-r^2}\left(\dot\pi(t)\tilde H(r)-\tilde
H(r)\dot\pi(t)\right) \dot\pi\matc{0}{1/r}{-r}{0}\\
&=\textstyle -2r\frac{d\tilde
H}{dr}\dot\pi(H_1)-\frac{2}{1-r^2}\left(\dot\pi(J)\tilde H(r)-\tilde
H(r)\dot\pi(J)\right)\dot\pi(X_1)\\
&\quad  -\tfrac{2r^2}{1-r^2}\left(\dot\pi(J)\tilde H(r)-\tilde
H(r)\dot\pi(J)\right)\dot\pi(X_2) \\
& \quad +\tfrac{2}{1-r^2}\left(\dot\pi(T)\tilde H(r)-\tilde
H(r)\dot\pi(T)\right)\dot\pi(X_1)\\
&\quad -\tfrac{2r^2}{1-r^2}\left(\dot\pi(T)\tilde H(r)-\tilde
H(r)\dot\pi(T)\right)\dot\pi(X_2).
\end{align*}
Now using that $J=X_1-X_2$ and $T=X_1+X_2$ the theorem follows.
\end{proof}

\

The theorems above are given in terms of linear transformations. Now we
will
give the corresponding statements in terms of
matrices by choosing an appropriate basis. If $\pi=\pi_{\ell}$ it is
well
known (see \cite{Hu}, p. 32) that there exists a
basis $\vz{v_i}_{i=0}^\ell$ of $V_\pi$ such that
\begin{align*}
&\dot\pi(H_1)v_i=(\ell-2i)v_i, &\displaybreak[0]\\
&\dot\pi(X_1)v_i=(\ell-i+1)v_{i-1}, \quad &
(v_{-1}=0), \displaybreak[0]\\
&\dot\pi(X_{2})v_i=(i+1)v_{i+1}, \quad &
(v_{\ell+1}=0),\\
\end{align*}
where  $H_1=\begin{pmatrix}1&0\\0&-1\end{pmatrix}$,
$X_1=\begin{pmatrix}0&1\\0&0\end{pmatrix}$ and
$X_2=\begin{pmatrix}0&0\\1&0\end{pmatrix}$.

\begin{lem}\label{Hdiagonal}
  The function $\tilde H$ is diagonalizable.
\end{lem}
\begin{proof}
The subgroup $M$ of $K$,
$$M=\vzm{m=\left(\begin{smallmatrix} e^{it}& 0 \\0&
e^{-it}\end{smallmatrix}\right)}{t\in \RR},$$
  fixes the points $(0,r)$ in the hyperbolic space $\HH=\vzm{(z,r)}{z\in
\CC, r\in\RR{^+}}$. We have
$H(kg)=\pi(k)H(g)\pi(k^{-1})$, therefore $\tilde
H(r)=\pi(m)\tilde H(r)\pi(m^{-1})$. Now since all finite
dimensional
irreducible $K$-modules are multiplicity free as
$M$-modules, $\tilde H(r)$, $r>0$, and $\dot\pi(H_1)$ diagonalize
simultaneously.
\end{proof}

\medskip
\noindent We introduce the functions $\tilde h_j(r)$ by means of the
relations
\begin{equation}
\tilde H(r)v_j=\tilde h_j(r)v_j.
\end{equation}

\smallskip
\begin{prop}
The functions $\tilde h_j$ defined above satisfy
$$\tilde h_j(r)=\tilde h_{\ell-j}(r^{-1}),$$
for all  $0<r<\infty$.
\end{prop}

\begin{proof} From $m=\matc{e^{it}}{0}{0}{e^{-it}}=\exp tiH_1$ it follows
that $\pi(m)v_j=e^{it(\ell-2j)}v_j$. Let
$b=\matc{0}{1}{-1}{0}\in K. $
If $m=\exp tiH_1$ we have that $mb=bm^{-1}$, thus
$\pi(m)\pi(b)v_j=e^{-it(\ell-2j)}\pi(b) v_j$. Therefore
$\pi(b)v_j=c_jv_{\ell-j}.$
Since  $b^2=\exp \pi iH_1$ it follows that $\pi(b^2)=(-1)^\ell I$, hence
$(-1)^\ell v_j=c_jc_{\ell-j}v_j$, so we have
$$c_jc_{\ell-j}=(-1)^\ell\, , \quad  \text{ for } j=0, \dots ,\ell.$$
Let $a_s=\begin{pmatrix}s&0\\0&s^{-1}\end{pmatrix}\in G$, $s>0$. Then
$$H(a_s)=H(p(a_s))=H(0, s^{-2})= \tilde H(s^{-2}).$$
On the other hand $H(a_{s^{-1}})=H(b^{-1}a_sb)= \pi(b^{-1})H(a_s)\pi(b)$,
thus
\begin{align*}
\tilde H(s^2)v_j=& \pi(b^{-1})\tilde H(s^{-2})\pi(b) v_j= c_j \tilde
h_{\ell-j}(s^{-2})\pi(b^{-1})v_{\ell-j}\\
=& (-1)^\ell
c_jc_{\ell-j} \tilde h_{\ell-j}(s^{-2})v_j,
\end{align*}
therefore
$$\tilde h_j(s^2)= (-1)^\ell c_jc_{\ell-j} \tilde h_{\ell-j}(s^{-2})=
\tilde h_{\ell-j}(s^{-2})$$
\end{proof}

\begin{cor}\label{sistema}
The function $\tilde H(r)=(\tilde h_0(r), \cdots, \tilde h_\ell(r) )$,
$(0<r<1)$,
satisfies $(\tfrac 14\tilde D\tilde H)(r)=\ld \tilde H(r)$ if
and
only if
\begin{align*}
  r^2 &\tilde h_i''-\frac{r(1+3\,r^2)}{1-r^2} \tilde h_i'
\displaybreak[0]\\
  & +\frac{4\, r^2}{(1-r^2)^2}\bigl(i(\ell-i+1)(\tilde h_{i-1}-\tilde h_{i})+
(i+1)(\ell-i)(\tilde h_{i+1}-\tilde h_i) \bigr) =\ld\tilde h_i
\end{align*}
for all $i=0,\cdots, \ell$.
\end{cor}

\

\begin{cor}\label{sistema2}
The function $\tilde H(r)=(\tilde h_0(r), \cdots, \tilde h_\ell(r) )$,
$(0<r<1)$,
satisfies $(\tilde E\tilde H)(r)=\mu \tilde H(r)$ if
and only if
\begin{align*}
& -2r(\ell-2i) \tilde h_i'
  +\frac{4\, i(\ell-i+1) }{1-r^2}(\tilde h_{i-1}-\tilde h_{i})
-\frac{4\,r^2(i+1)(\ell-i)}{1-r^2} (\tilde h_{i+1}-\tilde h_i) \\& \quad
=\mu\tilde h_i
\end{align*}
for all $i=0,\cdots, \ell$.
\end{cor}

We introduce the change of variable $t=r^2$ and we put $H(t)=\tilde H(\sqrt
t)$ and $h_i(t)=\tilde h_i(\sqrt t)$. The differential operators
$\tilde D$ and $\tilde E$ are converted into new differential operators $D$
and $E$. At this point there is a slight abuse of notation, since $D$ and
$E$
were used earlier to denote operators on $\HH$.
We get the following expressions for these new differential operators $D$
and $ E$.

\begin{cor}\label{DandEint}
For $0<t<1$, the function $\tilde H$ satisfies $(D H)(t)=4\ld  H(t)$ and $(
E H)(t)=\mu  H(t)$ if and
only if
\begin{align*}
4t^2 h_i''-\frac{8t^2}{1-t} h_i'
  &+\frac{4t}{(1-t)^2}i(\ell-i+1)(h_{i-1}-h_{i})\\
&+\frac{4t}{(1-t)^2}(i+1)(\ell-i)(h_{i+1}-h_i)=\ld h_i,\\
4(\ell-2i) t h_i'&+\frac{4 }{1-t}i(\ell-i+1)(h_{i-1}-h_{i})\\
&-\frac{4\,t}{1-t} (i+1)(\ell-i)(h_{i+1}-h_i)=\mu h_i,
\end{align*}
respectively.
\end{cor}

In  matrix notation, the differential operators $D$ and $E$ are given by
\begin{align}
DH&=4\,t^2 \,H''-\frac{8\,t^2}{1-t}\,H'+\frac{4\,t}{(1-t)^2}(C_0+C_1)H.
\label{DD}\displaybreak[0]\\
EH&=-4\,t\, A_0 H'+\frac{4}{1-t}C_0 H-\frac{4\,t}{1-t}C_1 H. \label{EE}
\end{align}
Here $H$ denotes the column vector
$H=(h_0,\dots,h_\ell)$ and the  matrices are given by:
\begin{align*}
A_0&=\textstyle\sum_{i=0}^\ell(\ell-2i)E_{ii},\displaybreak[0]\\
C_0&=\textstyle\sum_{i=0}^\ell
i(\ell-i+1)(E_{i,i-1}-E_{i,i}),\displaybreak [0]\\
C_1&=\textstyle\sum_{i=0}^\ell(i+1)(\ell-i)(E_{i,i+1}-E_{i,i}).
\end{align*}

\begin{remark}\label{relacionautovalores2}
The function $\Phi\in C^\infty(G)\otimes\End(V_\pi)$
satisfies  $(4\overline \Omega) \Phi=\tilde\lambda \Phi$ and
$(\Omega-\overline\Omega)\Phi=\tilde\mu \Phi$
  if and only if the function $H=H(t)$ associated with $\Phi$ satisfies
$DH=\lambda H$ and $E H=\mu H$, with
$$\tilde\lambda =4\lambda\quad \text{and}
\quad \tilde \mu =\mu +\ell(\ell+2),$$
(see Proposition \ref{relacionautovalores}).
\end{remark}

\section{Eigenfunctions of $D$}\label{D}

The goal of Sections \ref{D} and \ref{DE} is:

\noindent(i) to describe all solutions to $DH=\lambda H$,
$EH=\mu H$, $0<t<1$,

\noindent(ii) to identify among these, those $H$ that have a limit as
$t\rightarrow1^-$. These correspond precisely to the spherical
functions, as seen in Section \ref{onevariable}.

We are interested in considering  first  the functions
$H:(0,1)\longrightarrow \CC^{\ell+1}$ such that $DH=\ld H$, $\lambda\in\CC$.
It is well
known that such eigenfunctions are analytic functions on the open interval
$(0,1)$ and that the dimension of the corresponding eigenspace $V_\ld$
is $2(\ell+1)$. This space will be decomposed in different ways as a direct
sum of $\ell+1$ dimensional spaces stable under $E$. See Propositions
\ref{Vlambda}, \ref{complementoV(p)} and \ref{complementoV(p)2}.
The action of $E$ on these spaces will be taken up in Section \ref{DE}.

\

\noindent To determine the eigenfunctions of the differential operator
$$DH=4\,t^2
\,H''-\frac{8\,t^2}{1-t}\,H'+\frac{4\,t}{(1-t)^2}(C_0+C_1)H$$
for $0<t<1$ it is necessary to take a close look at it and this is the
purpose of the next proposition.

Since the $(\ell+1)\times (\ell+1)$ matrix $C_0+C_1$ is not a diagonal
matrix, the equation $DH=\ld H$ is a coupled system of $\ell+1$ second
order
differential equations in the components $(h_0, \dots , h_\ell)$ of
$H$.
But fortunately the matrix $C_0+C_1$ is a symmetric one, thus
diagonalizable.

\begin{prop}\label{Hahn} The matrix $C_0+C_1$ is diagonalizable.
Moreover
the eigenvalues are $-j(j+1)$ for $0\le j\le \ell$
and the corresponding eigenvectors are given by
$u_j=(u_{0,j},\dots,u_{\ell,j})$ where
\begin{equation}\label{Ucolumnas}
u_{i,j}= \lw{3}F_2\left( \begin{smallmatrix}
-j,\;-i,\;j+1 \\ 1,\;-\ell \end{smallmatrix}; 1 \right),
\end{equation}
an instance of the Hahn orthogonal polynomials.
\end{prop}

\begin{proof}
We shall see that
\begin{equation}\label{CC}
(C_0+C_1)u_j=-j(j+1) u_j\quad\text{for all}\quad 0\leq j\leq \ell.
\end{equation}

The matrix $C_0+C_1$ is given by
$$\sum_{i=0}^\ell
i(\ell-i+1)E_{i,i-1}- \left(i(\ell-i+1)+ (i+1)(\ell-i)\right) E_{i,i}
+(i+1)(\ell-i) E_{i,i+1}.$$
Then the $i$th-component of the equation \eqref{CC} is given by
\begin{align*}
i(\ell-i+1) u_{i-1,j}&- \left(i(\ell-i+1)+ (i+1)(\ell-i)\right) u_{i,j}
\\& + (i+1)(\ell-i) u_{i+1,j}=-j(j+1)u_{i,j}.
\end{align*}
This allows us to recognize the eigenvectors of $C_0+C_1$ as an instance of
the Hahn polynomials.

For real numbers $\al, \beta >-1$,  and for a positive integer $N$ the Hahn
polynomials $Q_n(x)=Q_n(x;\al, \beta, N)$ are defined by
$$Q_n(x)=\lw{3}F_2\left( \begin{smallmatrix}
-n,\;-x,\;n+\al+\beta+1 \\ \al+1,\;-N+1 \end{smallmatrix}; 1 \right),
\qquad \text{ for } n=0,1,\dots , N-1.$$
These Hahn polynomials are examples of orthogonal polynomials and hence
satisfy a three term recursion relation, see \cite{KMcG}, (1.3) or
\cite{AAR}.
If one takes $\al=\beta=0$, $x=i$, $N=\ell+1$, $n=j$  we obtain
$$ u_{i,j}= Q_j(i)= \lw{3}F_2\left( \begin{smallmatrix}
-j,\;-i,\;j+1 \\ 1,\;-\ell \end{smallmatrix}; 1 \right).$$

Since all the eigenvalues are different we have that $C_0+C_1$ is
diagonalizable.
\end{proof}

\

Let $U=(u_{i,j})$ be the matrix defined by \eqref{Ucolumnas}. Then
\begin{equation}\label{UU-1}
U\Bigg(-\sum_{j=0}^\ell j(j+1)E_{jj}\Bigg) U^{-1}=C_0+C_1.
\end{equation}
If we define $\check H(t)=U^{-1}H(t)$, we see that the $i$th-component
$\check h_i(t)$ of $\check H(t)$ satisfies
\begin{equation}\label{eqtilde}
4\,t^2 \,\check h_i''(t)-\frac{8\,t^2}{1-t}\,\check
h_i'(t)-4i(i+1)\frac{t}{(1-t)^2}\check h_i(t)-\ld \check h_i(t)=0
\end{equation}
if and only if $DH=\ld H$.

\

Recall that
$$V_\ld=\vzm{H=H(t)}{DH=\ld H(t), \, 0<t<1}.$$
Suppose that $0\neq H\in V_\ld$ is of the form $H(t)=t^p F(t)$ with $F$
analytic at $t=0$ and $p\in\CC$. In such a case we may assume that $F(0)\neq
0$. We
first note the following relation between $\ld$ and $p$.

\begin{lem} If $H=t^p F(t)$ with $F$ analytic at $t=0$, and $F(0)\neq
0$
satisfies $DH=\lambda H$, then $\lambda=4p(p-1)$.
\end{lem}
\begin{proof}
Let $F(t)=\sum_{j=0}^\infty\,t^jF_j$, then
$F$ satisfies
\begin{align*}
4t^2F''+\frac{8t}{1-t}&\left(p-(p+1)t\right)F'+\frac{4t}{(1-t)^2}(-2p(1-t)+C
_0+C_1)F\\
& +4p(p-1)F -\ld F=0.
\end{align*}

If we  multiply the above expression  by $(1-t)^2$ and set $t=0$ we
obtain
$(4p(p-1)-\ld)F(0)=0$ and the lemma follows.
\end{proof}

\

Observe that the function $\ld(p)=4p(p-1)$ satisfies $\ld(p)=\ld(1-p)$.
Given $\ld\in \CC$, if $p\in \CC$ satisfies $\ld=4p(p-1)$ let
\begin{equation}\label{V(p)def}
V(p)=\vzm{H\in V_\ld}{H(t)=t^pF(t), \, F \text{ analytic at } t=0}.
\end{equation}

We note that if $H(t)=t^pF(t)\in V(p)$, with $F(t)=\sum_{j=0}^\infty t^jF_j$
then the vector coefficients $F_j$ satisfy the three term recursion
relation
\begin{equation}\label{recur1}
\begin{split}
&j(j-1+2p)F_j-\bigl( 2p(2j-1)+2(j-1)^2- C_0-C_1)\bigr) \,F_{j-1}\\
&\quad +(j-1)(2p+j-2)\,F_{j-2}=0,\quad\quad j\ge1.
\end{split}
\end{equation}

\medskip

\begin{prop} \label{Hforma} Let $H\in V_\ld$, $\ld=4p(p-1)$ and let
$\check H(t)=U^{-1}H(t)$. If we write
$$\check H(t)=\frac{t^p}{(1-t)^{\ell+1}}\check  P(t),$$
then the $i$th-component $\check  P_i=\check  P_i(t)$  of $\check  P$  is
of the form $\check  P_i(t)=(1-t)^{\ell -i}R_i(t)$ where $R_i$ is a
solution of
the differential equation
\begin{equation}\label{ecuacionR}
t(1-t)R_i''+(2p-(2(p-i)t) R_i'+ i(2p-i-1)R_i=0.
\end{equation}
\end{prop}

\rmk{Remark} Notice that \eqref{ecuacionR} is nothing but the Gauss'
hypergeometric equation with parameters $a=-i$, $b=2p-i-1$, $c=2p$.

\medskip
\begin{proof} If $H\in V_\ld$ we can write
$H(t)=\frac{t^p}{(1-t)^{\ell+1}}
P(t)$ with $P$ analytic on the interval $0<t<1$. Then it is easy to
verify
that $P$ satisfies
\begin{equation}\label{ecuacionP}
\begin{split}
4t^2P''&+\frac{4t}{1-t}\left(2p-2(p-\ell) t\right)P'\\
&+\frac{4t}{(1-t)^2}
\left(\ell(\ell+1)t+2p\ell(1-t)+C_0+C_1 \right)P=0.
\end{split}
\end{equation}
If we define the function  $\check P=U^{-1}P$ it  follows from
\eqref{ecuacionP} that the $i$th-component $\check P_i$  of $\check P$
satisfies the following differential equation
\begin{equation*}
\begin{split}
t(1-t)\check P_i''&+\left(2p-2(p-\ell) t\right)\check P_i'\\
&+\frac{1}{(1-t)}
\left(\ell(\ell+1)t+2p\ell(1-t)-i(i+1) \right)\check  P_i=0.
\end{split}
\end{equation*}
We can also set $\check P_i=(1-t)^{\ell-i}R_i$ and then it is easy to
verify that $R_i$ satisfies
\begin{equation*}
t(1-t)R_i''+(2p-(2(p-i)t) R_i'+ i(2p-i-1)R_i=0.
\end{equation*}
This finishes the proof of the proposition.
  \end{proof}

\begin{prop}\label{2pnoentero} If $\lambda=4p(p-1)$ with $2p\in \CC$ not an
integer, and
$H\in V(p)$, $H\ne0$, then
$$H(t)=\frac{t^p}{(1-t)^{\ell+1}}P(t),$$
where $P=P(t)$ is a polynomial function of degree $\ell$.
More precisely $P(t)=U\check P(t)$ where the $i$th-component of $\check P$
is
$$\check P_i=\alpha_i (1-t)^{\ell-i} \,\lw{2}F_1 \left( \begin{smallmatrix}
-i,\,\;2p-i-1 \\ 2p \end{smallmatrix}; t \right)$$
with  $\al_i\in \CC$.
\end{prop}

\begin{proof} By hypothesis, for any $H\in V(p)$,  we may  write
$H(t)=\frac{t^p}{(1-t)^{\ell+1}} P(t)$ with  $P$  analytic at $t=0$.
Let $\check H=U^{-1}H$ and  $\check P=U^{-1}P$. Then by Proposition
\ref{Hforma} the $i$th-component $\check P_i$ of $\check  P$ is of the
form
$\check P_i=(1-t)^{\ell-i}R_i$ where $R_i$ satisfies the Gauss'
hypergeometric equation \eqref{ecuacionR}. Moreover $R_i$ is analytic
at $t=0$.

Since $2p\not\in \ZZ$ it is easy to verify that a basis of the
solutions of
\eqref{ecuacionR} is given  by the functions:
\begin{equation}\label{basegeneral}
\lw{2}F_1 \left( \begin{smallmatrix} -i,\,\;2p-i-1 \\
2p \end{smallmatrix}; t \right)
\, , \qquad  t^{1-2p} \lw{2}F_1\left( \begin{smallmatrix}
-i,\,\;-2p-i+1 \\
2-2p \end{smallmatrix}; t \right).
\end{equation}
Therefore
$$R_i=\al_i \,\lw{2}F_1 \left( \begin{smallmatrix} -i,\,\;2p-i-1 \\
2p \end{smallmatrix}; t \right)$$
for some $\al_i\in \CC$.
Hence each $\check P_i$ is a polynomial function of degree $\ell$,
whenever $\al_i\neq 0$.
Since $P=U\check  P$ the assertions made in the statement of the
proposition are now clear.
\end{proof}

\begin{cor}\label{dim2pnoentero} If $2p\in \CC$ is not an integer then
$\dim V(p)=\ell+1$.
\end{cor}

\medskip

\begin{prop} \label{Vlambda} If $2p\in \CC$ is not an integer and
$\ld=4p(p-1)$, then
$$V_\ld=V(p)\oplus V(1-p).$$
\end{prop}

\begin{proof}
Let $H\in V(p)\cap V(1-p)$. If $H\neq 0$ then $H=t^pF(t)=t^{1-p}G(t)$
with
$F$ and $G$ analytic at $t=0$. Since $2p\notin \ZZ$ from \eqref{recur1} we
get $F(0)\neq 0$.
Thus $t^{2p-1}=\frac{G(t)}{F(t)}$ is analytic at $t=0$, but $t^{2p-1}$
has a ramification point at $t=0$, because $2p\not\in \ZZ$. This
contradiction
proves that $V(p)\cap V(1-p)=0$.

Since $\dim V(p)=\dim V(1-p)=\ell+1$ (Corollary \ref{dim2pnoentero}) the
proof of the proposition is completed.
\end{proof}

\medskip
When $2p\in \ZZ$ the situation is completely different as shown in the next
proposition. The reader may choose to skip the computation of $\dim V(1-p)$
in i) and of $\dim V(p)$ in ii).

\begin{prop}\label{dimV(p)}
  Let $2p\in \ZZ$ and $\ld=4p(p-1)$.
\begin{enumerate}
\item[i)] If  $2p\geq 1$ then we have $V(p)\subset V(1-p)\subset
V_\ld$,
$\dim V(p)=\ell+1$  and  $\dim V(1-p)=\min \vz{2(\ell+1), \ell+2p}.$

\smallskip
\item[ii)] If  $2p\leq 0$ then $V(1-p)\subset V(p)\subset V_\ld$, $\dim
V(1-p)=\ell+1$ and $\dim V(p)=\min \vz{2(\ell+1), \ell+2-2p}$.
\end{enumerate}
\end{prop}

\begin{proof}
Let us first observe that i) and ii) are equivalent by changing $p$ by
$1-p$. Thus we can assume that  $2p\geq 1$.

If $H\in V(p)$ then $H(t)=t^pF(t)= t^{1-p}\left(t^{2p-1}F(t)\right)$
with $t^{2p-1}F(t)$ analytic at $t=0$. This proves that $V(p)\subset
V(1-p)$.

If $H\in V_\ld$ and
$\check H=U^{-1}H $, then by Proposition \ref{Hforma}, the
$i$th-component
$\check h_i$ of $\check H$ is of the form $\check
h_i=\frac{t^p}{(1-t)^{i+1}}R_i$ where $R_i$ is  a solution of the Gauss'
equation with parameters $a=-i$, $b=2p-i-1$, $c=2p$.

\noindent A convenient basis of  solutions of this Gauss' equation for
$0<t<1$, is given  by the functions (see \cite{Bat}, Section 2.2.2,
Case 22)
\begin{equation}\label{basebuena}
\lw{2}F_1 \left( \begin{smallmatrix} -i,\,\;2p-i-1 \\
2p \end{smallmatrix}; t \right)
\, , \qquad  t^{1-2p} \lw{2}F_1\left( \begin{smallmatrix}
-i,\,\;-2p-i+1 \\
2-2p \end{smallmatrix}; t \right),
\end{equation}
whenever $0\le i<2p-1$. On the other hand if $2p-1\leq i\leq\ell$  a basis
of
solutions is given  by the functions (see \cite{Bat}, Section 2.2.2,
Case 23)
\begin{equation}\label{basemala}
\lw{2}F_1 \left( \begin{smallmatrix} -i,\,\;2p-i-1 \\
2p \end{smallmatrix}; t \right)
\, , \qquad  (1-t)^{2i+1} \lw{2}F_1\left( \begin{smallmatrix}
i+1,\,\;2p+i
\\ 2i+2 \end{smallmatrix}; 1-t \right).
\end{equation}

In both cases it is easy to verify that these functions solve
Gauss' equation and are linearly independent.
In both cases the second elements of these bases are not analytic at
$t=0$ while the other one is. In fact, in  the case $2p>i+1$ the function
$t^{1-2p} \lw{2}F_1\left(-i,-2p-i+1, 2-2p; t\right)$ has a pole at
$t=0$ of
order $2p-1>0$.
When $1\leq 2p\leq i+1$ we have the following identity (\cite{Bat}, Ch.
2, (4)):
$$\lw{2}F_1\left( \begin{smallmatrix} i+1,\,\;2p+i \\ 2i+2
\end{smallmatrix}; z \right) =c
\left(\frac{d}{dz}\right)^{i+2p-1}\left(
(1-z)^i\left(\frac{d}{dz}\right)^{i+1-2p}\left( -z^{-1}\log
(1-z)\right)\right).$$
{}From this it is clear that the function $\lw{2}F_1(  i+1,2p+i,2i+2;
z)$ is
not analytic at $z=1$ because $\log (1-z)$ has a ramification point at
$z=1$. Thus  $(1-t)^{2i+1} \lw{2}F_1\left(  i+1,2p+i, 2i+2 ; 1-t
\right)$
is not analytic at $t=0$.

\smallskip
If $H\in V(p)$ then  $R_i$ is analytic at $t=0$  therefore
$$R_i=\al_i \,\lw{2}F_1 \left( \begin{smallmatrix} -i,\,\;2p-i-1 \\
2p \end{smallmatrix}; t \right)$$
for some $\al_i\in \CC$. Thus $\dim V(p)=\ell+1$.

Now we compute the dimension of $V(1-p)$. We observe that by
Proposition \ref{Hforma} for $H\in V(\ld)$,  the $i$th-component
$\check h_i$ of  $\check H$ is of the form $\check
h_i=\frac{t^p}{(1-t)^{i+1}}R_i$.
We can also obtain that $\check h_i(t)=\frac{t^{1-p}}{(1-t)^{i+1}}S_i$,
where
$S_i$ is a solution of the Gauss' equation with parameters $a=-i$,
$b=-2p-i+1$, $c=2-2p$. We have the relation $S_i=t^{2p-1}R_i$.
Therefore a
basis of the solutions of this Gauss' equation, for $0<t<1$ is given by
$$t^{2p-1}\lw{2}F_1 \left( \begin{smallmatrix} -i,\,\;2p-i-1 \\
2p \end{smallmatrix}; t \right)
\, , \qquad   \lw{2}F_1\left( \begin{smallmatrix} -i,\,\;-2p-i+1 \\
2-2p \end{smallmatrix}; t \right).$$
whenever $0\le i<2p-1$. On the other hand if $2p-1\leq i\leq\ell$ a basis of
solutions is given  by
$$t^{2p-1}\lw{2}F_1 \left( \begin{smallmatrix} -i,\,\;2p-i-1 \\
2p \end{smallmatrix}; t \right)
\, , \qquad  t^{2p-1}(1-t)^{2i+1} \lw{2}F_1\left( \begin{smallmatrix}
i+1,\,\;2p+i \\ 2i+2 \end{smallmatrix}; 1-t \right).$$
In the first case, both solutions are analytic at $t=0$, moreover they are
polynomials, while in the second case, the first solution is analytic and
the other one is not.

\smallskip
If $H\in V(1-p)$ then  $S_i$ is analytic at $t=0$  therefore we have
for $0\le i<2p-1$ that
$$  S_i=\al_i\,t^{2p-1}  \lw{2}F_1 \left( \begin{smallmatrix}
-i,\,\;2p-i-1
\\ 2p \end{smallmatrix}; t \right)  +\beta_i\,
\lw{2}F_1\left( \begin{smallmatrix} -i,\,\;-2p-i+1 \\
2-2p \end{smallmatrix}; t \right)$$
for some $\al_i, \beta_i\in \CC$. For $2p-1\leq i\leq\ell$  we have
$$S_i=\al_i \,t^{2p-1}\, \lw{2}F_1 \left( \begin{smallmatrix}
-i,\,\;2p-i-1 \\
2p \end{smallmatrix}; t \right),$$
for some $\al_i\in \CC$.
{}From the description of the solutions it follows that $\dim V(1-p)=\min
\vz{2(\ell+1),
\ell+2p}$.
\end{proof}

\

For any $p\in \CC$ such that dim $V(p)=\ell+1$ a useful complement
of $V(p)$ in $V_\lambda$ will turn out to be the space $W_\ld$ introduced
below.

A spherical function of type $\pi_\ell \in \hat K$ gives rise to a
function
$H\in V_\ld$ such that $\lim_{t\rightarrow 1^-}H(t)=(1, \dots , 1).$
This leads us to consider the linear space
$$W_\ld =\vzm{H\in V_\ld}{\lim_{t\rightarrow 1^-}H(t) \text{ is
finite}}.$$

\

\noindent We start with a useful lemma. The first two cases considered in
the hypothesis
will be used in Propositions \ref{complementoV(p)2} and
\ref{complementoV(p)},
respectively.

\begin{lem}\label{hipergeometrico}  If $2p\in\CC-\ZZ$,  or
$2p\in\ZZ$ and $2p>i+1$, or $2p\in\ZZ$ and $2p<-i+1$, then we have
\begin{equation}\label{lemhiper}
\begin{split}
\tfrac{(i+1)_{i+1}}{(1-2p)_{i+1}}\, \,\lw{2}F_1\left(
\begin{smallmatrix}
-i,\,\;2p-i-1 \\ 2p \end{smallmatrix}; t \right)
+ &\tfrac{(i+1)_{i+1}}{(2p-1)_{i+1}}\, \,t^{1-2p}\, \lw{2}F_1
\left(\begin{smallmatrix}
-i,\,\;-2p-i+1 \\ 2-2p\end{smallmatrix}; t \right)\\
& =  t^{1-2p}(1-t)^{2i+1}\,
\lw{2}F_1\left( \begin{smallmatrix}
i+1,\;i+2-2p \\ 2i +2 \end{smallmatrix}; 1-t \right).
\end{split}
\end{equation}
\end{lem}
\begin{proof} Let us first assume that $2p\in\CC-\ZZ$ and $\re(2p)>1$, or
that
$2p\in\ZZ$ and $2p>i+1$. Under this hypothesis
on $p$ the  hypergeometric
equation with parameters $a=-i$, $b=2p-i-1$ and $c=2p$, as we pointed out
before, has a basis of solutions given by \eqref{basebuena}.
In general we have that
$t^{1-c}(1-t)^{c-a-b}\lw{2}F_1(1-a,1-b, c+1-a-b; 1-t)$ is a solution of
the hypergeometric equation with parameters $a,b,c$. So, we get
\begin{equation*}
\begin{split}
(1-t)^{2i+1}\, &
\lw{2}F_1\left( \begin{smallmatrix}
i+1,\;i+2-2p \\ 2i +2 \end{smallmatrix}; 1-t \right) \\
&= At^{2p-1}\, \lw{2}F_1\left( \begin{smallmatrix}
-i,\,\;2p-i-1 \\ 2p \end{smallmatrix};   t \right)
+ B\, \lw{2}F_1\left(\begin{smallmatrix}
-i,\,\;-2p-i+1 \\ 2-2p\end{smallmatrix}; t \right)
\end{split}
\end{equation*}

We let  $t\rightarrow 0^+$. Since $\re(2p)>1$ we have that
$\lim_{t\rightarrow 0^+} t^{2p-1}=0$ and that the function $\lw{2}F_1
(i+1,i+2-2p, 2i+2; 1-t)$  converges when $t\rightarrow 0$ to the value
$\lw{2}F_1 (i+1,i+2-2p, 2i+2;
1)=\frac{\Gamma(2i+2)\Gamma(2p-1)}{\Gamma(i+1)\Gamma(i+2p)}$. (See, for
example, \cite{AAR}, Theorem 2.2.2).
Thus
$$ B=\lw{2}F_1\left( \begin{smallmatrix}
i+1,\;i+2-2p \\ 2i +2 \end{smallmatrix}; 1 \right)
=\frac{(i+1)_{i+1}}{(2p-1)_{i+1}}.$$

\noindent  Now by taking $\lim_{t\rightarrow 1}$ we have
$$0=A\, \lw{2}F_1\left( \begin{smallmatrix}
-i,\,\;2p-i-1 \\ 2p \end{smallmatrix}; 1 \right)+ B\,
\lw{2}F_1\left(\begin{smallmatrix}
-i,\,\;-i+1-2p \\ 2-2p\end{smallmatrix}; 1 \right)
=A \,\tfrac{(i+1)_{i}}{(2p)_{i}} +
B\, \tfrac{(i+1)_{i}}{(2-2p)_{i}}.$$
{}From here we can easily deduce that $A=\frac{(i+1)_{i+1}}{(1-2p)_{i+1}}$.

Now we realize that \eqref{lemhiper} is invariant under
the map $p\mapsto 1-p$,
by Euler's transformation (see \cite{AAR}, Theorem 2.2.5). Thus
\eqref{lemhiper}
also holds when $2p\in\CC-\ZZ$ and $\re(2p)<1$, or when
$2p\in\ZZ$ and $2p<-i+1$. By analytic continuation on $p$ the validity of
\eqref{lemhiper} will be extended to $\re(2p)=1$ but $2p\ne1$. In fact the
left
hand side of \eqref{lemhiper} is an analytic function of $p$ for
$2p\in\CC-\ZZ$
if $0<t<1$. On the other hand it is well known  that
$\frac1{\Gamma(c)}\ \lw{2}F_1(a,b,c;x)$ is an entire function of $a,b,c$ if
$|x|<1$
(see \cite{AAR}, p. 65). Thus the right hand side of \eqref{lemhiper} is an
entire
function of $p$ if $0<t<1$. This completes the proof of the lemma.
\end{proof}

\medskip

\begin{prop}\label{complementoV(p)} Let $2p$ be an integer and let
$\ld=4p(p-1)$. We have
\begin{enumerate}
\item [i)]  If $2p\geq 1$ then $V_\ld=V(p)\oplus W_\ld$.

\item [ii)] If $2p\leq 0$ then $V_\ld=V(1-p)\oplus W_\ld$.

\end{enumerate}
\end{prop}

\begin{proof}
Let us first observe that i) and ii) are equivalent by changing $p$ by
$1-p$. Thus let us assume that  $2p\geq 1$.

Let  $H\in V_\ld$, $\ld=4p(p-1)$ and
$\check H=U^{-1}H $, then by Proposition \ref{Hforma} the $i$th-component
$\check h_i$ of $\check H$ is of the form $\check
h_i=\frac{t^p}{(1-t)^{i+1}}R_i$ where $R_i$ is  a solution of the Gauss'
equation \eqref{ecuacionR}. We also have mentioned that a basis of
solutions of this equation is given by  \eqref{basebuena}, when $0\le
i<2p-1$
and
by \eqref{basemala} if $2p-1\le i\le\ell$.
Therefore for any $H \in V_\ld$  the $i$th-component of $\check H$ is
$$ \check h_i= \al_i \frac {t^{p}}{(1-t)^{i+1}}\, \lw{2}F_1 \left(
\begin{smallmatrix} -i,\,\;2p-i-1 \\
2p \end{smallmatrix}; t \right)
+ \beta_i   \frac {t^{1-p}}{(1-t)^{i+1}}\,\lw{2}F_1\left(
\begin{smallmatrix} -i,\,\;-2p-i+1 \\ 2-2p \end{smallmatrix}; t
\right).$$
whenever , $0\le i<2p-1$ and if $2p-1\le i\le\ell$ then
$$ \check h_i= \al_i \frac {t^{p}}{(1-t)^{i+1}}\, \lw{2}F_1 \left(
\begin{smallmatrix} -i,\,\;2p-i-1 \\
2p \end{smallmatrix}; t \right)
+\beta_i t^{p}\, (1-t)^{i}\, \lw{2}F_1\left( \begin{smallmatrix}
i+1,\,\;2p+i \\ 2i+2 \end{smallmatrix}; 1-t \right),$$
for some $\al_i, \beta_i \in \CC$.

In both cases we have seen in the proof of Proposition \ref{dimV(p)} that
the
functions  $H\in V(p)$ have
$\beta_i=0$, for all $i=0,\dots, \ell$. In particular $\check h_i$
diverges
when $t\rightarrow 1^-$.
Therefore we have that $V(p)\cap W_\ld=0$.

We recall that a function $H \in W_\ld$  corresponds to  functions
$\check h_i$ which have a finite limit when $t\rightarrow 1^-$.

In the case  $2p-1\le i\le\ell$ the second summand in $\check h_i$
\begin{equation}\label{betai}
\beta_it^{p}\, (1-t)^{i}\, \lw{2}F_1\left( \begin{smallmatrix} i+1,\,\;2p+i
\\2i+2 \end{smallmatrix}; 1-t \right)
\end{equation}
is convergent at $t=1$, while the first one is not if $\alpha_i\ne0$.
Therefore $\check h_i$ corresponds
to a  function $H\in W_\ld$ if and only if $\al_i=0$.

When $0\le i<2p-1$, by Lemma \ref{hipergeometrico}, if we
choose  $\al_i=a_i\tfrac{(i+1)_{i+1}}{(1-2p)_{i+1}}$ and $\beta_i=a_i
\tfrac{(i+1)_{i+1}}{(2p-1)_{i+1}}$, ($a_i\in \CC$),  we get
\begin{equation}\label{ai}
\check h_i(t)=a_i t^{1-p}(1-t)^{i}\,
\lw{2}F_1\left( \begin{smallmatrix}
i+1,\;i+2-2p \\ 2i +2 \end{smallmatrix}; 1-t \right)
\end{equation}
which is convergent at $t=1$.

Therefore $\dim W_\ld=\ell+1$ and since $\dim V(p)=\ell+1$ (Proposition
\ref{dimV(p)}) the proof is completed.
\end{proof}

\

Recall that for $2p\not\in\ZZ$ we had $V_\lambda=V(p)\oplus V(1-p)$.
When $2p\in\ZZ$ this decomposition fails but can be replaced by one
of the two in Proposition \ref{complementoV(p)}. Returning to the case
$2p\not\in\ZZ$ we have the following proposition.

\begin{prop}\label{complementoV(p)2}
  If $2p\in \CC$ is not an integer and  $\ld=4p(p-1)$ then
$$V_\ld=V(p)\oplus W_\ld=V(1-p)\oplus W_\ld.$$
\end{prop}

\begin{proof}
Let  $H\in V_\ld$, $\ld=4p(p-1)$ and
$\check H=U^{-1}H $, then by Proposition \ref{Hforma} the $i$th-component
$\check h_i$ of $\check H$ is of the form $\check
h_i=\frac{t^p}{(1-t)^{i+1}}R_i$ where $R_i$ is  a solution of the Gauss'
equation \eqref{ecuacionR}.
Since $2p\not\in \ZZ$ a basis of the solutions of this equation is
given by \eqref{basegeneral}.
Therefore for any $H \in V_\ld$  the $i$th-component of $\check H$ is
$$ \check h_i= \al_i \frac {t^{p}}{(1-t)^{i+1}}\, \lw{2}F_1 \left(
\begin{smallmatrix} -i,\,\;2p-i-1 \\
2p \end{smallmatrix}; t \right)
+ \beta_i   \frac {t^{1-p}}{(1-t)^{i+1}}\,\lw{2}F_1\left(
\begin{smallmatrix} -i,\,\;-2p-i+1 \\ 2-2p \end{smallmatrix}; t
\right).$$
Now the proof follows in the same way as in Proposition
\ref{complementoV(p)}.
\end{proof}

\smallskip
\begin{cor} \label{hiperg1-t}
  Let $H\in W_\ld$, $\ld=4p(p-1)$ with $\re(2p)\geq 1$
and let $\check H=U^{-1}H $. Then the $i$th-component  $\check h_i$ of
$\check H$ is of the following form:
\begin{enumerate}
\item [i)] If $2p\not\in \ZZ$ or $2p$ is an integer and $0\le i<2p-1$ there
exists $a_i\in \CC$ such that
$$\check h_i =a_i\, t^{1-p}(1-t)^{i}\,
\lw{2}F_1\left( \begin{smallmatrix}
i+1,\;i+2-2p \\ 2i +2 \end{smallmatrix}; 1-t \right).$$

\item [ii)] If $2p$ is an integer and $2p-1\leq i\le\ell$ then
$$\check h_i =a_i\, t^{p}\, (1-t)^{i}\, \lw{2}F_1\left(
\begin{smallmatrix}
i+1,\,\;2p+i \\ 2i+2 \end{smallmatrix}; 1-t \right).$$
\end{enumerate}
for some constant $a_i\in \CC$.
\end{cor}

\begin{proof} Instance i) consists of two cases. The first one follows
from Proposition \ref{complementoV(p)2} combined with Lemma
\ref{hipergeometrico}.
For the second see \eqref{ai}. Instance ii) was obtained in \eqref{betai}.
\end{proof}

\begin{cor}\label{dimWld} For all $p\in \CC$  we have that  $\dim
W_\ld=\ell+1$.
\end{cor}

\medskip
In the next section  the asymptotic behavior of $H\in V_\ld$ when
$t\rightarrow 0^+$ will play a fundamental role.

\begin{prop}\label{limiteV(p)}
  If $H\in V(p)$, $p\in \CC$,  then $\lim_{t\rightarrow 0^+} t^{-p}H(t)$
exists and is finite.
\end{prop}

\begin{proof} The assertion follows directly from the definition  of $
V(p)$ in \eqref{V(p)def}.
\end{proof}

\begin{prop} \label{limitet=0}
Let $H\in W_\ld$, $\ld=4p(p-1)$ with $\re(2p)> 1$. Then
$\lim_{t\rightarrow 0^+} t^{p-1}H(t)$ exists and is finite.
\end{prop}

\begin{proof} Observe that if  $\check H=U^{-1}H$ then it is equivalent
to
prove the existence of  $\lim_{t\rightarrow 0^+} t^{p-1}\check H(t)$.

In the case $2p\not\in \ZZ$ or $2p$ an integer with $2p > i+1$ we have,
by Corollary \ref{hiperg1-t}, i),  that
$$\check h_i =a_i\, t^{1-p}(1-t)^{i}\,
\lw{2}F_1\left( \begin{smallmatrix}
i+1,\;i+2-2p \\ 2i +2 \end{smallmatrix}; 1-t \right).$$
Since $\re(2p)>1$ Gauss' summation formula gives
$$\lim_{t\rightarrow 0^+} t^{p-1}\check h_i(t)=
\lim_{t\rightarrow 0^+} a_i \, \lw{2}F_1\left( \begin{smallmatrix}
i+1,\;i+2-2p \\ 2i +2 \end{smallmatrix}; 1 \right)= a_i
\tfrac{(i+1)_{i+1}}{(2p-1)_{i+1}}.$$

In the case $2p$ an integer $1\leq 2p\leq i+1$ we have
$$\check h_i =a_i\, t^{p}(1-t)^{i}\,
\lw{2}F_1\left( \begin{smallmatrix} i+1,\;2p+i \\ 2i +2
\end{smallmatrix};
1-t \right).$$
Therefore
$$\lim_{t\rightarrow 0^+} t^{p-1}\check h_i(t)=\lim_{t\rightarrow 0^+}a_i
t^{2p-1}\lw{2}F_1\left(
\begin{smallmatrix} i+1,\;i+2-2p \\ 2i +2 \end{smallmatrix}; 1-t \right)=
a_i\tfrac{(i+1)_{i+1}}{(2p-1)_{i+1}}$$
(see for example \cite{AAR}, Theorem 2.1.3). This completes the proof of the
proposition.
\end{proof}

\begin{prop} \label{limitet=02}
Let $2p=1$ and $H\in W_\ld$, $(\ld=-1)$. Then
$$\lim_{t\rightarrow 0^+} \;\frac{1}{t^{1/2}\log t} H(t)$$
exists and is finite.
\end{prop}
\begin{proof}
It is equivalent to prove that $\lim_{t\rightarrow 0^+}
\;\frac{1}{t^{1/2}\log t}\check H(t)$  exists, if $\check H=U^{-1}H $.
  By Corollary \ref{hiperg1-t} ii), for all $0\leq i\leq \ell$  we
have that
$$\check h_i =a_i\, t^{1/2}(1-t)^{i}\,
\lw{2}F_1\left( \begin{smallmatrix} i+1,\;i+1 \\ 2i +2
\end{smallmatrix};
1-t \right).$$
It is known that (Theorem 2.1.3 in \cite{AAR}):
$$\lim_{t\rightarrow 0^+} \frac{\lw{2}F_1\left( \begin{smallmatrix}
i+1,\;i+2-2p \\ 2i +2 \end{smallmatrix}; 1-t \right)} {\log(1/t)}=
\tfrac{\Gamma(2i+2)}{\Gamma(i+1)\Gamma(i+1)}.$$
Therefore
\begin{equation}\label{valorlimite}
\lim_{t\rightarrow 0^+} \;\frac{1}{t^{1/2}\log t}\check h_i(t)=-a_i
\,\tfrac{\Gamma(2i+2)}{\Gamma(i+1)\Gamma(i+1)}.
\end{equation}
This completes the proof of the proposition.
\end{proof}

\section{Eigenfunctions of $D$ and $E$}\label{DE}

Here we exploit the decompositions of $V_\lambda$, $\lambda=4p(p-1)$,
established in Section \ref{D}:

\noindent i) if $2p\not\in\ZZ$, $V_\lambda=V(p)\oplus V(1-p)=V(p)\oplus
W_\lambda
=V(1-p)\oplus W_\lambda$,

\noindent ii) if $2p\in\ZZ$, $2p\ge1$, $V_\lambda=V(p)\oplus W_\lambda$,

\noindent iii) if $2p\in\ZZ$, $2p\le0$, $V_\lambda=V(1-p)\oplus W_\lambda$.

\begin{prop}\label{Eestables} For all $p\in \CC$ the linear spaces $V_\ld$,
$W_\ld$ and
$V(p)$ are stable under the differential operator $E$. Therefore $E$
restricts to a linear map of each space into itself.
\end{prop}

\begin{proof} Since the differential operators $D$ and $E$ commute $E$
preserves the space  $V_\ld$.

If $H\in V(p)$ then $H=t^p F$, with $F$ analytic at $t=0$. Therefore
$EH\in V(p)$ since
\begin{align}
EH&=-4\,t\, A_0 H'+\frac{4}{1-t}C_0 H-\frac{4\,t}{1-t}C_1 H \nonumber
\\
& = t^p\, \left( -4\,t\, A_0 F'-4\,p\,A_0 F +\frac{4}{1-t}C_0
F-\frac{4\,t}{1-t}C_1 F \right)\label{aux1}.
\end{align}

In considering $W_\lambda$ we can assume that $\re(2p)\ge1$.
{}From Corollary \ref{hiperg1-t} we get that $H\in W_\ld$ if and only if
$H\in V_\ld$ and $H$ is analytic at $t=1$. Thus to prove that $EH\in
W_\ld$
if
$H\in W_\ld$ we only  have to see that
$\frac{1}{1-t}\left(C_0 -C_1\right) H $ is analytic at $t=1$, since we
have
$$\frac{4}{1-t}C_0 H-\frac{4\,t}{1-t}C_1 H= 4C_0+\frac{4\,
t}{1-t}(C_0-C_1)H.$$
Now from $DH=\ld H$ and $H\in W_\lambda$ it follows that
$(C_0+C_1)H(1)=0$.
{}From Proposition  \ref{Hahn} we get that an eigenvector of $C_0+C_1$
of
eigenvalue zero is a scalar multiple of  $u_0=(1,\dots,1 )$. Then
it is easy to check that $C_0u_0= C_1 u_0=0$ and therefore
$ \frac{1}{1-t}(C_0-C_1) H$ is analytic at $t=1$.
\end{proof}

\

\rmk{Remark} The argument above shows that any eigenfunction $H$ of $D$
which has a finite limit $L$ as $t\rightarrow 1^-$ has $L=c(1,\dots,1)$
for some scalar $c$.

\begin{prop}\label{etadef1}
(i) Let $p\in \CC$ , $2p\not \in \ZZ$  and let $\eta:V(p)\longrightarrow
\CC^{\ell+1}$ be the map defined by
$$\eta:H\longmapsto \lim_{t\rightarrow 0^+}(t^{-p}\,H(t)).$$
Then $\eta$ is a linear isomorphism. Moreover the following is a
commutative diagram
$$\begin{CD}
V(p) @ >E >>V(p)\\ @ V \eta VV @ VV \eta V \\ \CC^{\ell+1} @> L(p) >>
\CC^{\ell+1}
\end{CD}$$
where $L(p)$ is the $(\ell+1)\times(\ell+1)$ matrix given by
\begin{equation}\label{Lp}
L(p)=4C_0-4\,pA_0.
\end{equation}
\noindent (ii) If $2p\in\ZZ$, $2p\ge1$, the same result holds.

\noindent (iii) If $2p\in\ZZ$, $2p\le0$, the result holds by interchanging
$p$ into $1-p$.
\end{prop}

\begin{proof}
{}From Proposition \ref{limiteV(p)} it follows that $\eta$ is  well
defined.
Moreover it is clear that $\eta$ is a linear map between two vector
spaces
of dimension $\ell+1$. 
A function $H=t^p\, F(t)\in V(p)$, with $F(t)=\sum_{j=0}^\infty t^j
F_j$
is given by solving the three term recursion relation \eqref{recur1}:
\begin{equation*}
\begin{split}
&j(j-1+2p)F_j-\bigl( 2p(2j+1)+2(j-1)^2- C_0-C_1)\bigr) \,F_{j-1}\\
&\quad +(j-1)(2p+j-2)\,F_{j-2}=0.
\end{split}
\end{equation*}
This recurrence relation  shows that the coefficient vector $F_0$
determines $H(t)$.  Therefore $\eta$ is an injective map. The
nonvanishing of the factor $j(j-1+2p)$ is handled differently in cases
(i), (ii) and (iii).

\noindent To prove the second assertion let $H(t)=t^p F\in V(p)$ with
$F(t)=\sum_{j=0}^\infty\,t^jF_j$. Then by \eqref{aux1} we get
$$\eta(EH)=(t^{-p}(EH))(0)=-4pA_0F_0+4 C_0F_0=L(p)\,F_0=L(p)(\eta (H)).$$
\end{proof}

In Proposition \ref{Eautovalores} we shall prove that $L(p)$ is a
diagonalizable matrix with eigenvalues
$$\mu_k=\mu_k(p)=-4p(\ell-2k)-4k(\ell-k+1), \quad \text{for } 0\le k\le
\ell,$$
all with multiplicity one.
  Note that the sets of eigenvalues of $L(p)$ and $L(1-p)$ are the same, in
fact we have
$$\mu_k(1-p)=\mu_{\ell-k}(p)\, , \quad \text{ for all }\,  1\leq k \leq
\ell.$$

\

The following theorem gives us all simultaneous solutions of the system
$DH=\ld H$ and $EH=\mu H$, thus accomplishing task (i) at the beginning of
Section \ref{D} for $2p\not\in\ZZ$.

\begin{thm}\label{mainnoentero}
Let $\ld=4p(p-1)$, $2p\not\in \ZZ$. A function $H=H(t)$, $0<t<1$ is a
simultaneous eigenfunction of the system $DH=\ld H$ and $EH=\mu H$ if and
only if $H$ is of the following form:
  $$H(t)=\frac{t^p}{(1-t)^{\ell+1}} P(t)+ \frac{t^{1-p}}{(1-t)^{\ell+1}}
Q(t)$$ where $P(t)=U\check P(t)$ and  $Q(t)=U\check Q(t)$ and the
$i$th-component of $\check P$ and $\check Q$ are
\begin{align*}
\check P_i(t)&= \al_i(1-t)^{\ell-i} \lw{2}F_1\left( \begin{smallmatrix}
-i,\;2p-i-1 \\ 2p \end{smallmatrix}; t \right),\displaybreak[0]\\
\check Q_i(t)&= \beta_i(1-t)^{\ell-i} \lw{2}F_1\left( \begin{smallmatrix}
-i,\;-2p-i+1 \\ 2-2p \end{smallmatrix}; t \right)
\end{align*}
  with  $P(0)$ and  $Q(0)$, respectively,  eigenvectors of $L(p)$ and
$L(1-p)$,  of eigenvalue $\mu$.
\end{thm}

\begin{proof}
{}From Propositions \ref{2pnoentero} and \ref{Vlambda} it follows that
$H\in V_\lambda$ if and only if it is of the form stated in the first part
of the theorem.

Since $V(p)$ and $V(1-p)$ are $E$-stable and $V_\ld=V(p)\oplus V(1-p)$ it
follows that $H\in V_\ld $ is an eigenfunction of $E$
if and only if each
summand of $H$ is an eigenfunction of $E$ with the same eigenvalue, which is
equivalent by Proposition
\ref{etadef1} to the fact that $P(0)$ and $Q(0)$ be, respectively, an
eigenvector of $L(p)$ and $L(1-p)$ for the same eigenvalue.
\end{proof}

\

We will see in Section \ref{ejemplos} that even when $2p\in\ZZ$, $H(t)$
allows
a decomposition in the form above except for a limited range of values of
$p$. Under the condition $2p\notin\ZZ$ this decomposition is unique in the
sense of Proposition \ref{Vlambda}. This uniqueness fails when $2p\in\ZZ$.

\

To study the joint eigenfunctions $H=H(t)$ of the differential operators
$D$ and $E$ when $2p\in \ZZ$ we can use the decomposition of $V_\ld$ given
in Proposition \ref{complementoV(p)}.

The following three propositions describe the action of $E$ in $W_\lambda$
for any $p\in\CC$.

\medskip

\begin{prop}\label{etadiagr}
Let $p\in \CC$, $\re(2p)>1$ and let $\eta:W_\ld\longrightarrow
\CC^{\ell+1}$, $\ld=4p(p-1)$, be the map defined by
$$\eta:H\longmapsto \lim_{t\rightarrow 0^+}t^{p-1}\,H(t).$$
Then $\eta$ is a linear isomorphism. Moreover the following is a
commutative diagram
$$\begin{CD}
W_\ld @ >E >>W_\ld \\ @ V \eta VV @ VV \eta V \\ \CC^{\ell+1} @> L(1-p)
>>\CC^{\ell+1}
\end{CD}$$
where $L(1-p)$ is the $(\ell+1)\times(\ell+1)$ matrix given by
$$L(1-p)=4C_0-4\,(1-p)A_0.$$
A similar statement holds for $\re(2p)<1$, by changing $p$ into $1-p$.
\end{prop}

\begin{proof} {}From Proposition \ref{limitet=0} it follows that $\eta$
is well defined. Moreover it is clear that $\eta$ is a linear map between
two vector spaces of dimension $\ell+1$. So, we only have to verify that
$\eta$
is injective.

Let $H\in W_\ld$ and let $\check H=U^{-1}H =(\check h_0, \dots , \check
h_\ell)$. In the proof of Proposition \ref{limitet=0}
we have seen that in the case $2p\not\in \ZZ$ or $2p$ an integer $0\le
i<2p-1$
$$\check h_i =a_i\, t^{1-p}(1-t)^{i}\,
\lw{2}F_1\left( \begin{smallmatrix} i+1,\;i+2-2p \\ 2i +2
\end{smallmatrix}; 1-t \right),$$
and in the case $2p$ an integer $2p-1\leq i\leq\ell$ we have
$$\check h_i =a_i\, t^{p}\, (1-t)^{i}\, \lw{2}F_1\left(
\begin{smallmatrix}
i+1,\,\;2p+i \\ 2i+2 \end{smallmatrix}; 1-t \right).$$
We also  proved that in both cases
  $$\lim_{t\rightarrow 0^+} t^{p-1}\check h_i(t)=  a_i
\tfrac{(i+1)_{i+1}}{(2p-1)_{i+1}}.$$
Now  if $H\in W_\ld$ satisfies $\eta(H)=0$ then $\lim_{t\rightarrow
0^+}t^{p-1}\,\check h_i(t)=0$.
So we have  $a_i=0$ for all $0\leq i \leq \ell$ and thus $H=0$.
Therefore $\eta$ is an injective map.

To prove the second assertion let $H\in W_\ld$. {}From \eqref{EE} we have
$$\eta(EH)=\lim_{t\rightarrow 0^+}\left(-4A_0 t^{p}H'
+\tfrac{4}{1-t}C_0
t^{p-1}H-\tfrac{4}{1-t}C_1 t^pH\right).$$
We first note that $\lim_{t\rightarrow 0^+}t^pH(t)=0$ because
$\lim_{t\rightarrow 0^+}t^{p-1}H(t)$ exists and is finite.
By L'Hospital rule we have
$$\eta(H)=\lim_{t\rightarrow 0^+} \, \frac{t^pH(t)}{t}=p\eta(H)+
\lim_{t\rightarrow 0^+}\, t^pH'(t)$$
thus $\lim_{t\rightarrow 0^+}\, t^pH'(t)=(1-p)\eta(H)$. Therefore
$$\eta(EH) = -4(1-p) A_0 \eta(H)+4C_0\eta(H)= L(1-p)\, \eta(H).$$
\end{proof}

\medskip

\begin{prop} \label{etadiagr1}
Let $2p=1$, that is $\ld=-1$, and let $\eta:W_\ld\longrightarrow
\CC^{\ell+1}$ be the map defined by
$$\eta:H\longmapsto \lim_{t\rightarrow 0^+}\frac{1}{t^{1/2}\log t }H(t).$$
Then $\eta$ is a linear isomorphism. Moreover the following is a
commutative diagram
$$\begin{CD}
W_\ld @ >E >>W_\ld \\ @ V \eta VV @ VV \eta V \\ \CC^{\ell+1} @> L(1/2)
>>\CC^{\ell+1}
\end{CD}$$
\end{prop}

\begin{proof}
By Proposition \ref{limitet=02} $\eta$ is a well defined linear map between
two vector spaces of dimension $\ell+1$ (Corollary \ref{dimWld}).
Let $H\in W_\ld$ be such that $\eta(H)=0$. {}From \eqref{valorlimite} we get
$a_i=0$ for all $0\leq i \leq \ell$. Thus we have that $H=0$, proving that
$\eta$ is an isomorphism.

To prove that the diagram is commutative take $H\in W_\ld$. By definition
of the differential operator $E$ we have
$$\eta(EH)=\lim_{t\rightarrow 0^+} \left( -\frac{4t^{1/2}}{\log t} \,A_0
H'+ \frac{4}{t^{1/2}\log t}\,C_0 H -\frac{4t^{1/2}}{\log t}\, C_1 H
\right).$$

\noindent We observe that $\lim_{t\rightarrow 0^+} \frac{4}{t^{1/2}\log
t}\, H= 4\eta(H)$ and that
$\lim_{t\rightarrow 0^+} \frac{4\,t^{1/2}}{\log t}\, H= \lim_{t\rightarrow
0^+} \frac{4\,t}{t^{1/2}\log t}\, H=0.$
  Finally by L'Hospital rule we have
$$\eta(H)= \lim_{t\rightarrow 0^+}\frac{H(t)}{t^{1/2}\log t }=
\lim_{t\rightarrow 0^+}\frac{H'(t)}{\tfrac 12 t^{-1/2}\log t + t^{-1/2}}=
2\lim_{t\rightarrow 0^+}\frac{t^{1/2} H'(t)}{\log t }.$$

\noindent Therefore $\eta(EH)=-2A_0\eta(H)+ 4C_0\eta(H)= L(1/2)\eta(H)$.
This completes the proof of the proposition.
\end{proof}

\begin{prop}\label{etadiagr2}
Let $p\in \CC$, $\re(2p)=1$, $2p\neq 1$ and let $\eta:W_\ld\longrightarrow
\CC^{\ell+1}$, $\ld=4p(p-1)$, be the map defined by
$$\eta:H\longmapsto \lim_{t\rightarrow 0^+}t^{p-1}\,P(H)(t),$$
where $P$ denotes the projection of $V_\ld$ onto $V(1-p)$ along the
subspace $V(p)$.
Then $\eta$ is a linear isomorphism. Moreover the following is a
commutative diagram
$$\begin{CD}
W_\ld @ >E >>W_\ld \\ @ V \eta VV @ VV \eta V \\ \CC^{\ell+1} @> L(1-p)
>>\CC^{\ell+1}.
\end{CD}$$
\end{prop}

\begin{proof}
Let us observe that $\eta=\eta_1P$ where $\eta_1:V(1-p)\longrightarrow
\CC^{\ell+1}$ is the isomorphism introduced in Proposition \ref{etadef1},
exchanging $p$ by $1-p$. Let $H\in W_\ld$ such that $\eta(H)=0$. Then
$P(H)=0$, i.e. $H\in V(p)$, but $V(p)\cap W_\lambda=0$. Hence $H=0$ proving
that
$\eta$ is an isomorphism since $\dim W_\ld=\ell+1$ (see Corollary
\ref{dimWld}).

To prove the commutativity of the  diagram we observe that given $H\in
W_\ld$ we have
$$\eta(EH)=\eta_1 P(EH)=\eta_1 EP(H)=L(1-p)\eta_1(P(H))=L(1-p)\eta(H),$$
by Propositions \ref{Eestables} and \ref{etadef1}. This completes the proof
of the Proposition.
\end{proof}

\

\rmk{Remark} The construction above, with the introduction of the
projection $P$ may appear different from that in Proposition \ref{etadiagr}
and \ref{etadiagr1}. We explain now the need and the motivation behind
this choice of $P$.

In this case $p=\frac12+i\sigma$, $\sigma\ne0$, real. Proposition
\ref{Vlambda}
shows that any $H\ne0$ in $W_\lambda$ has two components: one behaves like
$t^{\frac12+i\sigma}$ and the other one as $t^{\frac12-i\sigma}$ as
$t\rightarrow 0^+$. There is no way to multiply both components by a common
function of $t$ to describe the asymptotic behaviour of $H$ at $t=0$. The
only
way out is to eliminate one of the two components first. This is the effect
of the projection. In the cases considered in Propositions \ref{etadiagr}
and \ref{etadiagr1} this was not necessary, but we could have done it too.

This phenomenon is reminiscent of the presence of ``Stokes lines¥¥ in the
discussion of ``dominant and recessive solutions¥¥ of an ordinary
differential
equation at an irregular singular point.

\

Now is time to take the study of the matrix $L(p)$. In the case $\ell=0$
we have $L=E=0$ for all $p$, thus we assume below $\ell>0$.

\begin{prop}\label{Eautovalores} If $2p\in\CC$ is not an
integer in the interval $-\ell+2\le 2p\le\ell$ then
$L(p)=4C_0-4pA_0$ is diagonalizable with eigenvalues
\begin{equation}\label{muk}
\mu_k=-4p(\ell-2k)-4k(\ell-k+1), \quad \text{for } 0\le k\le \ell,
\end{equation}
all with multiplicity one.

If $2p$ is an integer in the interval $-\ell+2\le 2p\le\ell$ then
$L(p)=4C_0-4pA_0$ is not diagonalizable. More precisely:

\begin{enumerate}
\item[i)] If $1\le 2p\le\ell$ then $\mu_k$, $0\le k<\frac{\ell+1-2p}2$,
are eigenvalues with multiplicity two,
but with geometric multiplicity one; $\mu_k$, $\ell+1-2p<k\le\ell$ or
$k=\frac{\ell+1-2p}2$ (when $\ell+1-2p$ is even), are
eigenvalues with  multiplicity one.

\item [ii)] If $-\ell+2\le 2p\le 0$ then $\mu_k$, $1-2p\le
k<\frac{\ell+1-2p}2$, are eigenvalues with multiplicity two,
but with geometric multiplicity one; $\mu_k$, $0\le k<1-2p$ or
$k=\frac{\ell+1-2p}2$ (when $\ell+1-2p$ is even), are
eigenvalues with  multiplicity one.
\end{enumerate}

The eigenspace corresponding
to the eigenvalue $\mu_k$ is generated by
$v_k=(v_{0,k},\dots,v_{\ell,k})$,
where $v_{i,k}=0$, for $0\leq i<k$, $v_{k,k}=1$
and
\begin{equation}\label{autovector}
v_{k+j,k}=\prod_{i=1}^{j} \frac{(k+i)(\ell+1-i-k)}{i(\ell+1-2p-2k-i)},
\end{equation}
for  $1\leq j \leq \ell-k$. If $\mu_k$ is an eigenvalue with
multiplicity
two and $\mu_k=\mu_{k'}$ we assume that $k>k'$.
\end{prop}

\begin{proof}
The matrix $L(p)=4C_0-4pA_0$ is a lower triangular matrix given by
$$L(p)=\sum_{i=0}^\ell
4i(\ell-i+1)E_{i,i-1}-4(i(\ell-i+1)+p(\ell-2i))E_{i,i}.$$
Therefore the eigenvalues of $L(p)$ are $\mu_k= -4p(\ell-2k)-
4k(\ell-k+1)$,
for each $k=0,\dots ,\ell$.

As a function of $k$, $\mu_k$ is symmetric around the point
$\frac{\ell+1-2p}2$, i.e. $\mu_k=\mu_{\ell+1-2p-k}$.
Therefore there are eigenvalues with multiplicity two if and only if
$2p\in\ZZ$ and $0<\frac{\ell+1-2p}2<\ell$, or
equivalently $2p\in\ZZ$ and $-\ell+2\le 2p\le\ell$. This proves the
first
assertion.

To prove i) it is enough to visualize in the real line the points
$0<\frac{\ell+1-2p}2<\ell+1-2p<\ell$ and exploit their relative
ordering
and symmetry around $\frac{\ell+1-2p}2$. Similarly to establish
ii) it is enough to visualize in the real line the points
$0<1-2p<\frac{\ell+1-2p}2<\ell$.

To prove the last assertion we first observe that
$x=(x_0,\dots,x_\ell)$ is
an eigenvector of $L(p)$ corresponding to the
eigenvalue $\mu_k$ if and only if the following equations hold
\begin{equation}\label{auto}
4i(\ell+1-i)x_{i-1}=(\mu_k-\mu_i)x_i\quad\quad{\text {for
all }}\quad
1\le i\le\ell.
\end{equation}
If $\mu_k$ is an eigenvalue with multiplicity two and $\mu_k=\mu_{k'}$
we
assume that $k>k'$. Then from (\ref{auto}) it
follows right away that $x_i=0$ for $0\le i\le k-1$, and that $x_k$
determines $x_i$ for $k+1\le i\le\ell$. Moreover a
closer look at (\ref{auto}) shows that $v_k$, see (\ref{autovector}),
generates the eigenspace associated to the eigenvalue
$\mu_k$.
\end{proof}

\

By collecting our previous results we arrive at our punch line:
for each $\ell\in\ZZ_{\ge0}$ and $p\in\CC$ there are exactly $\ell+1$
functions $H_k(t,p)$, $0\le\k\le\ell$, associated to spherical functions.
Their explicit expressions are given in the next proposition.
By the usual symmetry $p\mapsto 1-p$ we can assume $\re(2p)\ge1$.

\begin{thm}\label{punch}
Let $\ld=4p(p-1)$ with $\re(2p)\geq 1$. A function $H=H(t)$, $0<t<1$,
corresponds to a spherical function $\Phi$ of type $\pi=\pi_\ell$ if and
only if $H=U\check H$ where:
\begin{enumerate}
\item [i)] If $2p\not\in \ZZ$
$$\check h_i =a_i\, t^{1-p}(1-t)^{i}\,
\lw{2}F_1\left( \begin{smallmatrix}
i+1,\;i+2-2p \\ 2i +2 \end{smallmatrix}; 1-t \right)\quad 0\leq i\leq
\ell.$$

\item [ii)] If $2p\in \ZZ $
$$\check h_i= \begin{cases}
a_i\, t^{1-p}(1-t)^{i}\,
\lw{2}F_1\left( \begin{smallmatrix}
i+1,\;i+2-2p \\ 2i +2 \end{smallmatrix}; 1-t \right)
  & \quad 0\leq i \leq 2p-2, \\
a_i\, t^{p}\, (1-t)^{i}\, \lw{2}F_1\left(
\begin{smallmatrix}
i+1,\,\;2p+i \\ 2i+2 \end{smallmatrix}; 1-t \right) & \quad 2p-1\leq i\leq
\ell.
\end{cases}
$$
\end{enumerate}

\noindent When $2p\notin\ZZ$ or $i=0,\dots,2p-2$,
$a_i=\al_i\frac{(2p-1)_{i+1}}{(i+1)_{i+1}}$ and when $i=2p-1,\dots,\ell$,
$a_i=\al_i\frac{(1-2p)_{i+1}}{(i+1)_{i+1}}$. In both cases
$\al_0,\dots,\al_\ell$ is such
that $U(\al_0, \dots, \al_\ell)^T$ is an eigenvector of the matrix $L(1-p)$
of eigenvalue $\mu_k$, $0\le k\le\ell$, and $a_0$=1.

\noindent Moreover $\Omega \Phi=(\frac\lambda4+\mu_k+\ell(\ell+2)) \Phi$ and
$\overline \Omega\Phi=\frac\lambda4\Phi$.
\end{thm}

The functions $H_k(t,p)$ resulting from the $\ell+1$ choices of $\mu_k$
above
will be the main characters in Sections \ref{ejemplos} and
\ref{biespectral}.
We refer to these functions as {\em the families} associated to a given
$\ell$.

In the following proposition we clarify the correspondence
$(p,k)\mapsto\Phi_{(p,k)}$
which assigns to the pair $(p,k)\in\CC\times\ZZ$, $0\le k\le\ell$, the
spherical function
$\Phi_{(p,k)}$ of type $\ell$ associated to the eigenvalues
$\mu_k=-4p(\ell-2k)-4k(\ell-k+1)$ and $\lambda=4p(p-1)$.

\begin{prop}\label{fibra} If  $\Phi$ is a spherical function on $G$ of type
$\ell$ then
$\Phi=\Phi_{(p,k)}$ for some pair $(p,k)$. Moreover
$\Phi_{(p,k)}=\Phi_{(p',k')}$
if and only if
\begin{equation}
(p',k')=
\begin{cases}
(p,k),\\
(1-p,\ell-k),\\
(p,\ell+1-k-2p)\quad &\text{if}\;\; \ell+1-k-2p\in\{0,\dots,\ell\},\\
(1-p,k-1+2p)\quad &\text{if}\;\; k-1+2p\in\{0,\dots,\ell\}.
\end{cases}
\end{equation}
\end{prop}
\begin{proof} Let $\lambda(p,k)=4p(p-1)$ and
$\mu(p,k)=-4p(\ell-2k)-4k(\ell-k+1)$. By
Proposition \ref{unicidad} below the problem reduces to finding the pairs
$(p',k')$ such that
$\lambda(p,k)=\lambda(p',k')$ and $\mu(p,k)=\mu(p',k')$.

Now $\lambda(p,k)=\lambda(p',k')$ if and only if $(p-p')(p+p'-1)=0$, that
is, if and only
if $p'=p$ or $p'=1-p$.

Also it is easy to check that $\mu(p,k)=\mu(p,k')$ if and only if $k=k'$ or
$k'=\ell+1-k-2p$, and that $\mu(p,k)=\mu(1-p,\ell-k)$.

Therefore for $(p',k')$ we have the following possibilities: $(p,k)$ and
$(p,\ell+1-k-2p)$,
or $(1-p,\ell-k)$ and $(1-p,\ell+1-(\ell-k)-2(1-p))=(1-p,k-1+2p)$. This
completes the proof
of the proposition.
\end{proof}

\medskip

Throughout this paper we have carried out a bare hands construction of the
spherical functions starting from their definition. The rest of this section
relates
such a construction to the representation theory of $G$.

Spherical functions of type $\delta\in\hat K$ arise in a natural way upon
considering representations of $G$.
If $g\mapsto U(g)$ is a continuous representation of $G$ , say on a
complete, locally convex, Hausdorff topological vector
space $E$, then
$$P(\delta)=\int_K \chi_\delta(k^{-1})U(k)\, dk$$
is a continuous projection of $E$ onto
$P(\delta)E=E(\delta)$; $E(\delta)$ consists of those vectors in $E$,
the linear span of whose $K$-orbit is finite dimensional
and splits into irreducible $K$-subrepresentations of type $\delta$.
Whenever $E(\delta)$ is finite dimensional and not zero, the function
$\Ph_\delta:G\longrightarrow \End(E(\delta))$ defined by
$\Ph_\delta(g)a=P(\delta)U(g)a$,
$g\in G, a\in E(\delta)$ is a spherical function
of type $\delta$. Moreover any irreducible spherical function arises in this
way from a topologically irreducible representation of $G$,
(see \cite{GO}, \cite{GV} and \cite{T}).

If a spherical function $\Phi$ is associated to a Banach representation
of $G$ then it is quasi-bounded,
in the sense that there exists a semi-norm $\rho$ on $G$ and $M\in\RR$ such
that $\|\Phi(g)\|\le M\rho(g)$
for all $g\in G$. Conversely, if $\Phi$ is an irreducible quasi-bounded
spherical function on $G$, then
it is associated to a topologically irreducible Banach representation
of $G$ (see \cite{GO}, \cite{GV} and \cite{T}). Thus if
$G$ is compact any irreducible spherical function on $G$ is associated
to a Banach representation of $G$,
which is finite dimensional by the Peter-Weyl theorem.

When $G$ is a connected Lie group any spherical function is analytic, hence
it is
determined by $(D\Phi)(e)$ for all left invariant differential operators $D$
on $G$.
We also need to quote the following fact (cf. \cite{T}, Remark 4.7).

\begin{prop}\label{unicidad} Let $\Phi,\Psi:G\longrightarrow \End(V)$ be
two
spherical functions on a connected
Lie group $G$ of the same type $\delta\in K$. Then $\Phi=\Psi$ if and
only if
$(D\Phi)(e)=(D\Psi)(e)$ for all $D\in D(G)^K$.
\end{prop}

If $G$ is a noncompact connected semisimple Lie group with finite center,
and
$K$  is a maximal compact subgroup of $G$, then, from Casselman's
subrepresentation theorem (see \cite{Kn} p. 238), we know that any
irreducible $(\lieg,K)$
module can be  realized inside a generalized (also known as nonunitary)
principal series
representation of
$G$. Thus, in  particular, any irreducible spherical function of $(G,K)$ is
associated to a
generalized principal series representation of $G$, and thus it is
quasi-bounded.

This can be seen directly for $(G,K)=$(SL$(2,\CC)$,SU$(2)$) from our
explicit
construction of all irreducible spherical functions. Let $MAN$ be the upper
triangular subgroup of $G$. Define $\sigma$ on $M$ and $\nu$ on $\liea$ by
\begin{equation*}
\begin{split}
\sigma\left(\begin{matrix}e^{i\theta}&0 \\ 0&e^{-i\theta}\end{matrix}\right)
&=e^{ir\theta},\\
\nu\left(\begin{matrix}t&0 \\ 0&-t\end{matrix}\right)&=ivt,
\end{split}
\end{equation*}
for $r\in\ZZ$ and $v\in\CC$. We also point out that the half-sum of the
positive restricted roots is given by
\begin{equation*}
\rho\left(\begin{matrix}t&0 \\ 0&-t\end{matrix}\right)=2t.
\end{equation*}
Then $man\mapsto e^{\nu\log a}\sigma(m)$ is a representation of $MAN$, and
it is this
representation that we induce to $G$ to construct its generalized principal
series
representation. Thus we put $U^{\sigma,\nu}=Ind_{MAN}^G$.

A dense subspace of the representation space of $U=U^{\sigma,\nu}$ is
$$\{F:G\longrightarrow \CC \text{\; continuous}: F(xman)=e^{-(\nu+\rho)\log
a}\sigma(m)^{-1}F(x)\}$$ with norm
$$||F||^2=\int_K|F(k)|^2\, dk,$$
and $G$ acts by
$$U^{\sigma,\nu}(g)F(x)=F(g^{-1}x).$$
The actual Hilbert space and representation are then obtain by completion.

A $K$ type $\pi_\ell$ occurs in $U^{\sigma,\nu}$ if and only if $r=\ell-2j$
for some
$0\le j\le\ell$ (Frobenius Reciprocity Theorem). If this is the case, let
$\Phi=\Phi_{v,r}$
be the spherical function (irreducible) of type $\pi_\ell$ associated to
$U^{\sigma,\nu}$, and let $P_\ell$ be the corresponding projection onto the
$K$-isotypic component. To identify $\Phi_{v,r}$ we just need to compute
$[\Omega\Phi_{v,r}](e)$ and
$[\overline\Omega\Phi_{v,r}](e)$, (Proposition \ref{unicidad}). We start by
observing
that
$$4\Omega=H_1^2-H_2^2-4H_1+4iH_2-2iH_1H_2+4X_1Y_1-4X_2Y_2-4iX_1Y_2-4iX_2Y_1,
$$
$$4\overline\Omega=H_1^2-H_2^2-4H_1-4iH_2+2iH_1H_2+4X_1Y_1-4X_2Y_2+4iX_1Y_2+
4iX_2Y_1,$$
where $X_1=\frac12(V_1+W_1)$, $Y_1=\frac12(V_1-W_1)$, $X_2=\frac12(W_2+V_2)$
and
$Y_2=\frac12(W_2-V_2)$.

If $X\in\lien$, $Y\in\lieg$, and $F$ is a smooth function in the induced
space, then
$$\dot U(X)\dot U(Y)F(e)=\frac{d}{dt}\,\dot U(Y)F((\exp
tX)^{-1})\Big|_{t=0}=0,$$
because $\dot U(Y)F$ is right invariant under $N$. On the other hand
if $X=H+T\in \liea+\liem$, then
\begin{equation*}
\begin{split}
\dot U(X)F(e)&=\frac{d}{dt}\,\dot UF((\exp tX)^{-1})\Big
|_{t=0}=\frac{d}{dt}\,(e^{\nu+\rho}
(\exp tH)\sigma(\exp tT))F(e)\Big |_{t=0}\\
         &=\big((\nu+\rho)(H)+\dot\sigma(T)\big)F(e).
\end{split}
\end{equation*}
Since $X_1, X_2\in \lien$, $H_1\in\liea$ and
$H_2\in\liem$ we obtain
\begin{equation*}
\begin{split}
4\dot
U(\Omega)F(e)&=\big((iv+2)^2-(ir)^2-4(iv+2)+4i(ir)-2i(iv+2)(ir)\big)F(e)\\
                       &=(iv+2+r)(iv-2+r)F(e),\\
4\dot U(\overline\Omega)F(e)
&=\big((iv+2)^2-(ir)^2-4(iv+2)-4i(ir)+2i(iv+2)(ir)\big)F(e)\\
                       &=(iv+2-r)(iv-2-r)F(e).
\end{split}
\end{equation*}
Since $\Omega,\overline\Omega\in D(G)^G$, $\dot U(\Omega)$ and $\dot
U(\overline\Omega)$
commute with $U(g)$ for all $g\in G$. Therefore $\dot
U(\Omega)=(iv+2+r)(iv-2+r)I$ and
$\dot U(\overline\Omega)=(iv+2-r)(iv-2-r)I$.

Now $\Omega\Phi(g)P_\ell=P_\ell U(g)\dot
U(\Omega)P_\ell= 4(iv+2+r)(iv-2+r)\Phi(g)P_\ell$.
A similar argument holds for $\overline\Omega$. Thus
$$\Omega\Phi=4(iv+2+r)(iv-2+r)\Phi,\quad \text{and}\quad
\overline\Omega\Phi=4(iv+2-r)(iv-2-r)\Phi.$$

\begin{prop}
For $v\in\CC$ and $r=\ell-2j$, $0\le j\le\ell$, let $p=\frac14(iv+2-r)$ and
$k=\frac12(r+\ell)$. Then the spherical function $\Phi_{v,r}$ of type
$\pi_\ell$
associated to $U^{\sigma,\nu}$ is equivalent to $\Phi_{(p,k)}$.
\end{prop}
\begin{proof}
As we said before, the only thing we have to verify is that the eigenvalues
of $\Omega$ and $\overline\Omega$ corresponding to $\Phi_{v,r}$ and
$\Phi_{(p,k)}$
are respectively the same. Taking into account Remark
\ref{relacionautovalores2},
Theorem \ref{punch} and \eqref{muk} this in turn amounts to checking that
$$(iv+2-r)(iv-2-r)=4p(4p-4),$$
and
$$ivr=-4p(\ell-2k)-4k(\ell-k+1)+\ell(\ell+2),$$
which is straightforward.
\end{proof}

\begin{remark}\label{onto} The resulting relations $v=i(\ell-2k+2-4p)$,
$r=2k-\ell$,
with $p\in\CC$ and $0\le k\le\ell$ show that the map $(v,r)\mapsto(p,k)$ is
one to one and onto. Thus any irreducible spherical function is equivalent
to one
associated to a principal series representation as advertised above.
Moreover the
map $(v,r)\mapsto(-v,-r)$ corresponds to the map $(p,k)\mapsto(1-p,\ell-k)$.
This
explains the equivalence between $\Phi_{(p,k)}$ and $\Phi_{(1-p,\ell-k)}$
(Proposition
\ref{fibra}) in terms of the  equivalence between $U^{\sigma,\nu}$ and
$U^{w\cdot\sigma,w\cdot\nu}$, where $w$ is the nontrivial  element of the
restricted Weyl group
$W(\lieg,\liea)$.
\end{remark}

\section{Some examples and a conjecture}\label{ejemplos}

The results in Sections \ref{D} and \ref{DE} include a number of different
expressions for the components of the vector $H(t)$ corresponding to a
spherical function. Some of these expressions are rather elaborate; for
instance those given in Theorem \ref{punch} involve Gauss'
hypergeometric functions as well as Hahn polynomials. Others involve simple
linear combinations of $t^p$ and $t^{1-p}$ with rational coefficients, as
in Theorem \ref{mainnoentero}. This ''embarrassment of riches" results from
different decompositions of $V_\ld$ as a direct sum of subspaces. These
decompositions are dictated by the nature of $p$ as summarized at the
beginning of Section \ref{DE}. The case $2p\not\in \ZZ$ is simpler to deal
with and this is reflected in the simpler expressions appearing in Theorem
\ref{mainnoentero}.

The relative merit of these or other explicit representations of $H(t)$
depends on the task at hand. In this section we identify those cases when
$H(t)$ contains logarithmic terms. For this task an expression as in
Theorem \ref{mainnoentero} is most convenient and we show that its validity
goes beyond the condition $2p\not\in \ZZ$.

The first part of this section is done entirely in terms of explicit
examples given in full detail. This should allow the reader to check many
points of the theoretical treatment in the rest of the paper.

We close with a conjecture of the form of $H(t)$ which is a ''first cousin"
of a conjecture in \cite{GPT}.

\medskip

\noindent
{\bf The case $\ell=0$.}
In the case $\ell=0$ there is only one equation, namely
\[
4t^2 h''(t) - \frac {8t^2}{1-t} h'(t) = 4p(p-1)h(t)
\]
with a solution
\[
h(t)=\frac {t^p - t^{1-p}}{(2p-1)(t-1)}
\]
as in Theorem \ref{mainnoentero}.
This expression is valid as soon as $p$ is not $\frac12$. In the case
$p=\frac12$
the spherical function is given by taking the limit in the expression above.
One gets

\[
h(t)=\frac {\sqrt{t} \log (t)}{t-1}.
\]

This, as well as the examples below, illustrates the fact that the result
in Theorem \ref{mainnoentero} is valid as long as $2p$ is not an integer in
the range $-\ell+1,-\ell+2,\dots,-1,0,1,\dots,\ell+1$. In fact it is only
for these values of $p$ that  spherical functions will contain a
logarithmic term.

We display this in the cases $\ell=1$ and $\ell=2$ below. Afterwards we
describe the pattern in the case of $\ell=5$.

\medskip
\noindent {\bf The case $\ell=1$.}
We have one family with $\mu = -4 p$. (corresponding to $k=0$).
\[
\begin{split}
\ h_0(t) &= \frac {t^{2-p}}{(p-1)(2p-1)(t-1)^2} + \frac
{t^p(2(p-1)t-2p+1)}{(p-1)(2p-1)(t-1)^2} ,  \\
h_1(t) &= \frac {t^p}{(p-1)(2p-1)(t-1)^2} - \frac
{t^{1-p}((2p-1)t-2p+2)}{(p-1)(2p-1)(t-1)^2}.
\end{split}
\]
The other family (corresponding to $k=1$) is obtained by exchanging $p$ and
$1-p$.

We have an exceptional form for $p=\frac12$, namely
\[
h_0(t) = \frac {2\sqrt{t} (t \log (t) - t+1)}{(t-1)^2}, \qquad
h_1(t) = -\frac {2\sqrt{t} (\log (t) - t + 1)}{(t-1)^2}
\]
and for $p=1$ where we get
\[
h_0(t) = -\frac {2t(\log (t) - t + 1)}{(t-1)^2},\qquad
h_1(t) = \frac {2(t \log (t) - t + 1)}{(t-1)^2} .
\]
As usual (for reasons of symmetry) it is enough to consider $p$ to
one side of the
value $\frac12$.

In summary we have, for $\ell=1$, a total of two spherical functions with
logarithmic terms.

\

\noindent {\bf The case $\ell=2$.}
We have three families corresponding to $k=0,1,2$.

\

The first one has $\mu=-8p$
and the vectors $P(t)$ and $Q(t)$ in Theorem \ref{mainnoentero} can be taken
to be

\begin{align*}
P_0(t)&= \frac
{3(2p^2t^2-5pt^2+3t^2-4p^2t+8pt-3t+2p^2-3p+1)}{(p-1)(2p-3)(2p-1)},
\displaybreak[0] \\
P_1(t)&=\frac {3(2pt-3t-2p+1)}{(p-1)(2p-3)(2p-1)},\displaybreak[0] \\
P_2(t)&=\frac {3}{(p-1)(2p-3)(2p-1)},
\end{align*}

and
\begin{align*}
Q_0(t)&= -\frac {3t^2}{(p-1)(2p-3)(2p-1)}, \displaybreak[0] \\
Q_1(t)&=
\frac {3t(2pt-t-2p+3)}{(p-1)(2p-3)(2p-1)},\displaybreak[0] \\
Q_2(t)&= -\frac
{3(2p^2t^2-3pt^2+t^2-4p^2t+8pt-3t+2p^2-5p+3)}{(p-1)(2p-3)(2p-1)},
\end{align*}
as long as $p\ne \frac12,1,\frac32$.

The second family of solutions (well defined for $p$
different from $0,\frac12,1$) has $\mu = -8$. The vector $P(t)$ can be taken
as
\begin{align*}
P_0(t)&=  \frac {3t(pt-t-p)}{(p-1)p(2p-1)},\displaybreak[0]\\
P_1(t)&=\frac
{3(p^2t^2-2pt^2+t^2-2p^2t+2pt+t+p^2)}{(p-1)p(2p-1)},\displaybreak[0] \\
P_2(t)&=\frac {3(pt-t-p)}{(p-1)p(2p-1)},
\end{align*}
and the vector $Q(t)$ can be taken as
\begin{align*}
Q_0(t)&=\frac {3t(pt-p+1)}{(p-1)p(2p-1)},\displaybreak[0] \\
Q_1(t)&=-\frac
{3(p^2t^2-2p^2t+2pt+t+p^2-2p+1)}{(p-1)p(2p-1)}, \displaybreak[0] \\
Q_2(t)&=\frac {3(pt-p+1)}{(p-1)p(2p-1)}.
\end{align*}

By changing $p$ into $1-p$ in the first family we get the third family of
spherical functions, well defined for $p$ not in the set
$-\frac12,0,\frac12$.
We have  $\mu = 8p-8$, corresponding to $k=2$. The vector $P(t)$ can be
taken as
\begin{align*}
P_0(t)&= \frac {3t^2}{p(2p-1)(2p+1)}, \displaybreak[0]\\
P_1(t)&=\displaystyle{\frac {3t(2pt-t-2p-1)}{p(2p-1)(2p+1)},}
\displaybreak[0]\\
P_2(t)&= \frac {3(2p^2t^2-pt^2-4p^2t+t+2p^2+p)}{p(2p-1)(2p+1)},
\end{align*}
and the vector $Q(t)$ can be taken as
\begin{align*}
Q_0(t)&=-\frac {3(2p^2t^2+pt^2+4p^2t+t+2p^2-p)}{p(2p-1)(2p+1)},
\displaybreak[0]\\
Q_1(t)&= \frac {3(2pt+t-2p+1)}{p(2p-1)(2p+1)},\displaybreak[0]\\
Q_2(t)&=-\frac {3}{p(2p-1)(2p+1)}.
\end{align*}

Now we analyze each one of the exceptional points.

\medskip
\noindent Exceptional point $p=-\frac12$. The third family gives
\[
\begin{split}
h_0(t) &= \frac {3t^{3/2}(2\log(t) + t^2 - 4t + 3)}{2(t-1)^3}, \\
h1(t) &= -\frac {3t^{1/2}(2t\log(t)-t^2+1)}{(t-1)^3}, \\
h_2(t) &= \frac {3t^{-1/2}(2t^2\log(t) - 3t^2+4t-1)}{2(t-1)^3}.
\end{split}
\]

\medskip
\noindent Exceptional point $p=0$. The third family gives
\[
\begin{split}
h_0(t) &= -\frac {3t(2t\log(t)-t^2+1)}{(t-1)^3}, \\
h_1(t) &= \frac {6t(t\log(t) + \log(t) - 2t+2)}{(t-1)^3}, \\
h_2(t) &= -\frac {3(2t\log(t)-t^2+1)}{(t-1)^3}.
\end{split}
\]

The second family gives the same spherical function, since in both cases we
have
$\lambda = 0$, $\mu=-8$.

\medskip
\noindent Exceptional point $p=\frac12$. The third family gives
\[
\begin{split}
h_0(t) &= \frac {3t^{1/2}(2t^2\log(t) - 3t^2+4t-1)}{2(t-1)^3}, \\
h_1(t) &= -\frac {3t^{1/2}(2t\log(t)-t^2+1)}{(t-1)^3}, \\
h_2(t) &= \frac {3t^{1/2}(2\log(t)+t^2-4t+3)}{2(t-1)^3}.
\end{split}
\]
The second family gives
\[
\begin{split}
h_0(t) &= \frac {6t^{3/2}(t\log(t)+\log(t)-2t+2)}{(t-1)^3}, \\
h_1(t) &= -\frac {3t^{1/2}(t^2\log(t)+6t
\log(t)+\log(t)-4t^2+4)}{(t-1)^3}, \\
h_2(t) &= \frac {6t^{1/2}(t\log(t)+\log(t)-2t+2)}{(t-1)^3}.
\end{split}
\]

The first family gives the same spherical function as the one obtained from
the
third family, since they share the value of $(\lambda,\mu)$.

Here we find something {\em different}. The first and third family give
(unsurprisingly) the same limiting spherical function, but the second family
gives something else. This is as it should be since the values of
$(\lambda,\mu)$ for the first and third families are $\lambda =-1$, $\mu=-4$
while for the second family we get $\lambda =-1$, $\mu=-8$.

The exceptional points $p=$1 and $p=\frac32$ reproduce the results of the
points
$p=0$ and $p=-\frac12$, respectively.

In summary, when $\ell=2$ we have a total of {\em four}
spherical functions with logarithmic terms, one each from the pairs
$p=-\frac12,\frac32$
and $p=0,1$ and two from the value $p=\frac12$.

In Figure 1 we display the three families in the $(p,\mu)$-plane.
The exceptional points are marked on the $p$ axis. On each one of the
lines giving the different families we mark the points giving rise to
spherical functions with logarithmic terms.

\begin{center}
\begin{figure}
\hspace{0.25in}
\begin{picture}(200,150)(10,20)

\put(50,75){\vector(1,0){200}}
\put(100,0){\vector(0,1){150}}

\put(100,25){\circle*{2}}

\put(70,25){\line(1,0){140}}

\put(65,110){\line(1,-1){150}}
\put(65,-10){\line(1,1){150}}
\multiput(75,0)(25,25){3}{\circle*{4}}
\multiput(125,25)(0,25){2}{\circle*{4}}
\multiput(175,0)(0,0){1}{\circle*{4}}
\multiput(125,50)(25,-25){3}{\circle*{4}}

\multiput(75,75)(25,0){5}{\circle*{2}}


\put(95,80){0}
\put(67,80){$-\tfrac12$}
\put(120,80){$\tfrac12$}
\put(148,80){1}
\put(170,80){$\tfrac32$}

\put(225,25){$\mu$=-8}
\put(225,140){$\mu$=8$p$-8}
\put(225,-40){$\mu$=-8$p$}

\put(230,67){$p$}
\put(85,130){$\mu$}

\end{picture}
\vspace{.18in}
\caption{}

\medskip

\end{figure}

\end{center}

\noindent {\bf The case $\ell=5$.}
It is now clear what the general pattern is, as illustrated here for
$\ell=5$.

There are six families of spherical functions and each one of
them has six exceptional values of $p$.
The table below gives the exceptional values of $p$ for each family.

\[
\begin{array}{rrrrrrr}
1st&        -2,&  -\tfrac32,&         -1,&   -\tfrac12,&         0,&
\tfrac12.\\
2nd& -\tfrac32,&         -1,&  -\tfrac12,&           0,&  \tfrac12,&
1.\\
3rd&        -1,&  -\tfrac12,&          0,&    \tfrac12,&         1,&
\tfrac32.\\
4th& -\tfrac12,&          0,&   \tfrac12,&           1,&  \tfrac32,&
2.\\
5th&         0,&   \tfrac12,&          1,&    \tfrac32,&         2,&
\tfrac52.\\
6th&  \tfrac12,&          1,&   \tfrac32,&           2,&  \tfrac52,&
3.
\end{array}
\]

and the corresponding values of $\mu$ are listed below
$$\mu = -20p,\; -4(3p+5),\; -4(p+8),\; 4(p-9),\; 4(3p-8),\; 20(p-1).$$

The first and sixth family are related by the exchange of $p$ into $1-p$,
and the same relation connects the second and fifth families, and the
third and fourth ones. This relation is clearly manifested in the
corresponding exceptional sets as well as in the values of $\mu$ above.

For $p=-2$ (and the corresponding symmetric value $p=3$) the first
and sixth families give {\em one} spherical function with $\lambda =24$,
$\mu=-60$.

For $p=-\tfrac32$ (and its companion $p=\tfrac52$) the first and sixth
families
contribute one spherical function, with $(\lambda,\mu)=(15,-50)$.
The second and fifth families give a common spherical function
with $(\lambda,\mu)=(15,-50)$. In this way we get {\em one} more spherical
function with logarithmic terms.

For $p=-1$ (and its companion value $p=2$) the first and sixth families
as well as the third and fourth families give a common spherical function
with $(\lambda,\mu)= (8,20)$.
The second and fifth families give a common spherical function with
$(\lambda,\mu)=(8,-44)$. This gives us {\em two} more spherical functions
with logarithmic terms.

For $p=-\tfrac12$ (and its companion value $p=\tfrac32$) the first and sixth
families
give a spherical function which is common with the one that the third
family gives for $p=3/2$ and the fourth gives for $p=-\tfrac12$. The values
of
$(\lambda,\mu)$ are $(3,-30)$.
Now $p=-\tfrac12$ is also an exceptional value for the third family and
$p=\tfrac32$
is one for the fourth family. This produces another spherical function
which coincides with the one that corresponds to the second and fifth
family. The value of $(\lambda,\mu)$ is $(3,-38)$. So, for the pair
$p=-\tfrac12,\tfrac32$ there are  {\em two} spherical functions with
logarithmic terms.

For the pair $p=0,1$ there are {\em three} different spherical functions
corresponding
to the first and sixth families, the second and fifth, and the third and
fourth. The corresponding values for the eigenvalues of $D$ and $E$ are
$(0,-20)$, $(0,-32)$ and $(0,-36)$.

For $p=\tfrac12$ there are {\em three} different spherical functions
corresponding to the first and sixth families, the second and fifth
families, and the third and fourth ones. The values for the parameters
$(\lambda,\mu)$ are
$(-1,-10)$, $(-1,-26)$ and $( -1,-34)$, respectively.

In summary, for $\ell=5$, we have a total of {\em twelve} spherical
functions with logarithmic terms.

\medskip

We close the section by stating the conjecture advertised above. From
Theorem \ref{punch} we get a completely explicit form for the spherical
functions in terms of Gauss' hypergeometric functions and a matrix $U$
whose columns are given by Hahn polynomials (Proposition \ref{Hahn}). This
gives our functions $H(t)$ as a linear combination of the classical
functions.

A reader familiar with the results in \cite{GPT} may wonder if it is
possible to give an expression for $H(t)$
in terms of a single {\em generalized} hypergeometric function.
Such a result is offered below in the form of a conjecture.

\begin{conj}  For a given $\ell \ge 0$, the functions $H(t)$ have
components that are expressed in terms of a special class of generalized
hypergeometric
functions of the form $\lw{p+2}F_{p+1}$, namely
\[
\lw{p+2}F_{p+1} \left( \begin{array}{c}
a,b,s_1+1,\dots,s_p+1 \\
c,s_1,s_2,\dots,s_p \end{array} ;z \right).
\]
\end{conj}

\medskip

More explicitly, the $k$-th family of spherical functions, characterized by
$\lambda = 4p(p-1)$, $\mu_k = -4p(\ell -2k)-4k(\ell-k+1)$, has
components $h_i(t)$ given by
\[
h_i(t) = t^{p+k-i} {}_{a(\ell,i,k)}F_{a(\ell,i,k)-1} \left( \begin{array}{c}
\ 1,2p+k-i,\{s_j+1\} \\
\ell+2,\{s_j\} \end{array} ; 1-t\right)
\]
with $a(\ell,i,k) =  \min\{i+k+2,\ell+2\}$. Here $\{s_j\}$ denotes
a row vector
of $a(\ell,i,k)-2$ denominator parameters and $\{s_j+1\}$ are the
same coefficients
shifted up by one.

\

\section{Unitary spherical functions}\label{unit}

Given a function $\Phi:G\longrightarrow \End(V)$, where $V$ is a finite
dimensional complex
vector space with an inner product, we define
$\check\Phi:G\longrightarrow
\End(V)$
by $\check\Phi(g)=\Phi(g^{-1})^*$, where * denotes the operation of
taking
adjoint.

\begin{prop} The function $\Phi$ is spherical of type $\delta\in\hat K$
if
and only if
$\check\Phi$ is spherical of type $\delta$.
\end{prop}
\begin{proof} Let us assume that $\Phi$ is spherical of type $\delta$.
Then
\begin{equation*}
\begin{split}
\check\Phi(x)\check\Phi(y)&=\left(\Phi(y^{-1})\Phi(x^{-1})\right)^*=
\left(\int_K\chi_\delta(k^{-1})\Phi(y^{-1}kx^{-1})dk\right)^*\\
&=\int_K\chi_\delta(k)\check\Phi(xk^{-1}y)dk.
\end{split}
\end{equation*}
Moreover $\check\Phi(e)=\Phi(e)^*=I$, hence $\check\Phi$ is a spherical
function
of type $\delta$.
On the other hand since $\check{\check\Phi}=\Phi$ the proposition
follows.
\end{proof}

\begin{defn} A spherical function $\Phi:G\longrightarrow \End(V)$ is
said
to be {\em unitarizable} if there exists an inner product on $V$ such that
$\Phi(g)^*=\Phi(g^{-1})$, for all $g\in G$. In such a case we also say
that
$\Phi$ is a unitary spherical function. In other words $\Phi$ is
unitary
if and only if $\check\Phi=\Phi$.
\end{defn}

\begin{prop} If $\Phi:G\longrightarrow \End(V)$ is a unitarizable
irreducible
spherical function then there exists a unique
inner product on $V$, up to a positive multiplicative constant, which
turns
$\Phi$ unitary.
\end{prop}
\begin{proof} Let $(\;,\;)$ and $\langle \; ,\;\rangle$ be two inner
products on $V$
which make $\Phi$ unitary. Then there exists a linear
operator $A$ of $V$
such that $(u,v)=\langle Au,v\rangle$, for all $u,v\in V$. By
hypothesis we
have
$$\langle A\Phi(g)u,v\rangle=(\Phi(g)u,v)=(u,\Phi(g^{-1})v)=
\langle Au,\Phi(g^{-1})v\rangle=\langle\Phi(g)Au,v\rangle.$$
Thus $A\Phi(g)=\Phi(g)A$ for all $g\in G$. If $\lambda$ is an
eigenvalue of
$A$
then the corresponding eigenspace $V_\lambda$ is a non zero subspace of
$V$
which is stable by the set of linear transformations $\Phi(G)$. Hence,
by
the
irreducibility assumption, it follows that $V_\lambda=V$, which proves
the
proposition.
\end{proof}

\smallskip

Before tackling the next proposition it is best to recall the definitions
of $P(\delta)$, $E(\delta)$ and $\Phi(\delta)$ given following Proposition
\ref{fibra}.

\begin{prop}\label{unitaria} If $U$ is a unitary representation of $G$
on a
Hilbert space $E$
and $E(\delta)$ is finite dimensional and not zero, then the spherical
function
$\Ph_\delta$ is unitary.
\end{prop}
\begin{proof} Let $u,v\in E(\delta)$ and $g\in G$. Since $P(\delta)$ is
self
adjoint we have
\begin{equation*}
\begin{split}
(\Ph_\delta(g)u,v)&=(P(\delta)U(g)u,v)=(U(g)u,v)=(u,U(g^{-1})v)\\
                   &=(u,P(\delta)U(g^{-1})v)=(u,\Ph_\delta(g^{-1})v).
\end{split}
\end{equation*}
\end{proof}

\begin{cor} If $G$ is a compact group then any irreducible spherical
function
on $G$ is unitarizable.
\end{cor}
\begin{proof} If $\Phi$ is an irreducible spherical function on $G$ of
type
$\delta$
then there exists an irreducible finite dimensional representation
$(E,U)$ of $G$
such that $\Phi=P(\delta)U$ on $E(\delta)$. Since $U$ is unitarizable
from Proposition \ref{unitaria} it follows that $\Phi$ is unitarizable.
\end{proof}

This is the reason why the issue of unitarizability was not raised in
\cite{GPT}.

\medskip

Now we go on to determine which are the unitarizable irreducible
spherical functions among all of those associated to
the pair (SL(2,$\CC$),SU(2)).

\begin{lem}\label{Omega} Let $V$ be a finite dimensional complex vector
space with an inner product. If
$\Phi:G\longrightarrow \End(V)$ is an irreducible spherical function on
$G=\SL(2,\CC)$ and
\begin{equation*}
(\Omega+\overline\Omega)\Phi=a\Phi,\quad\quad
(\Omega-\overline\Omega)\Phi=b\Phi,
\end{equation*}
then
\begin{equation*}
(\Omega+\overline\Omega)\check\Phi=\overline a\check\Phi,\quad\quad
(\Omega-\overline\Omega)\check\Phi
=-\overline b\check\Phi.
\end{equation*}
\end{lem}
\begin{proof} For $X,Y\in\lieg$ and $F:G\longrightarrow \End(V)$ a
smooth
function we have
$(XY\,\check F)(e)=(YX\,F)^*(e)$. In fact
\begin{equation*}
\begin{split}
(XY\,\check F)(e)&=\left(\frac{\partial}{\partial
t}\frac{\partial}{\partial
s}\check F(\exp tX\exp sY)\right)_{s=t=0}\\
                  &=\left(\frac{\partial}{\partial
t}\frac{\partial}{\partial
s}F(\exp(-sY)\exp(-tX)\right)^*_{s=t=0}
=(YX\,F)^*(e).
\end{split}
\end{equation*}
{}From Proposition \ref{generadoresD(G)G} we get
\begin{equation*}
\begin{split}
2(\Omega+\overline\Omega)&=H_1^2+V_1^2+V_2^2-H_2^2-W_1^2-W_2^2,\\
\Omega-\overline\Omega&=-i(H_1H_2+V_1W_2-V_2W_1).
\end{split}
\end{equation*}
Now we can compute the eigenvalue of $\Omega+\overline\Omega$
corresponding
to $\check\Phi$,
$$((\Omega+\overline\Omega)\check\Phi)(e)=((\Omega+\overline\Omega)\Phi)^*
(e)=\overline a\Phi(e)^*=\overline a I.$$
Similarly, since $[H_1,H_2]=[V_1,W_2]=[W_1,V_2]=0$ we have
$$((\Omega-\overline\Omega)\check\Phi)(e)=-i(i(\Omega-\overline\Omega)\Phi)^
*
(e)=-((\Omega-\overline\Omega)\Phi)^*(e)=-\overline b\Phi(e)^*=-\overline b
I.$$
This completes the proof of the lemma.
\end{proof}

\medskip

No we come to the main result of this section.
\smallskip
\begin{thm}\label{unitarizabilidad} An irreducible spherical function
$\Phi$
on $G$ is
unitarizable if and only if the eigenvalues
$$a=((\Omega+\overline\Omega)\Phi)(e) \quad{\text and}\quad
b=((\Omega-\overline\Omega)\Phi)(e)$$
are, respectively, real and purely imaginary.
\end{thm}
\begin{proof} If $\Phi$ is unitarizable there exists an inner product
on $V$
such that $\Phi=\check\Phi$.
Thus from Lemma \ref{Omega} it follows that $a=\overline a$ and
$b=-\overline b$.

Conversely, let us take on $V$ an inner product such that
$\pi(k)=\Phi(k)$
becomes a unitary
representation of $K$. This implies that $(D\Phi)(e)=(D\check\Phi)(e)$
for
all $D\in D(K)^K$.
Moreover, from the same Lemma \ref{Omega} and from the hypothesis we
obtain
that
$$((\Omega+\overline\Omega)\Phi)(e)=((\Omega+\overline\Omega)\check\Phi)(e)
\quad\text{and}\quad
((\Omega-\overline\Omega)\Phi)(e)=((\Omega-\overline\Omega)\check\Phi)(e).$$
Since $D(G)^K=D(K)^K\otimes D(G)^G$ and $D(G)^G$ is generated by the
algebraically independent
elements $\Omega$ and $\overline\Omega$, from Proposition
\ref{unicidad} it
follows that
$\Phi=\check\Phi$, as we wanted to prove.
\end{proof}

\medskip

We recall that to each irreducible spherical function $\Phi$ on $G$ of type
$\pi=\pi_\ell$ we
associated a function $H$ on $G$ defined by $H(g)=\Phi(g)\Phi_\pi(g)^{-1}$
which can be also viewed
as a function on $\HH=G/K$ since it is right invariant under $K$. In
Section \ref{onevariable} we
also introduced the functions $\tilde H(r)=H(0,r)$ and $H(t)=\tilde H(\sqrt
t)$.
This last one is an
eigenfunction of the differential operators $D$ and $E$.

{}From Remark \ref{relacionautovalores2} we easily get the following
relations
between the eigenvalues $a$ and $b$ of $\Omega+\overline\Omega$ and
$\Omega-\overline\Omega$
corresponding to $\Phi$ and the eigenvalues $\lambda$ and $\mu$ of $D$ and
$E$
corresponding to $H=H(t)$:
\begin{equation}\label{relacionautovalores3}
2\lambda=a-b, \quad\text{and}\quad \mu=b-\ell(\ell+2).
\end{equation}

\begin{prop}\label{check} Let $H=(h_0(t),\dots,h_\ell(t))$ and
$K=(k_0(t),\dots,k_\ell(t))$ be the functions
associated, respectively, to the spherical  functions $\Phi$ and
$\check\Phi$ of type
$\pi=\pi_\ell$. If $DH=\lambda H$ and $EH=\mu H$ then
$$DK=\big(\overline\lambda+\overline\mu+\ell(\ell+2)\big)K\quad{\text
and}\quad
EK=\big(-\overline\mu-2\ell(\ell+2)\big)K.$$
Moreover
$$k_i(t)=t^{(\ell-2i)/2}\overline{h_{\ell-i}(t)}.$$
\end{prop}
\begin{proof} Let $\lambda_1=(DK)(e)$ and $\mu_1=(EK)(e)$. Then from Lemma
\ref{Omega} and
\eqref{relacionautovalores3} we get
\begin{equation*}
\begin{split}
2\lambda_1&=\overline a+\overline
b=2\big(\overline\lambda+\overline\mu+\ell(\ell+2)\big)\\
\mu_1&=-\overline b-\ell(\ell+2)=-\overline\mu-2\ell(\ell+2).
\end{split}
\end{equation*}
To prove the second assertion, for $s>0$ we let
$$a_s=\left(\begin{matrix}s&0 \\ 0&s{-1}\end{matrix}\right).$$
Then, using the basis $\{v_i\}$ of $V_\pi$ introduced in Section
\ref{onevariable}, we have
$\pi(a_s)v_i=s^{\ell-2i}v_i$. Also we have $p(a_s)=(0,s^{-2})\in\HH$, which
says that
$H(a_s)=H(0,s^{-2})=\tilde H(s^{-2})$. This implies that
\begin{equation}\label{Phi}
\Phi(a_s)v_i=s^{\ell-2i}\tilde h_i(s^{-2})v_i.
\end{equation}

If we consider on $V_\pi$ an inner product such that $\pi$ becomes a unitary
representation
of $K$, then the basis $\{v_i\}$ is orthogonal. Taking into account that the
adjoint of a
linear transformation defined by a diagonal matrix $T$ with respect to an
orthogonal basis
is the linear transformation defined by the conjugate matrix $\overline T$,
we
obtain
$$\check \Phi(a_s)v_i=\Phi(a_{s^{-1}})^*v_i=s^{-(\ell-2i)}\overline{\tilde
h_i(s^2)}v_i.$$
On the other hand changing $\Phi$ by $\check\Phi$ in \eqref{Phi} we get
$\check\Phi(a_s)v_i=s^{\ell-2i}\tilde k_i(s^{-2})v_i$. Therefore
$\tilde k_i(s^{-2})=s^{-2(\ell-2i)}\overline{\tilde h_i(s^2)}.$ If we put
$r=s^{-2}$ and use
Lemma \ref{Hdiagonal} we have
$$\tilde k_i(r)=r^{(\ell-2i)}\overline{\tilde
h_i(r^{-1})}=r^{(\ell-2i)}\overline{\tilde h_{\ell-i}(r)}.$$
Now replacing $r=\sqrt t$, and since $k_i(t)=\tilde k_i(\sqrt t)$ and
$h_{\ell-i}(t)=\tilde h_{\ell-i}(\sqrt t)$ we finally get
$k_i(t)=t^{(\ell-2i)/2}\overline{h_{\ell-i}(t)},$
which completes the proof of the proposition.
\end{proof}

\begin{cor}
In terms of the parameters $\lambda$ and $\mu$ the corresponding spherical
function
$\Phi$ of type $\pi=\pi_\ell$ is unitarizable if and only if
$$\lambda=\overline\lambda+\overline\mu+\ell(\ell+2)\quad{\textstyle
and}\quad
\mu=-\overline\mu-2\ell(\ell+2).$$
Moreover this happens precisely when
$$h_i(t)=t^{(\ell-2i)/2}\overline{h_{\ell-i}(t)}.$$
\end{cor}

\begin{prop}\label{check1} If $\Phi=\Phi_{(p,k)}$ is the irreducible
spherical function of
type $\ell$ associated to the parameters $(p,k)$ then
$\check\Phi=\Phi_{(p',k')}$
where $k'=k$ and $p'=(\ell-2k)/2+1-\bar p$.
\end{prop}
\begin{proof} The spherical function $\Phi_{(p,k)}$ of type $\ell$
corresponds
to the eigenvalues $\mu(p,k)=-4p(\ell-2k)-4k(\ell-k+1)$ and
$\lambda(p,k)=4p(p-1)$,
and the spherical function $\check\Phi_{(p,k)}$ corresponds
to the eigenvalues $\mu'=-\bar\mu-2\ell(\ell+2)$ and
$\lambda'=\bar\lambda+\bar\mu+\ell(\ell+2)$, see Proposition \ref{check}.
Therefore
the only thing we need to do is to check that $\lambda'=\lambda(p',k')$ and
that
$\mu'=\mu(p',k')$. This is a straightforward computation.
\end{proof}

\begin{cor}\label{unitarizables}
The spherical function $\Phi_{(p,k)}$ of type $\ell$ is unitarizable
if and only if

\noindent(i) when $\ell\neq 2k$, $\re(p)=\frac{\ell-2k}4+\frac12$,

\noindent(ii) when $\ell=2k$, $\re(p)=\frac12$ or $p\in\RR$.
\end{cor}
\begin{proof}
 From Proposition \ref{fibra} and Proposition \ref{check1} it follows that
$\Phi_{(p,k)}$ is
unitarizable if and  only if
\begin{equation}\label{pp'}
(p,k)=
\begin{cases}
(p',k),\\
(1-p',\ell-k),\\
(p',\ell+1-k-2p')\quad &\text{if}\;\; \ell+1-k-2p'\in\{0,\dots,\ell\},\\
(1-p',k-1+2p')\quad &\text{if}\;\; k-1+2p'\in\{0,\dots,\ell\}.
\end{cases}
\end{equation}
with $p'=(\ell-2k)/2+1-\bar p$.
\end{proof}
Now we just need to observe that the first and the second cases in
\eqref{pp'} are respectively
(i) and (ii),
and that the third and fourth cases occur when $\ell=2k$ and $2p=1$ which is
included in (ii).

\begin{remark}\label{real}
Recall that in terms of the parameters $v,r$ which parametrize naturally the
spherical
function of type $\ell$ associated to the principal series representation
$U^{\sigma,\nu}$ we have
$$v=4i\Big(\frac{\ell-2k}{4}+\frac12-p\Big)\quad\text{and}\quad r=2k-\ell,$$
and thus if $\ell\ne 2k$ the previous Corollary says that $\Phi_{v,r}$ is
unitarizable
exactly when $v$ is real. This corresponds, as is well known, to the cases
when
$U^{\sigma,\nu}$ is unitary (see \cite{Kn}). If $\ell=2k$ ($r=0$) we see
that
$\Phi_{v,0}$ is unitarizable  when $v$ is either real or purely imaginary;
when $w=iv$
satisfies $0<w<2$ we are precisely in the case of the complementary series
of SL$(2,\CC)$.
\end{remark}
\

\section{A matrix valued form of the bispectral property}\label{biespectral}

We close this paper by displaying a rather intriguing {\it matrix
valued} three term recursion relation that we have found
for a function put together from our spherical
functions. A similar relation was conjectured in \cite{GPT} and
proved in \cite{GPT1}.
In some sense the relation in those papers is more natural than the one
found here since there we were dealing with
a compact situation. Here, in spite of the fact that we are in a
non-compact case we obtain, once again, a three-term recursion. It may
appear
more natural to search for a differential equation in the spectral
parameter, and we start by discussing this point first.

In the classical case, with $\ell=0$, the spherical function is
given after setting $s=\log t$ by
\[ H(t,p) = H(s,p) = \frac{\sinh((2p-1)s/2)}{(2p-1)\sinh(s/2)}, \]
and the differential equation  $D H = \lambda H$ reads
\[ t^2 H_{tt} -\frac{2t^2}{1-t} H_t  = (p-1) p H,\]
or in the variable $s$,
\[ H_{ss} + \coth (s/2) H_s  = (p-1) p H.\]
In this case we get a differential equation in $p$, namely
\[ H_{pp}  + \frac4{2p-1} H  = s^2 H.\]

This result, and the fact that in this case we are dealing
with a special case of the Jacobi functions $(\alpha=1/2, \beta=1/2 )$
led one of us, see \cite{G}, to inquire if this type of ``bispectral
situation'' could hold in a more general setup. A tentative negative answer
to this question for other values of $\alpha,\beta$ is given in Section~10
of \cite{K}. A very detailed
analysis of this issue led to
\cite{DG}, from which it is clear that the question in \cite{G} has indeed
a negative answer except in very exceptional cases, such as the one that
leads to the three dimensional hyperbolic space. We now take up this
question in the context of matrix valued spherical functions.

For a given nonnegative integer  $\ell$  consider the matrix
whose rows are given by the vectors $H(t)$
corresponding to the values $k=0,1,2,\dots ,\ell$ discussed in the comments
preceding and following Theorem \ref{punch}. Denote
the corresponding matrix by
\[ \Phi (t,p). \]

In the spirit of \cite{GPT}, \cite{GPT1} one could look for a relation of
the type
\[
A(p) \Phi_{pp}(t,p) + B(p) \Phi_p(t,p) + C(p)  \Phi(t,p)  = M(t) \Phi(p,t)
\]
where $A(p), B(p), C(p)$ are matrices independent of $t$ and the matrix
$M(t)$
is independent of $p$.
Of course, we could try to place the coefficient matrices to the
right of the function $\Phi(t,p)$ and its derivatives, and this would lead
to another set of possible differential equations in the spectral parameter.

We have succeeded in finding such matrices (by allowing for the second
possibility only) but the results are rather
disappointing: the matrices in question are scalar functions multiplied
by the matrix made up of all ones.
When one looks at this result and observes the fact that the trace of each
one
of our (diagonal) matrix valued functions $H(t,p)$ equals
$(\ell+1)$ times the value of these functions for $\ell=0$, we see
that nothing new has been learned. As we see below the situation is very
different if one looks for three term recursions.

\subsection{ The bispectral property}
As a function of $t$, $\Phi(t,p)$ satisfies two differential equations
\[
D \Phi(t,p)^t = \Phi(t,p)^t\Lambda , \quad E \Phi(t,p)^t=\Phi(t,p)^tM,
\]
$D$ and $E$ are the differential operators introduced earlier. Moreover we
have

\begin{conj} There exist matrices $A_p$, $B_p$, $C_p$, independent of $t$,
such that
\[
A_p\Phi(t,p-1)+B_p\Phi(t,p)+C_p \Phi(t,p+1) = \ \frac{t^2+1}{t}\Phi(t,p).
\]
\end{conj}

The matrices $A_p$ (upper) and $C_p$ (lower)  consist, in general,
    of three
diagonals each and $B_p$ is
tridiagonal. Recall that for $t=1$ the matrix
$\Phi(1,p)$ consists of all ones.

\medskip
\noindent We illustrate these results below for $\ell=0,1,3,4$.

a) $\ell =0$. Here we have
\[
\Phi(t,p) =\frac{t^{p}-t^{1-p}}{(2p-1)(t-1)},
\]
and our property reads
\[
\frac{2p-3}{2p-1} \,\Phi(t,p-1) +
\frac{2p+1}{2p-1}\,\Phi(t,p+1) \ = \ \frac{t^2+1}{t}\, \Phi(t,p).
\]

b) $\ell =1$. Here we have $\mu =-4p,4p-4$, and for the
matrix $\Phi$ we get
\[
\Phi(t,p)=\left( \begin{array}{cc}
\frac{(2p-2)t^{2p+1}-(2p-1)t^{2p}+t^2}{(p-1)(2p-1)(t-1)^2t^p}
&
\frac{t^{2p}+(1-2p)t^2+2(p-1)t}{(p-1)(2p-1)(t-1)^2t^p}\\ {}\\
\frac{t^{1-p}(t^{2p}-2pt+2p-1)}{p(2p-1)(t-1)^2}
&
\frac{(2p-1)t^{2p+1}-2pt^{2p}+t}{p(2p-1)(t-1)^2t^p}
\end{array}\right),
\]
and our property reads
\begin{eqnarray*}
&& \left(\begin{array}{cc}
\frac{p-2}{p-1} & \frac{1}{(p-1)(2p-1)}\\ {}\\
0 & \frac{2p-3}{2p-1}\end{array}\right)\,
\Phi(t,p-1) +
\left(\begin{array}{cc}
\frac{1}{(1-p)(2p-1)} & \frac{1}{(p-1)(2p-1)}\\ {}\\
\frac{1}{p(2p-1)} & \frac{1}{p(1-2p)}\end{array}\right)\,
\Phi(t,p)   \\ {}\\{}\\ && \qquad +
\left(\begin{array}{cc}
\frac{2p+1}{2p-1} & 0\\ {}\\
\frac{1}{p(2p-1)} & \frac{p+1}{p}\end{array}\right)
\,\Phi(t,p+1) \ = \ \frac{t^2+1}{t} \, \Phi(t,p).
\end{eqnarray*}

c) $\ell =3$. Here we have
    $\mu =-12p,-4p-12, 4p-16, 12p-12$, and
with $\Phi(t,p)$ defined as above we get
\[
A_p \Phi(t,p-1)+B_p\Phi(t,p) + C_p\Phi(t,p+1)
=\frac{t^2+1}{t} \,\Phi(t,p)
\]
with

\[
A_p = \left(
\begin{array}{cccc}
\frac {p-3}{p-2} & \frac {3}{(p-1)(2p-3)} & \frac
{3}{(p-2)(p-1)(2p-3)(2p-1)} &
0 \\ \vspace{2\jot}
0 & \frac {(p-2)p(2p-5)}{(p-1)^2(2p-3)} & \frac
{16(p-2)p}{(p-1)(2p-3)(2p-1)^2} & \frac {3}{(p-1)^2(2p-1)^2} \\
\vspace{2\jot}
0 & 0 & \frac
{(p-2)(2p-3)(2p+1)}{(p-1)(2p-1)^2} & \frac {3(2p-3)}{(p-1)(2p-1)^2} \\
\vspace{2\jot}
0 & 0 & 0 &
\frac {2p-3}{2p-1}
\end{array}
\right)
\]

\[
B_p = \left(
\begin{array}{cccc}
-\frac {3}{(p-1)(2p-3)} & \frac {3}{(p-1)(2p-3)} & 0 & 0 \\
\vspace{2\jot}
\frac {3(p-2)}{(p-1)^2(2p-3)} & -\frac {14p^2 - 21p -
9}{(p-1)p(2p-3)(2p-1)} & \frac {(2p-3)(2p+1)}{(p-1)^2p(2p-1)} & 0 \\
\vspace{2\jot}
0 & \frac
{(2p-3)(2p+1)}{(p-1)p^2(2p-1)} & -\frac
{14p^2-7p-16}{(p-1)p(2p-1)(2p+1)}
& \frac {3(p+1)}{p^2(2p+1)} \\ \vspace{2\jot}
0 & 0 & \frac {3}{p(2p+1)} & -\frac {3}{p(2p+1)}
\end{array}
\right)
\]

\[
C_p = \left(
\begin{array}{cccc}
\frac {2p+1}{2p-1} & 0 & 0 & 0 \\
\vspace{2\jot}
\frac {3(2p+1)}{p(2p-1)^2} & \frac {(p+1)(2p-3)(2p+1)}{p(2p-1)^2} & 0 &
0 \\ \vspace{2\jot}
\frac {3}{p^2(2p-1)^2} & \frac
{16(p-1)(p+1)}{p(2p-1)^2(2p+1)} & \frac {(p-1)(p+1)(2p+3)}{p^2(2p+1)} & 0 \\
\vspace{2\jot}
0 & \frac {3}{p(p+1)(2p-1)(2p+1)} & \frac
{3}{p(2p+1)} & \frac {p+2}{p+1}
\end{array}
\right)
\]

d)  $\ell =4$. Here we have $\mu =-16p,-8p-16, -24, 8p-24, 16p-16$,
and with $\Phi(t,p)$ defined as above we have
\[
A_p \Phi(t,p-1)+B_p\Phi(t,p) + C_p\Phi(t,p+1)
=\frac{t^2+1}{t} \,\Phi(t,p)
\]
where $A_p$ is of the form
\[
A_p = \left(
\begin{array}{ccccc}
a_{00}& a_{01}&a_{02}&0&0\\
0& a_{11}& a_{12}&a_{13}&0 \\
0&0& a_{22}& a_{23}&a_{24}\\
0&0&0& a_{33}& a_{34}\\
0&0&0&0& a_{44}
\end{array} \right)
\]
with
\[
\begin{array}{lll}
a_{00}=\frac {2p-7}{2p-5}, & a_{01}=\frac {2}{(p-2)(p-1)}, &  a_{02}=\frac
{6(p-2)^{-1}(p-1)^{-1}}{(2p-5)(2p-1)},  \\ &&\\
a_{11}=\frac {(p-3)p(2p-5)}{(p-2)(p-1)(2p-3)},  & a_{12}=\frac
{3p(2p-5)}{(p-2)(p-1)^2(2p-1)}, & a_{13}= \frac
{9}{(p-1)^{2}(2p-3)(2p-1)},\\ && \\
a_{22}=\frac
{p(p-2)(2p-5)(2p+1)}{(p-1)^2(2p-1)^2},  & a_{23}=\frac
{3(p-2)(2p+1)}{(p-1)^2p(2p-1)}, & a_{24}=\frac {6}{(p-1)p(2p-1)^2}, \\& & \\
a_{33}=\frac
{(p-2)(p+1)(2p-3)}{(p-1)p(2p-1)}, & a_{34}=\frac {2(2p-3)}{(p-1)p(2p-1)},
  & a_{44}=\frac {2p-3}{2p-1}.
\end{array}
\]

\

\noindent The matrix $B_p$ is of the form
\[
B_p = \left(
\begin{array}{ccccc}
b_{00}& b_{01}&0&0&0\\
b_{10}& b_{11}& b_{12}&0&0 \\
0&b_{21}& b_{22}& b_{23}&0 \\
0&0&b_{32}& b_{33}& b_{34}\\
0&0&0&b_{43}& b_{44}
\end{array} \right)
\]
with
\[
\begin{array}{lll}
b_{00}=-\frac {2}{(p-2)(p-1)}, & b_{01}=\frac {2}{(p-2)(p-1)}, &\\&&\\
b_{10}=\frac {2(2p-5)}{(p-2)(p-1)(2p-3)}, & b_{11}=-\frac
{5p^2-10p-4}{(p-2)(p-1)^2p}, & b_{12}=\frac{3(p-2)(2p+1)}{(p-1)^2p(2p-3)},
\\&&\\
b_{21}=\frac {3(p-2)(2p+1)}{(p-1)^2p(2p-1)}, & b_{22}= -\frac
{3(2p^2-2p-3)}{(p-1)^2p^2}, & b_{23}=\frac {3(p+1)(2p-3)}{(p-1)p^2(2p-1)},\\
&&\\
b_{32}=\frac
{3(p+1)(2p-3)}{(p-1)p^2(2p+1)}, & b_{33}= -\frac
{5p^2-9}{(p-1)p^2(p+1)}, & b_{34}= \frac {2(2p+3)}{p(p+1)(2p+1)},\\  && \\
  b_{43}=\frac {2}{p(p+1)}, & b_{44}=-\frac {2}{p(p+1)}. &
\end{array}
\]

\

\noindent The matrix $C_p$ is of the form
\[
C_p = \left(
\begin{array}{ccccc}
c_{00}&0&0&0&0\\
c_{10}& c_{11}&0&0&0 \\
c_{20}&c_{21}& c_{22}&0&0 \\
0&c_{31}&c_{32}& c_{33}&0\\
0&0&c_{42}&c_{43}& c_{44}
\end{array} \right)
\]
with
\[
\begin{array}{lll}
c_{00}=\frac {2p+1}{2p-1}, \\ && \\
c_{10}=\frac {2(2p+1)}{(p-1)p(2p-1)}, &c_{11}=\frac
{(p-2)(p+1)(2p+1)}{(p-1)p(2p-1)}, \\ && \\
c_{20}=\frac{6}{(p-1)p(2p-1)^2} ,&
c_{21}=\frac{3(p+1)(2p-3)}{(p-1)p^2(2p-1)}, & c_{22}=\frac
{(p-1)(p+1)(2p-3)(2p+3)}{p^2(2p-1)^2}, \\ && \\
c_{31}=\frac {9}{p^2(2p-1)(2p+1)}, & c_{32}=\frac
{3(p-1)(2p+3)}{p^2(p+1)(2p-1)}, & c_{33}= \frac
{(p-1)(p+2)(2p+3)}{p(p+1)(2p+1)}, \\ && \\
c_{42}=\frac{6(p+1)^{-1}}{p(2p-1)(2p+3)}, & c_{43}=\frac {2}{p(p+1)}, &
c_{44}=\frac {2p+5}{2p+3}.
\end{array}
\]

\

\

\section{Appendix}

For completeness we include here the proof of those  propositions and lemmas
stated in the paper in Sections \ref{redG/K} and
\ref{onevariable}. We start with the first order differential operator $E$.

\begin{propE} For any $H\in C^\infty(\HH)\otimes
\End(V_\pi)$
we have
\begin{align*}
EH&= H_x \;\dot\pi\matc{z-\overline
z}{\;\;\frac{2(1+|z|^2)}r}{2\,r}{\overline z-z} +
H_y \;\dot\pi\matc{i(z+\overline
z)}{\;\;\frac{2i(1+|z|^2)}r}{-2i\,r}{i(z+\overline z)}\\
&\quad +H_r \;\dot\pi\matc{-2r}{\;\;4z}{0}{\;\;2r}.
\end{align*}
\end{propE}

\begin{proof} By equation \eqref{Edefuniv} we have
\begin{align}\label{menem}
EH=& H_1(H)\, \Phi_\pi \dot\pi(H_1)\Phi_\pi^{-1}+V_1(H)\, \Phi_\pi
\dot\pi(V_1)\Phi_\pi^{-1}+V_2(H)\, \Phi_\pi
\dot\pi(V_2)\Phi_\pi^{-1}\\
\nonumber =& H_1(H)\, \Phi_\pi \dot\pi(H_1)\Phi_\pi^{-1}+ (V_1(H)+iV_2(H))\,
\Phi_\pi
\dot\pi(E_{12})\Phi_\pi^{-1}\\ \label{expre}
\nonumber &+ (V_1(H)-iV_2(H))\, \Phi_\pi \dot\pi(E_{21})\Phi_\pi^{-1}.
\end{align}

We start by computing $H_1(H)$. For $g=\matc{a}{b}{c}{d}\in G$ we have
\begin{align*}
p(g\exp tH_1)&=\textstyle p \left(g\matc{e^t}{0}{0}{e^{-t}}\right)= g
\matc{e^{2t}}{0}{0}{e^{-2t}}g^* \left
(\begin{smallmatrix} 0\\1\end{smallmatrix}\right)\\
&=\left (\begin{matrix} a\overline c e^{2t}+b\overline d e^{-2t}\\
|c|^2e^{2t}+|d|^2e^{-2t}\end{matrix}\right)\\
&=\left (\begin{matrix}u(t) \\ v(t)\end{matrix}\right)=\left
(\begin{matrix}u_1(t)+i\,u_2(t) \\ v(t)\end{matrix}\right).
\end{align*}

Also we have
\begin{align*}
  \left(\frac{d \,u}{dt}\right)_{t=0}
&=\left(\frac{d \,u_1}{dt}\right)_{t=0}
+ i\,\left(\frac{d \,u_2}{dt}\right)_{t=0} = 2(a\overline c-b\overline
d),
\\  \left(\frac{d \,v}{dt}\right)_{t=0}&=2(|c|^2-|d|^2).
\end{align*}
By the chain rule we get
$$H_1(H)(g)=2 H_x \re(a\overline c-b\overline d)+2 H_y \im(a\overline
c-b\overline d)+ 2 H_r(|c|^2-|d|^2).$$

\noindent For $V_1$ we have
\begin{align*}
p(g\exp tV_1)&=g \exp{tV_1}(\exp {tV_1})^* g^*\left
(\begin{smallmatrix}
0\\1\end{smallmatrix}\right)\displaybreak[0]\\
& =\matc{a}{b}{c}{d} \matc{\cosh t}{\sinh t}{\sinh t}{\cosh
t}\matc{\cosh
t}{\sinh t}{\sinh t}{\cosh t} \matc{\overline
a}{\overline c}{\overline b}{\overline d}\left
(\begin{matrix}
0\\1\end{matrix}\right) \displaybreak[0]\\&= \left(
\begin{matrix}
(\cosh^2 t+\sinh^2 t)(a\overline c+b\overline d) + 2(b\overline
c+a\overline
d)\sinh t \cosh t \displaybreak[0]\\
(\cosh^2 t+\sinh^2 t) (|c|^2+|d|^2)+2(d\overline c+c\overline d)\sinh t
\cosh t
\end{matrix}\right )\displaybreak[0]\\
&=\left (\begin{matrix}\tilde u(t) \\ \tilde
v(t)\end{matrix}\right)=\left
(\begin{matrix}\tilde u_1(t)+i\,\tilde u_2(t) \displaybreak[0]\\
\tilde v(t)\end{matrix}\right).
\end{align*}
Thus,
\begin{align*}
  \left(\frac{d \,\tilde u}{dt}\right)_{t=0}
  = 2(b\overline c+a\overline d),\qquad
\left(\frac{d \,\tilde v}{dt}\right)_{t=0}&=2(d\overline c+c\overline d).
\end{align*}
We get
$$V_1(H)(g)=2 H_x \re(b\overline c+a\overline d)+2 H_y \im(b\overline
c+a\overline d)+ 2 H_r (d\overline c+c\overline d).$$

\noindent Similarly we get
$$V_2(H)(g)=2 H_x \im(-b\overline c+a\overline d)-2H_y
\re(b\overline
c-a\overline d)+ 2i H_r (c\overline d-d\overline c).$$
Therefore
\begin{align*}
\left(V_1(H)+iV_2(H)\right)(g)&=2H_x (b\overline c+ \overline a
d)-2iH_y
(b\overline c- \overline a d)  +4H_r d\overline c, \displaybreak[0] \\
\left(V_1(H)-iV_2(H)\right)(g)&=2H_x (\overline b c+ a \overline d)+
2iH_y
(\overline b c- a \overline d)  +4H_r c\overline d.
\end{align*}

\smallskip

\noindent On the other hand we compute
$$\Phi_\pi(g)\dot\pi(X)\Phi_\pi^{-1}(g)= \left( \frac{d}{dt}\right)
_{t=0}
\pi(g\exp{tX}g^{-1})= \dot\pi(gXg^{-1}).$$
Then
\begin{align*}
\Phi_\pi(g)\dot\pi(E_{12})\Phi_\pi^{-1}(g)&=\dot\pi\matc{-ac}{\;a^2}{-c^2}{\
;ac}
,\displaybreak[0]\\
\Phi_\pi(g)\dot\pi(E_{21})\Phi_\pi^{-1}(g)&=\dot\pi\matc{bd}{\;-b^2}{d^2}{\;
-bd}
,\displaybreak[0]\\
\Phi_\pi(g)\dot\pi(H_{1})\Phi_\pi^{-1}(g)&=\dot\pi\matc{ad+bc}{-2ab}{2cd}{-a
d-bc}.
\end{align*}
Therefore, by \eqref{menem} we have
$EH(g)=H_x\, A +H_y \, B+ H_r\, C$, where
\begin{align*}
A&=2(b\overline c+ \overline a d)\,\dot\pi\matc{-ac}{\;a^2}{-c^2}{\;ac}
+ 2(\overline b c+a\overline
d)\,\dot\pi\matc{bd}{\;-b^2}{d^2}{\;-bd} \\
& \qquad +2\re(a\overline c-b\overline
d)\,\dot\pi\matc{ad+bc}{-2ab}{2cd}{-ad-bc}\displaybreak[0],\\
B&= -2i(b\overline c- \overline a d)\,\dot\pi\matc{-ac}{\;a^2}{-c^2}{\;ac}
+ 2i(\overline b c-a\overline d)\,\dot\pi\matc{bd}{\;-b^2}{d^2}{\;-bd} \\ &
\qquad +2\im(a\overline c-b\overline
d)\,\dot\pi\matc{ad+bc}{-2ab}{2cd}{-ad-bc}\displaybreak[0],\\
C&=4 d\overline c\,\dot\pi\matc{-ac}{\;a^2}{-c^2}{\;ac}
+4\overline d c\,\dot\pi\matc{bd}{\;-b^2}{d^2}{\;-bd} \\
& \qquad +
2(|c|^2-|d|^2)\,\dot\pi\matc{ad+bc}{-2ab}{2cd}{-ad-bc}.
\end{align*}
By means of straightforward calculations we obtain
\begin{align*}
A&=\dot\pi\matc{z-\overline z}{\;\;\frac{2(1+|z|^2)}r}{2\,r}{\overline
z-z}
\, ,\qquad
B=\dot\pi\matc{-i(z+\overline
z)}{\;\;\frac{2i(1+|z|^2)}r}{-2i\,r}{i(z+\overline
z)}\displaybreak[0]\\
C&=\dot\pi\matc{-2r}{\;\;4z}{0}{\;\;2r}.
\end{align*}
This completes the proof of the proposition.
\end{proof}

\

\begin{propD} For any $H\in C^\infty(\HH)\otimes
\End(V_\pi)$
we have

\begin{align*}
\textstyle\frac 14 DH&= (\re^2 z +1)\,H_{xx}+(\im^2 z
+1)\,H_{yy}+r^2\,H_{rr}+2\re z\im z \, H_{xy} \\
& \quad
+2r\re z \,H_{xr} +2r\im z \, H_{yr}+3\re z\, H_x+3\im z \,H_y+3r
\,H_r.
\end{align*}
\end{propD}

\begin{proof} We have $DH=\left(H_1^2+V_1^2+V_2^2\right)H$. For any
$X\in
\lieg $ one has
$$X^2(H)(g)= \left(\frac{d}{ds}\,\frac{d}{dt} H\left(p(g\exp
(s+t)X)\right)
\right)_{s=t=0}.$$
Let us introduce the functions $u,v,u_1, u_2$ by means of
\begin{equation*}
p(g\exp (s+t)X)=(u(s,t)\, , \, v(s,t))=(u_1(s,t)\, , \, u_2(s,t)\, ,
\,v(s,t)),
\end{equation*}
where $u(s,t)=u_1(s,t)+iu_2(s,t)$ and $u_1, u_2, v$ are  real valued
functions. We note that if $u(t)=u_1(t)+iu_2(t)$, then
$\frac{du_1}{dt}=\re(\frac{du}{dt})$ and
$\frac{du_2}{dt}=\im(\frac{du}{dt})$.
Using the chain rule we have
\begin{align*}
X^2(H)(g)&= H_{xx}\, \frac{\partial u_1}{\partial s}\,\frac{\partial
u_1}{\partial t} + H_{yy}\, \frac{\partial u_2}{\partial
s}\,\frac{\partial u_2}{\partial t} + H_{rr}\, \frac{\partial
v}{\partial
s}\,\frac{\partial v}{\partial t}\displaybreak[0]\\&
+ H_{xy}\,
\left(\frac{\partial u_1}{\partial s}\,\frac{\partial u_2}{\partial
t}+\frac{\partial u_2}{\partial s}\,
\frac{\partial u_1}{\partial t}\right) + H_{xr}\,
\left(\frac{\partial v}{\partial s}\,\frac{\partial u_1}{\partial
t}+\frac{\partial u_1}{\partial s}\,
\frac{\partial v}{\partial t}\right) \displaybreak[0]\\&+ H_{yr}\,
\left(\frac{\partial v}{\partial s}\,\frac{\partial u_2}{\partial
t}+\frac{\partial u_2}{\partial s}\,
\frac{\partial v}{\partial t}\right) + H_{x}\frac{\partial^2
u_1}{\partial
t\partial s}
+H_{y}\frac{\partial^2 u_2}{\partial t\partial s}
+H_{r}\frac{\partial^2 v}{\partial t\partial s}.
\end{align*}

We need to handle separately the cases $X=H_1$, $X=V_1$ and $X=V_2$.
For $X=H_1$ we start with
\begin{align*}
p(g\exp (s+t)H_1)&=
\matc{a}{b}{c}{d}\matc{e^{2(t+s)}}{0}{0}{e^{-2(t+s)}}\matc{\overline
a}{\overline c}{\overline
b}{\overline d}
\left (\begin{matrix} 0\\1\end{matrix}\right)\displaybreak[0]\\
&=\left (\begin{matrix} a\overline c e^{2(t+s)}+b\overline d
e^{-2(t+s)}\\
|c|^2e^{2(t+s)}+|d|^2e^{-2(t+s)}\end{matrix}\right).
\end{align*}

\noindent Therefore at $s=t=0$ we have

$$ \frac{\partial u}{\partial s}
=\frac{\partial u}{\partial t}= 2a\overline c-2b\overline d\, ,\qquad
\frac{\partial v}{\partial s}=\frac{\partial v}{\partial t}
=2(|c|^2-|d|^2),$$
and
$$ \frac{\partial^2 u}{\partial s\partial t}
=4(a\overline c+b\overline d)\, ,\qquad
\frac{\partial^2 v}{\partial s\partial t}=4(|c|^2+|d|^2).$$

\noindent Thus we obtain
\begin{equation}\label{H1^2H}
\begin{split}
H_1^2(& H) = 4H_{xx}\re^2(a\overline c-b\overline d) +
4H_{yy}\im^2(a\overline c-b\overline d)
+4H_{rr}(|c|^2-|d|^2)^2\displaybreak[0]\\
&  +8H_{xy}\re (a\overline c-b\overline d)\im(a\overline c-b\overline
d)
+8H_{xr}\re(a\overline c-b\overline d)(|c|^2-|d|^2)\displaybreak[0]\\&
+ 8 H_{yr}\im(a\overline c-b\overline d)(|c|^2-|d|^2)+4H_x\re
(a\overline
c+b\overline d)\displaybreak[0]\\&
+4H_y\im (a\overline c+b\overline d) +4H_r (|c|^2+|d|^2).
\end{split}
\end{equation}

\smallskip

  For $X=V_1$ we have $\exp tV_1=\matc{\cosh t}{\sinh t}
{-\sinh t}{\cosh t}$ and
$ p(g\exp (s+t)V_1)=(u(s,t)\, , \, v(s,t))$
where
\begin{align*}
u(s,t)=& \left( \cosh^2(t+s)+\sinh^2(t+s)\right) \left(a\overline
c+b\overline d \right)\\
& \quad+2\sinh (s+t)\cosh(s+t)\left(b\overline c+a\overline d \right),
\\
v(s,t) =&\left(\cosh^2(t+s)+\sinh^2(t+s)\right)\left(|c|^2+|d|^2\right)
\\&
\quad +2\sinh (s+t)\cosh (s+t) \left(d\overline c+c\overline d\right).
\end{align*}

\noindent Therefore at $s=t=0$ we have

$$ \frac{\partial u}{\partial s}
=\frac{\partial u}{\partial t}= 2(b\overline c+a\overline d)
\, ,\qquad
\frac{\partial v}{\partial s}=\frac{\partial v}{\partial t}
=2(d \overline c+c\overline d),$$
$$ \frac{\partial^2 u}{\partial s\partial t}
=4(a\overline c+b\overline d)\, ,\qquad
\frac{\partial^2 v}{\partial s\partial t}=4(|c|^2+|d|^2),$$

\noindent and we obtain
\begin{equation}\label{V1^2H}
\begin{split}
V_1^2(H) &= 4H_{xx}\re^2(b\overline c+a\overline d) +
4H_{yy}\im^2(b\overline c+a\overline d) +4H_{rr}(d \overline
c+c\overline
d)^2\displaybreak[0]\\
& +8H_{xy}\re (b\overline c+a\overline d)\im(b\overline c+a\overline d)
+8H_{xr}(d \overline c+c\overline d)\re(b\overline c+a\overline d)\\& +
8H_{yr}(d \overline c+c\overline d)\im(b\overline
c+a\overline d)+4H_x\re (a\overline c+b\overline d)\displaybreak[0]\\
&+4H_y\im (a\overline c+b\overline d) +4H_r (|c|^2+|d|^2).
\end{split}
\end{equation}

\smallskip

\noindent Finally for $X=V_2$ we obtain
\begin{equation}\label{V2^2H}
\begin{split}
V_2^2(H)& = 4H_{xx}\im^2(a\overline d-b\overline c) +
4H_{yy}\re^2(a\overline d-b\overline c) -4H_{rr}(c\overline d-d\overline
c)^2\displaybreak[0]\\ & -8H_{xy}\re (a\overline d-b\overline
c)\im(a\overline d-b\overline c)
-8iH_{xr}(c\overline d-d\overline
c)\im(a\overline d-b\overline c)\\
& + 8 i H_{yr}(c\overline d-d\overline
c)\re(a\overline d-b\overline
c)+4H_x\re (a\overline c+b\overline d)\displaybreak[0]\\
&+4H_y\im (a\overline c+b\overline d) +4H_r (|c|^2+|d|^2).
\end{split}
\end{equation}

\noindent Now, the proposition follows from equations
\eqref{H1^2H}, \eqref{V1^2H} and \eqref{V2^2H}  and by making use of
the
expressions
$$ r= |c|^2+|d|^2 \quad \text{and}\quad z=a\bar c +b\bar d.$$
\end{proof}

\

Now we shall compute all first and second order partial derivatives that we
used in
Section \ref{onevariable}. Given $(z,r)\in \HH$, $z=z_1+iz_2$,
$z_1,z_2\in\RR$, $r\in \RR_{>0}$ we
define
\begin{equation}\label{sdef}
s=\frac 1{2r}\left(1+|z|^2+r^2 - \sqrt{(|z|^2+(r+1)^2)(|z|^2+
(1-r)^2)}\right),
\end{equation}
\begin{equation}\label{abdef}
a=\frac{1-rs}{\sqrt{(1-rs)^2+(s\,|z|)^2}}
\qquad \text{and}\qquad
b=\frac{s\,z}{\sqrt{(1-rs)^2+(s\,|z|)^2}}.
\end{equation}
Then the matrix
$A=A(z,r)= \matc{\overline a}{\;-b}{\bar b}{\;a}$ satisfies
$(z,r)= A\cdot(0,s)$. In fact in Section \ref{onevariable} we defined a
function $k=k(z,r)$
such that $k(z,r)\cdot(z,r)=(0,s)$, and $A=k^{-1}$. Therefore
  $$H(z,r)=\pi(A(z,r))\tilde H(s(z,r))\pi(A(z,r)^{-1}).$$

\smallskip
\begin{lem} \label{der2H}
At $(0,r) \in \HH$ we have
$$ \frac{\partial H}{\partial z_j}(0,r)=
\left[ \dot\pi\Bigl( \frac{\partial A}{\partial z_j}\Bigr)\, ,\, H(r)
\right],$$
and
\begin{align*}
\frac{\partial ^2 H}{\partial z_j^2}(0,r)=&
- \frac{2r}{1-r^2} \frac{d \tilde H}{d r}
+\frac{\partial ^2 (\pi\circ A)}{\partial z_j^2} \tilde H(r)
+ \tilde H(r) \frac{\partial ^2 (\pi\circ A^{-1})}{\partial z_j^2}  \\
& - 2\dot\pi\Bigl(\frac{\partial A}{\partial z_j}\Bigr)\, \tilde H(r) \,
\dot\pi\Bigl(\frac{\partial A}{\partial z_j}\Bigr).
\end{align*}
\end{lem}

\begin{proof} We have
$$H(z,r)=\pi(A(z,r))\tilde H(s(z,r))\pi(A(z,r)^{-1}).$$
Then for $j=1,2$
\begin{equation}\label{der1}
\begin{split}
\displaystyle\frac{\partial H}{\partial z_j} =&
  \frac{\partial (\pi\circ A)}{\partial z_j}  \,( \tilde H \circ s )\,
(\pi\circ A^{-1})
+ (\pi \circ A)\, \frac{\partial (\tilde H\circ s)}{\partial z_j}\,
(\pi\circ A^{-1}) \displaybreak[0] \\
&+ (\pi\circ A )\, (\tilde H\circ s )\, \frac{\partial (\pi\circ
A^{-1})}{\partial z_j}.
\end{split}
\end{equation}
We also have
$$ \frac{\partial (\tilde H\circ s)}{\partial z_j}=
\Bigl(\frac{d\tilde H}{dr}\circ s\Bigr) \frac{\partial s}{\partial
z_j}.$$

\noindent Note that $s(0,r)=r$, $A(r,0)=I$ and that at
$(r,0)$
$$\frac{\partial s}{\partial z_j}=0 ,\qquad
  \frac{\partial (\pi\circ A)}{\partial z_j} =
\dot\pi \Bigl(\frac{\partial  A}{\partial z_j} \Bigr)\, , \qquad
\frac{\partial (\pi\circ A^{-1})}{\partial z_j} =
-\dot\pi \Bigl(\frac{\partial  A}{\partial z_j} \Bigr).$$
Therefore $ \frac{\partial H}{\partial z_j}(0,r)=
\left[ \dot\pi\Bigl( \frac{\partial A}{\partial z_j}\Bigr)\, ,\, H(r)
\right]$.

To compute the second order derivatives we start by observing that at
$(r,0)$
we have
$$ \frac{\partial ^2 (\tilde H \circ s)}{\partial z_j^2}
=\frac{d^2 \tilde H}{dr^2}\left(\frac{\partial s}{\partial z_j}\right)^2
+\frac{d\tilde H}{dr}\frac{\partial ^2 s}{\partial z_j^2}
=-\frac{2r}{1-r^2} \frac{d\tilde H}{dr}.
$$

\noindent Now the lemma follows by differentiating the expression given in
\eqref{der1} and evaluating it at $(r,0)\in \CC^2$.
\end{proof}

\

\begin{lem}\label{lemgral}
At $(0,r) \in \HH$ we have
\begin{align*}
\frac{\partial ^2 (\pi\circ A)}{\partial z_j^2}=&
\dot\pi\Bigl ( \frac{\partial ^2  A}{\partial z_j^2}\Bigr )
- \dot\pi\Bigl ( \Bigl( \frac{\partial A}{\partial z_j}\Bigr)^2\Bigr )
+ \dot\pi\Bigl ( \frac{\partial A}{\partial z_j}\Bigr )^2
\end{align*}
and
\begin{align*}
\frac{\partial ^2 (\pi\circ A^{-1})}{\partial z_j^2}=&
-\dot\pi\Bigl ( \frac{\partial ^2  A}{\partial z_j^2}\Bigr )
+ \dot\pi\Bigl ( \Bigl( \frac{\partial A}{\partial z_j}\Bigr)^2\Bigr )
+ \dot\pi\Bigl ( \frac{\partial A}{\partial z_j}\Bigr )^2
\end{align*}
\end{lem}

For the proof see Proposition 13.3 in \cite{GPT}. Here $\pi$ refers to its
extension to GL$(2,\CC)$
which is trivial on the center.

\

\noindent Finally we give the proof of the following proposition, stated in
Section \ref{onevariable}.

\begin{propderxy}
Let $J=\begin{pmatrix}0&1\\-1&0\end{pmatrix}$ and $T=\begin{pmatrix}
0&1\\ 1&0\end{pmatrix}$. We have
\begin{align*}
H_{x}(0,0,r)&= -{\frac r{1-r^2}}\left( \dot\pi(J)\tilde H(r)-\tilde
H(r)\dot\pi(J)\right),\\
H_{y}(0,0,r)&= -{\frac {i\,r}{1-r^2}}\left( \dot\pi(T)\tilde H(r)-\tilde
H(r)\dot\pi(T)\right),
\end{align*}
and
\begin{equation*}
\begin{split}
H_{xx}(0,0,r)&= -\frac {2r}{1-r^2} \frac{d \tilde H}{d r} +
\frac{r^2}{(1-r^2)^2}\left(\dot\pi(J)^2\, \tilde H(r)+ \tilde
H(r)\,
\dot\pi(J)^2\right)\displaybreak[0]\\
&\quad -\frac{2r^2}{(1-r^2)^2}\dot\pi(J)\tilde H(r)\, \dot\pi(J), \\
H_{yy}(0,0,r)&= -\frac {2r}{1-r^2} \frac{d \tilde H}{d r}
- \frac{r^2}{(1-r^2)^2}\left(\dot\pi(J)^2\, \tilde H(r)+ \tilde H(r)\,
\dot\pi(J)^2\right)\displaybreak[0]\\
&\quad +\frac{2 r^2}{(1-r^2)^2}\dot\pi(J)\tilde H(r)\, \dot\pi(J).
\end{split}
\end{equation*}
\end{propderxy}

\begin{proof}
The matrix $A$ in the previous lemmas is defined by $A=A(z,r)=
\matc{\overline a}{\;-b}{\bar b}{\;a}$ with $a,b$ given in \eqref{abdef}.

We start by computing $ \frac{\partial A}{\partial z_j}$ and
$\frac{\partial^2 A}{\partial z_j^2}$. It is easy to verify that at
$(0,r)\in \HH$ we have
$$\frac{\partial a}{\partial z_j}=0, \quad
\frac{\partial^2 a}{\partial z_j^2}=-\frac{r^2}{(1-r^2)^2},\quad
\frac{\partial b}{\partial z_j}=\delta_j\frac{r}{1-r^2}, \quad
\frac{\partial^2 b}{\partial z_j^2}=0,$$
where
$\delta_j=1\;$ if $\;j=1$, and $\delta_j=i\;$ if $\;j=2$. Then at $(0,r)\in
\HH$
$$ \frac{\partial A}{\partial z_1}=-\frac{r}{1-r^2} J, \qquad
  \frac{\partial A}{\partial z_2}=-\frac{i r}{1-r^2} T , \qquad
  \frac{\partial^2 A}{\partial z_j^2}=-\frac{r^2}{(1-r^2)^2} I. $$
We observe that
$\left( \frac{\partial A}{\partial z_j}\right)^2=  \frac{\partial^2
A}{\partial z_j^2}$, then by Lemma \ref{lemgral} we have at $(0,r)$
$$\frac{\partial ^2 (\pi\circ A)}{\partial z_j^2}= \frac{\partial ^2
(\pi\circ A^{-1})}{\partial z_j^2}=
\dot \pi\left(\frac{\partial A}{\partial z_j}\right)^2 .$$
Now the proposition follows directly from Lemma \ref{der2H}.

\end{proof}

\

\

\newcommand\bibit[5]{\bibitem[#1]
{#2}#3, {\em #4,\!\! } #5}

\end{document}